\documentclass[12pt]{article}
\usepackage{amsmath}
\usepackage{amsfonts}
\usepackage{amssymb}
\usepackage{booktabs} 
\usepackage{graphicx}
\usepackage[hidelinks, pdfborder={0 0 0}]{hyperref}
\usepackage{natbib}
\usepackage[normalem]{ulem} 
\graphicspath{{./Fig/}}
	
\usepackage{color}

\usepackage{amsthm}
\usepackage{latexsym}
\usepackage{lscape}
\usepackage{natbib}

\newtheorem{theorem}{Theorem}[section]
\newtheorem{lemma}[theorem]{Lemma}
\newtheorem{definition}[theorem]{Definition}

\newtheorem{remark}[theorem]{Remark}

\newtheorem{assumption}[theorem]{Assumption}

\newcommand{\ignore}[1]{}

\newcommand{\Matrix}[1]{\ensuremath{\left[\begin{array}{rrrrrrrrrrrrrrrrrr} #1 \end{array}\right]}}

\newcommand{\Matrixc}[1]{\ensuremath{\left[\begin{array}{cccccccccccc} #1 \end{array}\right]}}
\newcommand{\Eqn}[1]{\ensuremath{\begin{array}{cccccccccc}#1\end{array}}}

\newcommand{\R}{{\mathbf R}}

\newcommand{\z}{{\mathbf z}}

\renewcommand{\epsilon}{\varepsilon}
{\begingroup

 \begin{enumerate}}
 {\end{enumerate}\endgroup}

\newcommand{\HH}{{\mathcal H}}
\newcommand{\II}{{\mathcal I}}

\newcommand{\LL}{{\mathcal L}}
\newcommand{\RR}{{\mathcal R}}

 \newcommand{\RED}[1]{{\color{black}{#1}}}

\usepackage{epsfig}
\usepackage{epstopdf}
\usepackage[usenames,dvipsnames]{xcolor}    

\newcommand{\mbF}{\mathbf{F}}

\newcommand{\mbn}{\mathbf{n}} 
\newcommand{\mbx}{\mathbf{x}}
\newcommand{\mbu}{\mathbf{u}}
\newcommand{\mbw}{\mathbf{w}}

\newcommand{\mby}{\mathbf{y}}

\newcommand{\norm}[1]{\left\lVert#1\right\rVert}

\title{Shape \emph{versus} timing: linear responses of a limit cycle with hard boundaries under instantaneous and static perturbation}

\author{Yangyang~Wang\thanks{Department of Mathematics, The University of Iowa, Iowa City, IA 52242, USA. Email: yangyang-wang@uiowa.edu}, Jeffrey P.~Gill\thanks{Department of Biology, Case Western Reserve University, Cleveland, OH 44106, USA. Email: jpg18@case.edu}, Hillel J.~Chiel\thanks{Departments of Biology, Neurosciences and Biomedical
		Engineering, Case Western Reserve University, Cleveland, OH 44106, USA. Email: hjc@case.edu}, Peter J.~Thomas\thanks{Departments of Biology, Mathematics, Applied Mathematics, and Statistics, Case Western Reserve University, Cleveland, OH 44106, USA. Email: pjthomas@case.edu}}
\date{} 

\begin{document}
	\maketitle
	\abstract{When dynamical systems that produce rhythmic behaviors operate within hard limits, they may exhibit limit cycles with sliding components, that is, closed isolated periodic orbits that make and break contact with a constraint surface.  Examples include heel-ground interaction in locomotion, firing rate rectification in neural networks, \RED{and stick-slip oscillators}. In many rhythmic systems, robustness against external perturbations involves response of both the shape and the timing of the limit cycle trajectory.  
		The existing methods of infinitesimal phase response curve (iPRC) and variational analysis are well established for quantifying changes in timing and shape, respectively, for smooth systems.  These tools have recently been extended to nonsmooth dynamics with transversal crossing boundaries. In this work, we further extend \RED{the iPRC method} to nonsmooth systems with sliding components, \RED{which enables us to make predictions about the synchronization properties of weakly coupled stick-slip oscillators}. We observe a new feature of the isochrons in a planar limit cycle with hard sliding boundaries: a nonsmooth kink in the asymptotic phase function, originating from the point at which the limit cycle smoothly departs the constraint surface, and propagating away from the hard boundary into the interior of the domain. Moreover, the classical variational analysis neglects timing information and is restricted to instantaneous perturbations. By defining the ``infinitesimal shape response curve" (iSRC), we incorporate timing sensitivity of an oscillator to describe the shape response of this oscillator to parametric perturbations. In order to extract timing information, we also develop a ``local timing response curve" (lTRC) that measures the timing sensitivity of a limit cycle within any given region. We demonstrate in a specific example that taking into account local timing sensitivity in a nonsmooth system greatly improves the accuracy of the iSRC over global timing analysis given by the iPRC.}
	
	
	\section{Introduction}\label{sec:intro}
	
	A \emph{limit cycle with sliding component} (LCSC) is a closed, isolated, periodic orbit of an $n$-dimensional, autonomous, deterministic nonsmooth dynamical system, in which the trajectory is constrained to move along a surface of dimension $k<n$ during some portion of the orbit. \RED{The motion of a trajectory sliding along a constraint surface is called a \textit{sliding mode}.  LCSCs appear naturally in dynamical systems models of physiological and robotic motor control systems \citep{barajon1992,gelfand1988,mortin1989,holmes2006,revzen2011,lyttle2017,guckenheimer2018} as well as mechanical stick-slip systems \citep{GB99,Galvanetto2001,DL11,LN2013}. In control theory, the sliding mode concept has been used to design controllers for nonlinear systems \citep{SS83,Slotine84,Lee2009}.}
	
	\RED{Both natural and engineered motor systems are robust to certain short and long term disturbances.  Studies of the robustness properties of these systems have relied on applying the \textit{infinitesimal phase response curve} (iPRC) and variational analysis to quantify changes in timing and shape of the motor trajectory to weak perturbations. For understanding the response to instantaneous perturbations, the two methods are ubiquitously used in the literature of smooth dynamical systems \citep{spardy2011a,spardy2011b,Park2017}. 
		\citep{filippov1988,bernardo2008,LN2013,DL11} show that variational analysis can be applied to nonsmooth systems including LCSC systems. Recently,} the iPRC has also been generalized to nonsmooth systems, provided the flow is always transverse to any switching surfaces at which nonsmooth transitions occur \citep{shirasaka2017,park2018,CGL18,wilson2019}. \RED{However, there have been fewer reported studies that analyze the model response (especially the shape response) to perturbations that are sustained over long times. Even fewer works have analyzed the response of LCSCs, in which the transverse flow condition fails, to both instantaneous and sustained perturbations.  Our goal in this paper is to bridge such a knowledge gap by providing a first description of the \textit{infinitesimal shape response curve} (iSRC) that can account for the shape response of an oscillator to sustained (e.g., parametric) perturbations and extending both iPRC and iSRC to LCSCs in nonsmooth systems. }

\begin{figure}[!t]
		\begin{center}
			\includegraphics[width=6in]{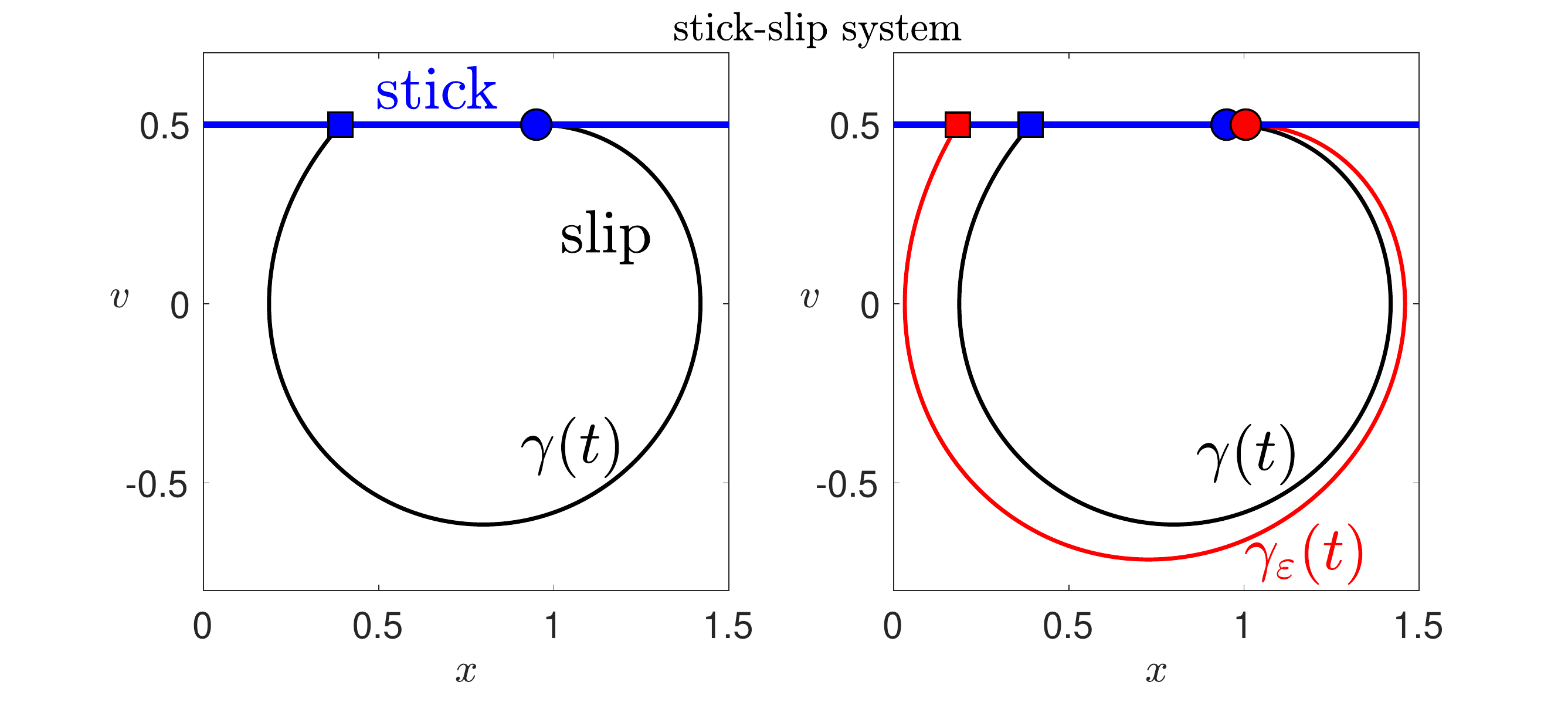}
		\end{center}
		\caption{\label{fig:ss-motivating} Periodic LCSC solution trajectories of a stick-slip system. Left: projection of the unperturbed trajectory $\gamma(t)$ onto the $(x,v)$ phase space (black), where $x$ is the displacement of the mass from the position in which the spring assumes its natural length and $v$ is the velocity of the mass. The trajectory enters and leaves the hard boundary (blue line) at $\mbx_{\rm land}$ (blue square) and $\mbx_{\rm liftoff}$ (blue circle), respectively. Right: Projections of the unperturbed trajectory $\gamma(t)$(black) and the trajectory under parametric perturbation $\gamma_{\varepsilon}(t)$(red). See \S \ref{sec:one-mass} for further details.}
	\end{figure}

	\RED{In this paper we consider the case of continuous LCSC solutions to nonsmooth systems with degree of smoothness one or higher; that is, systems with continuous trajectories, also known as \emph{Filippov} systems \citep{filippov1988,bernardo2008}. The simplest model of such a system can be written as follows 
		\begin{eqnarray} \label{eq:sliding_model}
		\frac{d\textbf{x}}{dt}=F(\mbx):=\left\{
		\Eqn{F^{\rm slide}(\mbx), &\quad\quad \mbx \in\RR^{\rm slide}\subset \Sigma \\
			F^{\rm interior}(\mbx), & \text{otherwise}\\}
		\right.
		\end{eqnarray}
		where $\mbx$ denotes the state variable.
		The trajectory is confined to travel within the closure of the domain $\RR$, the boundary of which is defined to be the hard boundary $\Sigma$. A trajectory entering the sliding region $\RR^{\rm slide}\subset\Sigma$ will slide along it with the vector field $F^{\rm slide}$ until it is allowed to reenter the interior. For points not in $\RR^{\rm slide}$, the dynamics is determined by $F^{\rm interior}$. 
		
		A stick-slip oscillator is one example of a system that exhibits a LCSC \citep{GB99}. A mass on a belt that moves at a constant velocity $u$ is connected to a fixed support by a linear elastic spring and by a linear dashpot. 
		Mechanical systems of this type are referred to as \textit{stick-slip} since there are times when the mass and the belt are moving together (stick phase) and others in which the mass slips relative to the belt (slip phase). 
		Such a solution trajectory alternating between stick and slip phases is a LCSC, as illustrated in Figure \ref{fig:ss-motivating}, left. Here the hard boundary is given by $v=u$, where $v$ is the velocity of the mass and $u$ is the driving velocity of the belt. Thus, for a stick-slip system, the ``sliding component" of the limit cycle corresponds to the ``stick" phase, during which the mass moves with fewer degrees of freedom than during the ``slip" phase.}
	
	\RED{While the response of an oscillator to an instantaneous perturbation is well understood, relatively little consideration has been given to studying the model response to sustained (e.g., parametric) perturbations even in the context of smooth systems.  
		Here,} we develop the mathematical framework required to analyze the changes in shape and timing wrought by parametric changes.
	\RED{In general, a small fixed change in a parameter gives rise to a new limit cycle, with different shape and timing than the original. See Figure \ref{fig:ss-motivating}, right. In order to account for the \textit{shape response} of an oscillator to parametric perturbations, we apply Lighthill's method \citep{JSS07} to derive an iSRC in \S\ref{sec:smooth-sust}. Later in \S\ref{sec:one-mass}, we use the iSRC to study the shape response of the stick-slip oscillator to a small parametric perturbation.} 
	In contrast to standard variational analysis, which neglects timing changes, the iSRC takes into account both timing and shape changes arising due to a parametric perturbation. 
	
	In many applications of LCSCs, the impact of a perturbation on local timing can be as important as the global effects.	For instance, any motor control or mechanical system that operates by making and breaking physical contact (e.g., walking, scratching, grasping \RED{or stick-slip}) would only experience perturbations limited to a discrete component of the limit cycle \RED{(e.g., the friction of the ground acts as a perturbation during the stance phase of locomotion and is absent during the swing phase)}. In these cases, one would need to compute the local timing changes of the trajectory during the phase when the perturbation exists to understand the robustness of this system. \RED{It is well known that the global change in timing of an oscillator due to a parametric perturbation can be captured using the iPRC, which, however, cannot be used to capture a local timing change induced by sustained perturbations in many LCSC systems.} In this paper, we develop a \emph{local timing response curve (lTRC)} that is analogous to the iPRC but measures the local timing sensitivity of a limit cycle within any given local region \RED{(e.g., the stance phase in locomotion)}.  Development of the lTRC allows us to compute local timing changes in an oscillator due to nonuniform sustained perturbations. \RED{Moreover, we show that the lTRC can be used to greatly improve the accuracy of the iSRC of a nonsmooth system. }

	\RED{Recently the iPRC has been extended to certain nonsmooth systems, but the theory does not extend directly to LCSCs in which the transverse flow condition fails. In this paper, we bridge such a knowledge gap by extending the iPRC to LCSC systems (see Theorem \ref{thm:main}). In contrast to the variational dynamics that exhibits discontinuities when a sliding motion begins \citep{filippov1988}, we find that the iPRC in a LCSC experiences a discontinuous jump when the trajectory leaves a sliding region and is continuous when the sliding mode begins.  }
	To our knowledge, ours is the first work considering iPRC for systems with sliding components.
	\footnote{\RED{For example, although the well-known monograph \textit{Hybrid Dynamical Systems} \citep{GST2009} addresses periodic orbits and synchronization in nonsmooth systems with \emph{transverse} boundary crossings, it avoids systems with \emph{sliding} components, which have non-transverse flow out of the constraint surface, nor does it discuss phase response curve methods.}}
	We illustrate the theory using a planar model consisting of four sliding components \RED{and a stick-slip system with one sliding component. 
		We also use the iPRC to predict the synchronization properties of two weakly coupled stick-slip oscillators. } 
	
	\RED{The rest of the paper is organized as follows. We consider smooth systems in \S\ref{sec:smooth-theory} and Filippov systems in \S\ref{sec:nonsmooth-theory}. In each of the two sections, we first review classical theory and methods to provide context, and then present our new results. The variational and phase response curve analysis for the responses of smooth dynamical systems to instantaneous and parametric perturbations are reviewed in \S\ref{sec:smooth-inst}. To account for shape responses to sustained perturbations, we define the iSRC in \S\ref{sec:smooth-sust} and define the lTRC in \S\ref{sec:smooth-local}. We present the classical two-zone Filippov system with transversal crossing boundary and define a Filippov system that produces a LCSC solution (see \eqref{eq:sliding_model}) in \S\ref{sec:filippov}.} While the applicability of the classical perturbative methods from \S\ref{sec:smooth-theory} is generally limited by the constraint that the dynamics of the system is smooth, some elements of the methods have already been generalized to nonsmooth systems, \RED{which are reviewed in \S\ref{sec:nonsmooth-review}. We extend the iPRC to the LCSC case in nonsmooth systems in \S\ref{sec:sliding}. The main result is summarized in Theorem \ref{thm:main}.} Appendix~\ref{ap:proof} gives a proof of the theorem. Numerical algorithms for implementing all the methods are presented in Appendix~\ref{sec:algorithm}. In \S\ref{sec:toy-model}, we illustrate both the theory and algorithms using a planar model, comprising a limit cycle with a linear vector field in the interior of a simply connected convex domain with four hard boundaries. In this example, we show that under certain circumstances (e.g.~non-uniform perturbation), the iSRC together with the lTRC provides a more accurate representation of the combined timing and shape responses to static perturbations than using the global iPRC alone. Surprisingly, we discover nondifferentiable ``kinks" in the isochron function that propagate backwards in time along an osculating trajectory that encounters the hard boundary exactly at the liftoff point (the point where the limit cycle trajectory smoothly departs the boundary). \RED{In \S\ref{sec:stick-slip}, we use the iSRC to understand the shape response of an actual mechanical system - a stick-slip oscillator - to a parametric perturbation and use the iPRC to study the synchronization of two weakly coupled stick-slip oscillators.  Lastly, we discuss limitations of our methods and possible future directions in \S\ref{sec:discussion}.} Appendix~\ref{ap:symbols} provides a table of symbols used in the paper.

	\section{Linear responses of smooth systems}\label{sec:smooth-theory}
	
	\RED{In this section we consider smooth dynamical systems.  We begin by reviewing the classical variational theory for limit cycles, and then derive new methods including the \textit{infinitesimal shape response curve} (iSRC) and the \textit{local timing response curve} (lTRC) for linear approximation of the effects of small perturbations on the timing and shape of a limit cycle trajectory in the smooth case. 
		
		Specifically, in \S\ref{sec:smooth-inst} we review the classical variational and infinitesimal phase response curve analyses that capture the shape response to small instantaneous perturbations and the effects of both instantaneous and sustained perturbations on the timing of trajectories near a limit cycle (LC) trajectory. However, little consideration has been given to the effect of sustained perturbations on the \emph{shape} of the orbit. 
		In \S\ref{sec:smooth-sust}, we derive the iSRC to account for the combined shape and timing response of a LC trajectory under \emph{sustained} (e.g., parametric) perturbations. 
		To obtain a more accurate iSRC when the limit cycle experiences different timing sensitivities in different regions within the domain, in \S\ref{sec:smooth-local}, we introduce the lTRC. In contrast to the iPRC, which measures the global shift in the period, the lTRC lets us compute the timing change of a LC trajectory within regions bounded between specified Poincar\'e sections.

		Consider a one-parameter family of $n$-dimensional dynamical systems 
		\begin{equation}\label{eq:dxdt=Feps}
		\frac{d\mbx}{dt}=F_\varepsilon(\mbx),
		\end{equation}
		indexed by a parameter $\varepsilon$ representing a static perturbation of  a reference system 
		\begin{equation}\label{eq:dxdt=F}
		\frac{d\textbf{x}}{dt}=F_0(\textbf{x}).
		\end{equation}
		
		\begin{assumption}\label{ass:smooth} Throughout this section, we make the following assumptions:
			\begin{itemize}
				\item The vector field $F_{\varepsilon}(\mbx): \Omega \times \II \to \mathbb{R}^n$ is $C^1$ in both the coordinates $\mbx$ in some open subset $\Omega\subset \mathbb{R}^n$ and the perturbation $\varepsilon\in \II\subset \mathbb{R}$, where $\II$ is a small open neighborhood of zero. 
				\item For $\varepsilon\in \II$, system \eqref{eq:dxdt=Feps} has a linearly asymptotically stable limit cycle $\gamma_\varepsilon(t)$, with a finite period $T_\varepsilon$ depending  (at least $C^1$) on $\varepsilon$.
			\end{itemize} 
		\end{assumption}
		
		It follows from Assumption \ref{ass:smooth} that when $\varepsilon=0$, $F_0(\mbx)$ is $C^1$ in $\mbx\in \Omega$ and the unperturbed system \eqref{eq:dxdt=F} exhibits a $T_0$-periodic asymptotically stable limit cycle solution $\gamma_0(t) = \gamma_0(t+T_0)$ with $0<T_0<\infty$.
		To simplify notation, we will drop the subscript $0$ and use $F(\mbx)$, $\gamma(t)$ and $T$ to denote the unperturbed vector field, limit cycle solution and period, except where required to avoid confusion.
		
		Moreover, Assumption \ref{ass:smooth} implies that we have the following approximations that will be needed for deriving the iSRC and lTRC: 
		\begin{equation} 
		F_{\varepsilon}(\mbx) = F_0(\mbx) + \varepsilon \frac{\partial F_\varepsilon}{\partial \varepsilon}(\mbx)\Big|_{\varepsilon=0} + O(\varepsilon^2),
		\end{equation}
		\begin{equation}\label{eq:T_eps_expansion} T_\varepsilon=T_0+\epsilon T_1+O(\epsilon^2),
		\end{equation}
		\begin{equation}\label{eq:x_epsilon_of_tau} \gamma_\varepsilon(\tau_\epsilon(t)) = \gamma_0(t)+\varepsilon\gamma_1(t)+O(\varepsilon^2)\quad (\text{uniformly in }t), 
		\end{equation}
		where $T_1$ is the linear shift in the limit cycle period $T_0$ in response to the static perturbation of size $\varepsilon$. 
		This global timing sensitivity $T_1>0$ if increasing $\varepsilon$ increases the period. 
		The perturbed time $\tau_\epsilon(t)$, which satisfies 
		$\tau_0(t)\equiv t$ and $\tau_\epsilon(t+T_0)-\tau_\epsilon(t)=T_\epsilon,$  will be described in 
		detail later (see \eqref{eq:tau_epsilon_conditions}); it allows the approximation \eqref{eq:x_epsilon_of_tau} to be uniform in time and permits us to compare perturbed and unperturbed trajectories at corresponding time points.  The vector function $\gamma_1(t)$ is a representative  belonging to an equivalence class that comprises the iSRC.}

	\subsection{Review of shape and timing response to perturbations (classical theory)}\label{sec:smooth-inst}  
	
	\subsubsection{Shape and timing response to \textit{instantaneous} perturbations}
	
	Suppose a small, brief perturbation is applied to \eqref{eq:dxdt=F} at time $t_0$ such that there is a small abrupt perturbation in the state space.
	We have
	\begin{equation}\label{eq:initial-inst-perturbation}
	\tilde{\gamma}(t_0)= \gamma(t_0)+\varepsilon P,
	\end{equation}
	where $\tilde{\gamma}$ indicates the trajectory subsequent to  the instantaneous perturbation, $\varepsilon$ is the magnitude of the perturbation, and $P$ is the unit vector in the direction of the perturbation in the state space. As is well known, the effects of the small brief perturbation $\varepsilon P$ on the shape and timing of the limit cycle trajectory are given, respectively, by the solution of the variational equation \eqref{eq:var}, and the iPRC which solves the adjoint equation \eqref{eq:prc}.
	
	The evolution of a trajectory $\tilde{\gamma}(t)$ close to the limit cycle $\gamma(t)$ may be approximated as $\tilde{\gamma}(t)=\gamma(t)+\mbu(t)+O(\epsilon^2)$, where $\mbu(t)$ satisfies the \emph{variational equation}
	\begin{equation}\label{eq:var}
	\frac{d\mbu}{dt}=DF(\gamma(t))\mbu
	\end{equation}
	with initial displacement $\mbu(t_0)=\epsilon P$ given by \eqref{eq:initial-inst-perturbation}, for small $\epsilon$.
	Here $DF(\gamma(t))$ is the Jacobian matrix evaluated along $\gamma(t)$.
	
	On the other hand, an iPRC of an oscillator measures the timing sensitivity of the limit cycle to infinitesimally small perturbations at every point along its cycle.
	It is defined as the shift in the oscillator phase $\theta\in [0, T)$ per size of the perturbation, in the limit of small perturbation size. The limit cycle solution takes each phase to a unique point on the limit cycle, $\textbf{x}=\gamma(\theta)$, and its inverse maps each point on the cycle to a unique phase, $\theta=\phi(\textbf{x})$. 
	One may extend the domain of $\phi(\textbf{x})$ to points in the basin of attraction $\mathcal{B}$ of the limit cycle by defining the \textit{asymptotic phase}: $\phi(\textbf{x}): \mathcal{B}\rightarrow [0,T)$
	with 
	\[
	\frac{d\phi(\mathbf{x}(t))}{dt}=1,\quad \phi(\mathbf{x}(t))=\phi(\mathbf{x}(t+T)).
	\]
	If $\textbf{x}_0\in \gamma(t)$ and $\textbf{y}_0 \in \mathcal{B}$, then we say that $\mby_0$ has the same asymptotic phase as $\mbx_0$ if $\norm{\mbx(t; \mbx_0)-\mby(t;\mby_0) }\to 0$, as $t\to \infty$.
	This means that $\phi(\mbx_0)=\phi(\mby_0)$.
	The set of all points off the limit cycle that have the same asymptotic phase as the point $\mbx_0$ on the limit cycle is the \textit{isochron} with phase $\phi(\mbx_0)$.
	\RED{The asymptotic phase function $\phi$ is defined up to an additive constant; this constant is of no consequence other than to define an arbitrary reference point as the ``zero phase" location on the limit cycle trajectory.}
	
	Suppose $\varepsilon P$ applied at phase $\theta$ results in a new state $\gamma(\theta)+\varepsilon P \in\mathcal{B}$, which corresponds to a new phase $\tilde{\theta}=\phi(\gamma(\theta)+\varepsilon P)$. The phase difference $\tilde{\theta}-\theta$ defines the phase response curve (PRC) of the oscillator. One defines the iPRC as the vector function $\z:[0,T)\to\R^n$ satisfying
	\begin{equation} \label{eq:definition_of_iPRC}
	\z(\theta)\cdot{P}=\lim_{\varepsilon\to 0}\frac{1}{\varepsilon}\left(\phi(\gamma(\theta)+\varepsilon{P})-\theta\right) = \nabla_{\mbx}\phi(\gamma(\theta))\cdot{P}
	\end{equation}
	for arbitrary unit perturbation ${P}$.
	The first equality serves as a definition, while the second follows from routine arguments \citep{brown2004,ET2010,SL2012,Park2017}.
	It follows directly that the vector iPRC is the gradient of the asymptotic phase and it captures the phase (or timing) response to perturbations in any direction $P$ in state space.
	\RED{By assumption \ref{ass:smooth}, the vector field $F$ is $C^1$.} It follows that the iPRC is a continuous $T$-periodic solution satisfying the \textit{adjoint equation} \citep{SL2012},
	\begin{equation}\label{eq:prc}
	\frac{d\z}{dt}=-DF(\gamma(t))^\intercal \z, 
	\end{equation}
	with the normalization condition
	\begin{equation}\label{eq:prc-normalization}
	F(\gamma(\theta))\cdot \z(\theta) = 1.
	\end{equation}
	
	\begin{remark}\label{rem:uTz=const}
		By direct calculation, one can show that the solutions to the variational equation and the adjoint equation satisfy $\mbu^\intercal\z=\text{constant}$:
		\begin{equation}\label{eq:var-prc}
		\frac{d(\mbu^\intercal \z)}{dt}=\frac{d\mbu^\intercal}{dt} \z +\mbu^\intercal\frac{d\z}{dt}=\mbu^\intercal DF^\intercal\z+\mbu^\intercal(-DF^\intercal\z)=0.
		\end{equation}
		This relation holds for both smooth and nonsmooth systems with transverse crossings \citep{park2018}.
	\end{remark}
	
	For completeness, we note that differences between phase variables, as in \eqref{eq:definition_of_iPRC}, will be interpreted as the 
	\emph{periodic difference}, $d_T(\phi(\mbx),\phi(\mby))$.
	That is, if two angular variables $\theta$ and $\psi$ are defined on the circle $\mathbb{S}\equiv [0,T)$, then we set
	\begin{equation}
	d_T(\theta,\psi)=\begin{cases}
	\theta-\psi+T,&\theta-\psi<-\frac{T}2\\
	\theta-\psi,&-\frac{T}2\le \theta-\psi\le \frac{T}2\\
	\theta-\psi-T,&\theta-\psi>\frac{T}2,
	\end{cases}
	\end{equation}
	which maps $d_T(\theta,\psi)$ to the range $[-T/2,T/2]$.
	In what follows we will simply write $\theta-\psi$ for clarity rather than $d_T(\theta,\psi)$.

	\subsubsection{Timing response to \textit{sustained} perturbations}\label{sec:smooth-time-sust}
	
	\RED{Next we review how the iPRC can be used to estimate the linear shift in the limit cycle period in response to a sustained perturbation (see \eqref{eq:T1}). }
	
	Note that the perturbed periodic solution $\gamma_{\varepsilon}(t)$ to system \eqref{eq:dxdt=Feps} can be represented, to leading order, by the single variable system 
	\begin{equation}\label{eq:phase-reduction}
	\frac{d\theta}{dt}=1+\z(\theta)^\intercal G(\mbx,t),
	\end{equation}
	where $G(\mbx,t)=\varepsilon\frac{\partial F_\varepsilon(\gamma(t))}{\partial \varepsilon}|_{\varepsilon=0}$ 
	represents the $O(\varepsilon)$ perturbation of the vector field, $\theta\in [0, T_0)$ is the asymptotic phase as defined above, and $\z:\theta\in[0,T_0)\to\R^n$ is the iPRC.
	Recall that for $0\le \epsilon \ll 1$ we can represent $T_\varepsilon$ with $T_\varepsilon=T_0+\varepsilon T_1+O(\varepsilon^2)$ \RED{(see \eqref{eq:T_eps_expansion})}. 
	From \eqref{eq:phase-reduction},  $T_1$ can be calculated using the iPRC as
	\begin{equation}\label{eq:T1}
	T_1=-\int_{0}^{T_0} \z(\theta)^\intercal\frac{\partial F_\varepsilon(\gamma(\theta))}{\partial \varepsilon}\Big|_{\varepsilon=0}d\theta.
	\end{equation} 		
	
	\subsection{Shape response to sustained perturbations: iSRC}\label{sec:smooth-sust} 
	
	\RED{In this section, we develop new tools to analyze the effects of sustained perturbations on the shape of a limit cycle solution.
		As the parameter $\epsilon$ in \eqref{eq:dxdt=Feps} varies, the family of limit cycles produced by the flow forms a $2$-dimensional  ``ribbon'' in the $(n+1)$-dimensional space parametrized by $(\mbx,\epsilon)$. The smoothness of this ribbon, when the vector field $F_\epsilon$ is smooth, follows immediately from the persistence of hyperbolic  over- and under-flowing invariant manifolds (\cite{Wiggins1994book}, \S 6.2).}
	
	In contrast to the instantaneous perturbation considered in the previous section, changes in each aspect of shape and timing can now influence the other, 
	and hence a variational analysis of the combined shape and timing response of limit cycles under constant perturbation is needed. 
	\RED{To this end, we develop a new method that we call the iSRC (see~\eqref{eq:src}). 
		In our analysis, we adapt Lighthill's method of coordinate perturbation (``strained coordinates") to simultaneously stretch the time coordinate so as to accommodate the effect of parameter changes on period.}

	\RED{Generically, introducing a change in a parameter will lead to a change in period, as well as a displacement of the set of points comprising the limit cycle's orbit.  In order to quantify the change in shape, we must first accommodate any change in period.  To this end we introduce a rescaled time coordinate $t\to\tau_\epsilon(t)$,} which satisfies the consistency and smoothness conditions 
	\begin{equation}
	\label{eq:tau_epsilon_conditions}
	\frac{d\tau_\epsilon}{dt}>0,\quad\text{ and }\quad\frac1\epsilon\left(\int_{t=0}^{T_0}\left(\frac{d\tau_\epsilon}{dt}\right)dt - T_0\right)=T_1+O(\epsilon),\text{ as }\epsilon\to 0,
	\end{equation}
	where $T_0$ is the unperturbed period and $T_1$ is the linear shift in the period as given by \eqref{eq:T1}.
	These conditions do not determine the value of the derivative of $\tau_\epsilon$, which we write as
	\begin{equation}\label{eq:tau_eps_derivative}
	d\tau_\epsilon/dt=1/\nu_\varepsilon(t).
	\end{equation}
	In general, the iSRC will depend on the choice of $\nu_\varepsilon(t)$.
	However, some natural choices are particularly well adapted to specific problems, as we will see.
	Initially, we will make the simple ansatz 
	\begin{equation}\label{eq:ansatz}
	\nu_\varepsilon(t)=\text{const};
	\end{equation}  that is, we will impose uniform local timing sensitivity.
	Later in \S\ref{sec:smooth-local} we will introduce the \textit{local timing response curve} (lTRC) to exploit alternative time rescalings for greater accuracy.
	
	As discussed above, to understand how the static perturbation changes the shape of the limit cycle $\gamma(t)$, we need to rescale the time coordinate of the perturbed solution so that $\gamma(t)$ and $\gamma_{\varepsilon}(t)$ may be compared  at corresponding time points.
	That is, for $\epsilon>0$ we wish to introduce a rescaled perturbed time coordinate $\tau_\varepsilon(t)$ so that \eqref{eq:x_epsilon_of_tau} holds, \RED{uniformly in time $-\infty<t<\infty$,} which we repeat below: 
	\[
	\gamma_\epsilon(\tau_\varepsilon(t))=\gamma_0(t)+\epsilon\gamma_1(t)+O(\epsilon^2).
	\]
	We define the $T_0$-periodic function $\gamma_1(t)$ to be the \emph{infinitesimal shape response curve (iSRC)}.
	
	We show next that $\gamma_1(t)$ obeys an inhomogeneous variational equation \eqref{eq:src}.
	This equation resembles \eqref{eq:var}, but has two additional non-homogeneous terms arising, respectively, from time rescaling $t\to\tau_\varepsilon(t)$, and directly from the constant perturbation acting on the vector field.
	
	It follows from \eqref{eq:tau_eps_derivative} and \eqref{eq:ansatz} that the scaling factor is $\nu_\varepsilon=\frac{T_0}{T_\varepsilon}$.
	Moreover, by \eqref{eq:T_eps_expansion}, $\nu_\varepsilon$ can be written as  
	\[\nu_\varepsilon=1-\varepsilon \nu_1+ O(\varepsilon^2),\]
	where $\nu_1=\frac{T_1}{T_0}$ represents the relative change in frequency.
	In terms of $\nu_\varepsilon$ that is time-independent, the rescaled time  for $\gamma_{\varepsilon}(\tau_\varepsilon(t))$ can be written as $\tau_\epsilon(t)=t/\nu_\varepsilon \in [0, T_\varepsilon]$ for $t\in [0, T_0]$ (see \eqref{eq:tau_eps_derivative}).
	Differentiating $\gamma_\epsilon(\tau_\epsilon(t))$ given in \eqref{eq:x_epsilon_of_tau} with respect to $t$ ($\frac{d\gamma_\varepsilon}{dt} = \frac{d\gamma_\varepsilon}{d\tau_\epsilon} \frac{d\tau_\epsilon}{dt}$), substituting the ansatz ($\frac{d\tau_\epsilon}{dt}=\frac{1}{\nu_\varepsilon}$) and rearranging lead to
	\begin{equation}\label{eq:dxdtau1}
	\begin{split}
	\Eqn{\frac{d \gamma_{\varepsilon}}{d\tau_\epsilon}&=& \nu_\varepsilon\,(\gamma'(t)+\varepsilon  \gamma_1'( t)+O(\varepsilon^2))\\
		&=& (1-\varepsilon \nu_1+ O(\varepsilon^2))\,(\gamma'(t)+\varepsilon  \gamma_1'( t)+O(\varepsilon^2))\\
		&=& \gamma'( t)-\varepsilon \nu_1 \gamma'( t)  +\varepsilon \gamma_1'( t) + O(\varepsilon^2) \\
		&=& F_0(\gamma(t)) + \varepsilon(-\nu_1 F_0(\gamma(t)) +\gamma_1'( t))+O(\varepsilon^2),}
	\end{split}
	\end{equation}
	where $'$ denotes the derivative with respect to $t$.
	On the other hand, expanding the right hand side of \eqref{eq:dxdt=Feps} gives
	\begin{equation}\label{eq:dxdtau2}
	\begin{split}
	\Eqn{\frac{d \gamma_{\varepsilon}}{d\tau_\epsilon} &=& F_\varepsilon(\gamma_{\varepsilon}(\tau_\epsilon))\\
		&=& F_0(\gamma(t))+\varepsilon \Big(DF_0(\gamma(t)) \gamma_1(t) +\frac{\partial F_\varepsilon(\gamma(t))}{\partial \varepsilon}\Big|_{\varepsilon=0}\Big)+O(\varepsilon^2).}
	\end{split}
	\end{equation}
	Equating \eqref{eq:dxdtau1} and \eqref{eq:dxdtau2} to first order, we find that the linear shift in shape produced by a static perturbation, i.e.~the iSRC, satisfies 
	\begin{eqnarray}\label{eq:src}
	\frac{d \gamma_1(t)}{dt}
	&=& DF_0(\gamma(t)) \gamma_1(t) +\nu_1 F_0(\gamma(t)) +\frac{\partial F_\varepsilon(\gamma(t))}{\partial \varepsilon}\Big|_{\varepsilon=0},
	\end{eqnarray}
	with period $T_0$, as claimed before.
	
	It remains to establish an initial condition $\gamma_1(0)$. \RED{For smooth systems, without loss of generality, we may choose the initial condition for \eqref{eq:src} by taking a Poincar\'{e} section orthogonal to the limit cycle at a chosen reference point $p_0=\gamma(0)$. Then the initial condition is given by $\gamma_1(0)=(p_\varepsilon-p_0)/\varepsilon$ where $p_\varepsilon$ is the intersection point where the perturbed limit cycle $\gamma_\varepsilon$ crosses the Poincar\'{e} section. For nonsmooth systems discussed in the balance of the paper, we choose the reference section $\Sigma$ to be one of the switching or contact boundaries. That such an arbitrary choice of an initial condition using an orthogonal Poincar\'{e} section or a switching boundary does not compromise generality is a consequence of the following Lemma, the proof of which is in given in Appendix \ref{ap:isrc-ini}. 
		\begin{lemma}\label{lem:isrc}
			Let $\gamma^\textbf{a}_1(t)$ and $\gamma^\textbf{b}_1(t)$ be
			two $T_0$-periodic solutions to the iSRC equation \eqref{eq:src} for a smooth vector field $F_0$ with a hyperbolically stable limit cycle $\gamma_0(t)$.  
			Then, their difference satisfies $\gamma^\textbf{b}_1(t)-\gamma^\textbf{a}_1(t)=\varphi F_0(\gamma_0(t))$, where $\varphi$ is a constant representing a fixed phase offset.
		\end{lemma}
		Thus, if $\gamma^\textbf{a}_1(t)$ is a representative iSRC specified by taking the orthogonal Poincar\'{e} section, and $\gamma^\textbf{b}_1(t)$ is another representative specified by a different transverse section, then the two solutions differ by a fixed offset -- namely a vector in the direction of the flow along the limit cycle -- indexed by an additive difference in phase. Hence the differences between distinct periodic solutions to \eqref{eq:src} have precisely the same degree of ambiguity -- and for the same reason -- as the familiar ambiguity of the phase of an oscillator. }

	The accuracy of the iSRC in approximating the linear change in the limit cycle shape evidently depends on its timing sensitivity, that is, the choice of the relative change in frequency $\nu_1$.
	In the preceding derivation, we chose $\nu_1$ to be the relative change in the full period by assuming the limit cycle has constant timing sensitivity.
	It is natural to expect that different choices of $\nu_1$ will be needed for systems with varying timing sensitivities along the limit cycle.
	This possibility motivates us to consider local timing surfaces which divide the limit cycle into a number of segments, each distinguished by its own timing sensitivity properties.
	For each segment, we show that the linear shift in the time that $\gamma(t)$ spends in that segment can be estimated using a \emph{local timing response curve} (lTRC) derived in \S\ref{sec:smooth-local}.
	The lTRC is analogous to the iPRC in the sense that they obey the same adjoint equation, but with different boundary conditions.

	\subsection{Local timing response to perturbations: lTRC}\label{sec:smooth-local}  
	
	The iPRC captures the net effect on timing of an oscillation -- the phase shift -- due to a transient perturbation \eqref{eq:definition_of_iPRC}, as well as the net change in period due to a sustained perturbation \eqref{eq:T1}.
	In order to study the impact of a perturbation on \emph{local} timing as opposed to global timing, we introduce the notion of \textit{local timing surfaces} that separate the limit cycle trajectory into segments with different timing sensitivities.
	Examples of local timing surfaces in smooth systems include the passage of neuronal voltage through its local maximum or through a predefined threshold voltage, and the point of maximal extension  of reach by a limb.
	In nonsmooth systems, switching surfaces at which dynamics changes can also serve as local timing surfaces.
	For instance, in the feeding system of \textit{Aplysia californica} \citep{shaw2015,lyttle2017}, the open-closed switching boundary of the grasper defines a local timing surface.\footnote{\RED{It is worth noting that the idea of exploiting the presence of timing surfaces to specify the lTRC, by which we are taking advantage of corresponding features present in both the perturbed and unperturbed limit cycles, may be seen as an example of the  general notion of \emph{bisimulation} \citep{HTP2003,HTP2005}.}}

\begin{figure}[!t]
		\begin{center}
			\includegraphics[width=3in]{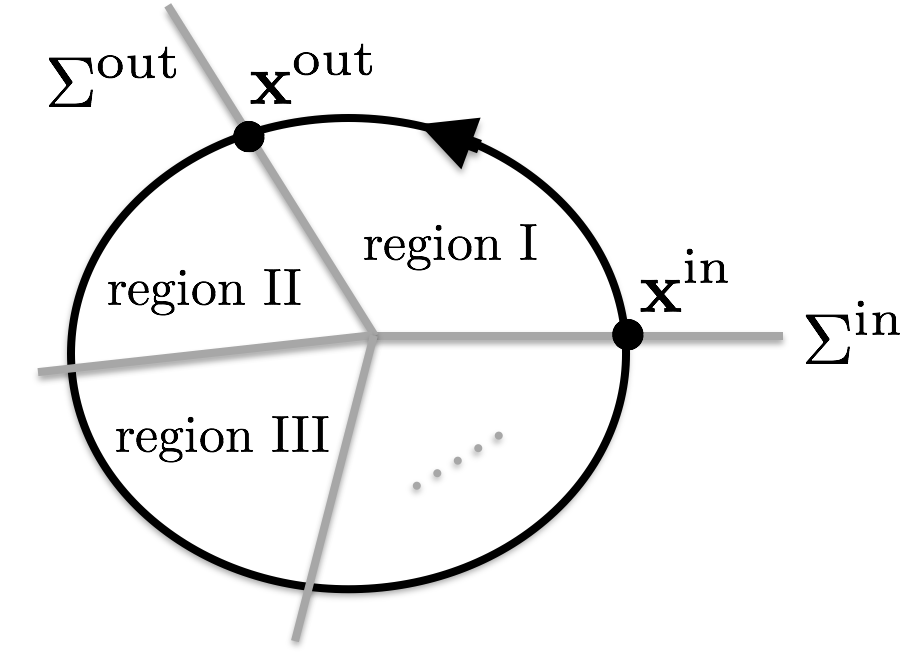}
		\end{center}
		\caption{\label{fig:local-time-response} Schematic illustration of a limit cycle solution for a system consisting of a number of regions, each with distinct constant timing sensitivities.
			$\Sigma^{\rm in}$ and $\Sigma^{\rm out}$ denote the local timing surfaces for region I.
			$\mbx^{\rm in}$ and $\mbx^{\rm out}$ denote the points where the limit cycle enters region I through $\Sigma^{\rm in}$ and exits region I through $\Sigma^{\rm out}$, respectively.}
	\end{figure}

	Whatever the origin of the local timing surface or surfaces of interest, it is natural to consider the phase space of a limit cycle as divided into multiple regions.
	Hence we may consider a smooth system $d\mbx/dt=F(\mbx)$ with a limit cycle solution $\gamma(t)$ passing through multiple regions in succession (see Figure \ref{fig:local-time-response}).
	In each region, we assume that $\gamma(t)$ has constant timing sensitivity.
	To compute the relative change in time in any given region, we define a lTRC to measure the timing shift of $\gamma(t)$ in response to perturbations delivered at different times in that region.
	Below, we illustrate the derivation of the lTRC in region I and show how it can be used to compute the relative change in time in this region, denoted by $\nu_1^{\rm I}$.
	
	Suppose that at time $t^{\rm in}$, $\gamma(t)$ enters region I upon crossing the surface $\Sigma^{\rm in}$ at the point $\mbx^{\rm in}$; at time $t^{\rm out}$, $\gamma(t)$ exits region I upon crossing the surface $\Sigma^{\rm out}$ at the point $\mbx^{\rm out}$ (see Figure \ref{fig:local-time-response}).
	Denote the vector field under a constant perturbation by $F_\varepsilon(\mbx)$ and let $\mbx_{\varepsilon}$ denote the coordinate of the perturbed trajectory.
	Let $T_0^{\rm I} = t^{\rm out}-t^{\rm in}$ denote the time $\gamma(t)$ spent in region I and let $T_{\varepsilon}^{\rm I}$ denote the time the perturbed trajectory spent in region I.
	Assume we can write $T_{\varepsilon}^{\rm I}=T_0^{\rm I}+\varepsilon T^{\rm I}_1+O(\varepsilon^2)$ as we did before.
	It follows that the relative change in time of $\gamma(t)$ in region I is given by 
	\[
	\nu_1^{\rm I}=\frac{T^{\rm I}_1}{T_0^{\rm I}} = \frac{T^{\rm I}_1}{t^{\rm out}-t^{\rm in}} .
	\] 
	The goal is to compute $\nu_1^{\rm I}$, which requires an estimate of $T^{\rm I}_1$.
	To this end, we define the local timing response curve $\eta^{\rm I}(t)$ associated with region I.
	We show that $\eta^{\rm I}(t)$ satisfies the adjoint equation \eqref{eq:ltrc} and the boundary condition \eqref{eq:ltrc0}.
	
	Let $\mathcal{T}^{\rm I}(\mbx)$ for $\mbx$ in region I be the time remaining until exiting region I through $\Sigma^{\rm out}$, under the unperturbed vector field.
	This function is at least defined in some open neighborhood around the reference limit cycle trajectory $\gamma(t)$ if not throughout region I.
	For the unperturbed system, $\mathcal{T}^{\rm I}$  satisfies 
	\begin{equation}\frac{d\mathcal{T}^{\rm I}(\mbx(t))}{dt}=-1\end{equation}
	along the limit cycle orbit $\gamma(t)$.
	Hence 
	\begin{equation}\label{eq:local-time-const}
	F(\mbx)\cdot \nabla\mathcal{T}^{\rm I}(\mbx)=-1
	\end{equation}
	for all $\mbx$ for which $\mathcal{T}^{\rm I}$ is defined.
	We define $\eta^{\rm I}(t):= \nabla \mathcal{T}^{\rm I}(\mbx(t))$ to be the local timing response curve (lTRC) for region I.
	It is defined for $t\in [t^{\rm in}, t^{\rm out}]$.
	We show in Appendix \ref{ap:T1} that $T^{\rm I}_1$ can be estimated as
	\begin{equation}\label{eq:local-time-shift}
	T^{\rm I}_{1} = \eta^{\rm I}(\mbx^{\rm in})\cdot \frac{\partial \mbx_{\varepsilon}^{\rm in}}{\partial \varepsilon}\Big|_{\varepsilon=0}+\int_{t^{\rm in}}^{t^{\rm out}}\eta^{\rm I}(\gamma(t))\cdot \frac{\partial F_\varepsilon(\gamma(t))}{\partial \varepsilon}\Big|_{\varepsilon=0}dt,
	\end{equation} 
	where $\mbx_\varepsilon^{\rm in}$ denotes the coordinate of the perturbed entry point into region I.
	We may naturally view $\eta^{\rm I}$ as either a function of space, as in \eqref{eq:local-time-shift}, or as a function of time, evaluated e.g.~along the limit cycle trajectory.
	Comparing \eqref{eq:local-time-shift} with \eqref{eq:T1}, the integral terms have the same form, albeit with opposite signs.
	In addition, \eqref{eq:local-time-shift} has an additional term arising from the impact of the perturbation on the point of entry to region I.
	On the other hand, the impact of the perturbation on the exit point, denoted by $\eta^{\rm I}(\mbx_\varepsilon^{\rm out})\cdot \frac{\partial \mbx_\varepsilon^{\rm out}}{\partial \varepsilon}\big|_{\varepsilon=0}$, is always zero because the exit boundary $\Sigma^{\rm out}$ is a level curve of $\mathcal{T}^{\rm I}$; in other words, $\mathcal{T}^{\rm I}\equiv 0$ at $\Sigma^{\rm out}$.
	This indicates that the lTRC vector $\eta^{\rm I}$ associated with a given region is always perpendicular to the exit boundary of that region.
	
	Similar to the iPRC, it follows from \eqref{eq:local-time-const} that $\eta^{\rm I}$ satisfies the adjoint equation
	\begin{eqnarray}\label{eq:ltrc}
	\frac{d\eta^{\rm I}}{dt}=-DF(\gamma(t))^\intercal \eta^{\rm I}
	\end{eqnarray}
	together with the boundary (normalization) condition at the exit point
	\begin{equation}\label{eq:ltrc0}
	\eta^{\rm I}(\mbx^{\rm out})=\frac{-n^{\rm out}}{n^{\rm out \intercal} F(\mbx^{\rm out})}
	\end{equation}
	where $n^{\rm out}$ is a normal vector of $\Sigma^{\rm out}$ at the unperturbed exit point $\mbx^{\rm out}$.
	The reason $\eta^{\rm I}$ at the exit point has the direction $n^{\rm out}$ is because $\eta^{\rm I}$ is normal to the exit boundary as discussed above.
	
	To summarize, in order to compute $\nu_1^{\rm I}=T^{\rm I}_{1}/({t^{\rm out}-t^{\rm in}})$, we need numerically to find $t^{\rm in},t^{\rm out}$ and evaluate~\eqref{eq:local-time-shift} to estimate $T^{\rm I}_{1}$, for which we need to solve the boundary problem of the adjoint equation~\eqref{eq:ltrc}-\eqref{eq:ltrc0} for the lTRC $\eta^{\rm I}$.
	The procedures to obtain the relative change in time in other regions $\nu_1^{\rm II}, \nu_1^{\rm III},\cdots$ are similar to computing $\nu_1^{\rm I}$ in region I and hence are omitted.
	The existence of different timing sensitivities of $\gamma(t)$ in different regions therefore leads to a piecewise-specified version of the iSRC \eqref{eq:src} with period $T_0$, 
	\begin{eqnarray}\label{eq:src-multiplescales}
	\frac{d \gamma_1^j(t)}{dt}
	&=& DF_0^j(\gamma(t)) \gamma_1^j(t) +\nu_1^{j} F_0^j(\gamma(t)) +\frac{\partial F_\varepsilon^j(\gamma(t))}{\partial \varepsilon}\Big|_{\varepsilon=0},
	\end{eqnarray}
	where $\gamma_1^j$, $F_0^j$, $F_{\varepsilon}^j$ and $\nu_1^j$ denote the iSRC, the  unperturbed vector field, the perturbed vector field, and the relative change in time in region $j$, respectively, with $j\in\{\rm I, II, III, \cdots\}$.
	Note that in a smooth system as concerned in this section, $F_0^j\equiv F_0$ for all $j$.
	
	In \S\ref{sec:toy-model} we will show in a specific example that the iSRC with piecewise-specified timing rescaling has much greater accuracy in approximating the linear shape response of the limit cycle to static perturbations than the iSRC using a global uniform rescaling.

	\begin{remark}
		The derivation of the lTRC in a given region still holds as long as the system is smooth in that region.
		Hence the assumption that $F(\mbx)$ is smooth everywhere can be relaxed to $F(\mbx)$ being piecewise smooth.
	\end{remark}
	
	\RED{
		\begin{remark}
			The lTRC is an intrinsic property of a limit cycle that is determined by the choice of local timing surfaces that are transverse to the limit cycle. In this work, we only use the lTRC in the situation where there are naturally occurring surfaces such as those corresponding to switching surfaces inherent in the geometry of a problem (e.g., the stick-slip system).  In general, the lTRC construction could be used to analyze smooth systems as well, for instance systems with heteroclinic cycling \citep{shaw2012}. 
	\end{remark}}

	\section{Linear responses of nonsmooth systems with continuous solutions}
	\label{sec:nonsmooth-theory}
	
	Nonsmooth dynamical systems arise in many areas of biology and engineering.
	However, methods developed for smooth systems (discussed in \S\ref{sec:smooth-theory}) do not extend directly to understanding the changes of periodic limit cycle orbits in nonsmooth systems, because their Jacobian matrices are not well defined \citep{CGL18,wilson2019}.
	Specifically, nonsmooth systems exhibit discontinuities in the time evolution of the solutions to the variational equations, $\mbu$ \eqref{eq:var} and $\gamma_1$ \eqref{eq:src}, and the solutions to the adjoint equations, $\z$ \eqref{eq:prc} and $\eta$ \eqref{eq:ltrc}.
	Following the terminology of \citet{park2018} and \citet{LN2013}, we call the discontinuities in $\z$ and $\eta$ ``jumps'' and call the discontinuities in $\mbu$ and $\gamma_1$ ``saltations''.
	Qualitatively, we use ``jumps" to refer to discontinuities in the \emph{timing} response of a trajectory, and ``saltations" in the \emph{shape} response.
	Since $\z$ and $\eta$ satisfy the same adjoint equation, they have the same discontinuities.
	Similarly, $\mbu$ and $\gamma_1$ obey versions of the variational equation  with the same  homogeneous term and different nonhomogeneous terms; since the jump conditions arise from the homogeneous terms (involving the Jacobian matrix in \eqref{eq:var}, \eqref{eq:src}) we will presume that $\mbu$ and $\gamma_1$ satisfy the same saltation conditions at the transition boundaries.
	In this section, we characterize the discontinuities in the solutions to the adjoint equation in terms of $\z$, and discuss nonsmoothness of the variational dynamics in terms of $\mbu$.

	\RED{In this paper, we consider nonsmooth systems with degree of smoothness one or higher (Filippov systems); that is, systems with continuous solutions.
		In such systems, the right-hand-side changes discontinuously as one or more \textit{switching surfaces} are crossed. 
		A trajectory reaching a switching surface has two behaviors: it may cross the surface transversally, or it may slide along it, in which case the motion is called \textit{sliding mode}.  
		In \citep{shirasaka2017,park2018,CGL18,wilson2019,bernardo2008,LN2013}, the solutions to the adjoint and variational equations have been studied in the case of solutions which cross the surface of discontinuity transversally. 
		The concept of saltations in the variational dynamics $\mbu$ has also been adapted by \citep{filippov1988,bernardo2008,LN2013,DL11} to the case of sliding motion on surfaces. 
		However, in the case of sliding mode motions, discontinuous jumps exhibited by the iPRC have not yet been characterized. 
		Bridging that gap is a principal goal of this paper. 
		
		In the remainder of this section, we first define Filippov systems with both types of discontinuities. 
		Then we review the existing  methods for computing the saltations in $\mbu$ (and $\gamma_1$) in both cases, and the jumps in $\z$ (and $\eta$) in the case of transversal crossing boundaries, following \cite{bernardo2008,LN2013} and \cite{park2013,park2018}. Lastly in \S\ref{sec:sliding}, we present our main results about the discontinuous behavior of the iPRC for nonsmooth systems with sliding motions. For completeness, we include both the old and new results in the case of sliding motions in the statement of Theorem \ref{thm:main}. }

	\begin{figure}[!t]
		\begin{center}
			\setlength{\unitlength}{1mm}
			\begin{picture}(145,60)(0,0)
			\put(0,0){\includegraphics[height=5cm]{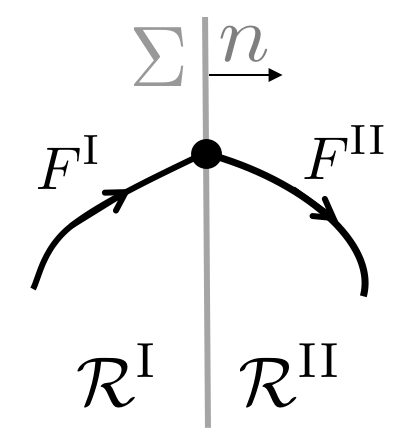}}
			
			\put(-3,45){\bf{\large{(A)}}}
			
			\put(80,0){\includegraphics[height=5cm]{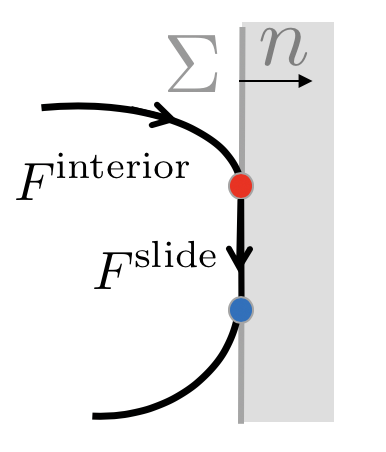}}
			\put(80,45){\bf{\large{(B)}}}
			\end{picture}
		\end{center}
		\caption{\label{fig:fillippov}
			Examples of trajectories of nonsmooth systems with a boundary $\Sigma$.
			(A) The trajectory of a two-zone nonsmooth system \eqref{eq:2zoneFS} intersects the boundary $\Sigma$ transversely at $\mbx_p$ (black dot).
			$F^{\rm I}$ and $F^{\rm II}$ denote the vector fields in the two regions $\RR^{\rm I}$ and $\RR^{\rm II}$.
			Components of $F^{\rm I}$ and $F^{\rm II}$ normal to $\Sigma$ at the crossing point have the same sign, allowing for the transversal crossing.
			(B) The trajectory of a nonsmooth system \eqref{eq:1zoneFP} hits the hard boundary $\Sigma$ at the landing point (red dot) and begins sliding along $\Sigma$ under the vector field $F^{\rm slide}$.
			At the liftoff point (blue dot), the trajectory naturally reenters the interior.
			$F^{\rm interior}$ denotes the vector field in the interior domain.
		}
	\end{figure}
	
	\RED{	\subsection{Filippov systems}\label{sec:filippov}
		
		It is sufficient to consider a Filippov system with a single boundary (switching surface) to illustrate the discontinuities of $\z$ and $\mbu$ at any boundary crossing point. Below we give the definition of a two-zone Filippov system with a transversal crossing boundary (see \eqref{eq:2zoneFS}) and a local representation of a Filippov system that exhibits the sliding mode (see \eqref{eq:1zoneFP}).  
		
		\subsubsection{Transversal crossing boundary}	
		\begin{definition}\label{def:2zoneFP}\rm 
			A \textit{two-zone system} with uniform degree of smoothness one (or higher) is described by
			\begin{eqnarray}\label{eq:2zoneFS}
			\frac{d\textbf{x}}{dt}=F(\mbx):=\left\{
			\Eqn{
				F^{\rm I}(\mbx),  & \mbx \in \RR^{\rm I} \\
				F^{\rm II}(\mbx), & \mbx \in \RR^{\rm II} \\
			}
			\right.
			\end{eqnarray}
			where $\RR^{\rm I}:=\{\mbx | H(\mbx)<0\}$ and $\RR^{\rm II}:=\{\mbx | H(\mbx)>0\}$ for a smooth function $H$ have non-empty interiors, and the vector fields $F^{\rm I,II}:\overline{\RR}^{\rm I,II}\to \mathbb{R}^n$ are at least $C^1$, where $\overline{\RR}$ denotes the closure of $\RR$ in $\mathbb{R}^n$.
		\end{definition}
		\begin{definition} \label{def:2zoneFP-boundary}\rm 
			The \emph{switching boundary} for the Filippov system \eqref{eq:2zoneFS} is the $\mathbb{R}^{n-1}$-dimensional manifold $\Sigma:= \overline{\RR}^{\rm I}\cap \overline{\RR}^{\rm II}=\{\mbx|H(\mbx)=0\}$. $\Sigma$ is called a \textit{transversal crossing boundary} if at $\mbx_p\in \Sigma$ the following holds
			\begin{eqnarray*}
				\Eqn{(n_{p}\cdot F^{\rm I}(\mbx_p)) (n_{p} \cdot F^{\rm II}(\mbx_p))> 0} 
			\end{eqnarray*}
			where $n_p=\nabla H(\mbx_p)$ refers to the vector normal to $\Sigma$ at $\mbx_p$.
		\end{definition}
		
		Without loss of generality, we assume at time $t=t_p$ a limit cycle solution $\gamma(t)$ of \eqref{eq:2zoneFS} crosses the boundary $\Sigma$ from $\RR^{\rm I}$ to $\RR^{\rm II}$ (see Figure \ref{fig:fillippov}A). In this case, we have
		\begin{eqnarray}\label{eq:FP-trans}
		\Eqn{n_{p}\cdot F^{\rm I}(\mbx_p)> 0\quad\quad n_{p} \cdot F^{\rm II}(\mbx_p)> 0, }
		\end{eqnarray}
		where the boundary crossing point can be defined as $\mbx_p:=\lim_{t\to t_p^-}\gamma(t)=\lim_{t\to t_p^+}\gamma(t)$. 
		
		\subsubsection{Sliding motion on a hard boundary} 
		
		Next we consider the second type of switching surface on which the transversal flow condition \eqref{eq:FP-trans} does not hold. That is, parts of the solution trajectory slide along a surface (e.g., Figure \ref{fig:fillippov}B). As an example of a hard boundary at which the transverse flow condition would break down, consider the requirement that firing rates in a neural network model be nonnegative.
		When a nerve cell ceases firing because of inhibition, its firing rate will be held at zero until the  balance of inhibition and excitation allow spiking to resume.
		At the point at which the firing rate first resumes positive values, the vector field describing the system lies tangent to the constraint surface rather than transverse to it.
		
		Below we present a model of a Filippov system with sliding motion along a hard boundary - system \eqref{eq:1zoneFP}. 
		We begin with precise definitions of \textit{hard boundary, sliding region, sliding vector field} and the \textit{liftoff condition}. In such systems, non-transversal crossing points include the \textit{landing point} at which a sliding motion begins, and the \textit{liftoff point} at which the sliding terminates (see Fig.~\ref{fig:fillippov}B).

		\begin{definition}\label{def:sigma}\rm 
			Consider a system with domain $\mathcal{R}$.
			We call a surface $\Sigma$ a \emph{hard boundary} if it is part of the boundary of the closure of $\mathcal{R}$.
		\end{definition}
		
		\begin{definition}\label{def:sliding}\rm
			The \textit{sliding region} ($\RR^{\rm slide}$) is defined as the portion of a hard boundary $\Sigma$ for which 
			\begin{eqnarray}\label{eq:slidingregion}
			\Eqn{\RR^{\rm slide}=\{\mbx\in\Sigma\,|\,n_\mbx\cdot F^{\rm interior}(\mbx)>0\},}
			\end{eqnarray}
			where $n_\mbx$ is a unit normal vector of $\Sigma$ at $\mbx$ that points away from the interior and $F^{\rm interior}$ denotes the vector field defined on the closure of the domain $\RR$.
		\end{definition}
		
		\begin{definition}\label{def:interior}\rm
			The \textit{interior domain} ($\RR^{\rm interior}$) is defined as the complement of $\RR^{\rm slide}$ in the closure of $\RR$. 
		\end{definition}
		
		When the trajectory enters the sliding region, the solution will continue along $\RR^{\rm slide}$ with time derivative $F^\text{slide}$ that is tangent to the hard boundary.
		While any vector field with vanishing normal component could be considered for $F^\text{slide}$, in this paper we adopt the natural choice of setting $F^\text{slide}$
		to be the following.
		
		\begin{definition}\label{def:slidingvf} \rm 
			The \textit{sliding vector field} $F^{\rm slide}$, defined on $\RR^{\rm slide}$, is given by the continuation of $F^{\rm interior}$ in the component tangential to $\RR^{\rm slide}$:
			\begin{eqnarray}\label{eq:slidingVF}
			\Eqn{F^{\rm slide}(\mbx)=F^{\rm interior}(\mbx)-\big(n_\mbx\cdot F^{\rm interior}(\mbx)\big)n_\mbx.}
			\end{eqnarray} 
		\end{definition}
		During the sliding motion, the flow will slide along $\Sigma$ with the sliding vector field until it is allowed to reenter the interior; that is
		\begin{definition}\label{def:liftoff}\rm 
			The flow exits the sliding region $\RR^{\rm slide}$ as the trajectory crosses the
			\textit{liftoff boundary} $\mathcal{L}$ defined as
			\begin{equation}\label{eq:cond-liftoff}
			\Eqn{\mathcal{L}=\{\mbx\in\Sigma  \,| \,n_\mbx \cdot F^{\rm interior}(\mbx)=0\}.}
			\end{equation}
		\end{definition}
		Thus the liftoff boundary constitutes the edge of the sliding region of the hard boundary. To identify the liftoff point at which the trajectory reenters the interior of the domain, we further require the nondegeneracy condition that the trajectory crosses the liftoff boundary $\mathcal{L}$ at a finite velocity.
		Specifically, the outward normal component of the interior velocity should switch from positive (outward) to negative (inward) at the liftoff boundary, as one moves in the direction of the flow (see Fig.~\ref{fig:sliding-motion-coordinates}).
		That is,
		\begin{eqnarray}\label{eq:nondege-1}
		\Eqn{\left.\left[\nabla(n_\mbx\cdot F^{\rm interior}(\mbx)) \cdot F^{\rm slide}(\mbx)\right] \right|_{\mbx \in \mathcal{L}}< 0}.
		\end{eqnarray} 
		Note that the liftoff condition \eqref{eq:cond-liftoff} together with the nondegeneracy condition \eqref{eq:nondege-1}  uniquely defines a liftoff point for the trajectory that slides along $\Sigma$.
		At the liftoff point, we have $F^{\rm slide}=F^{\rm interior}$.}
	
	\begin{remark}
		Our definition of the sliding region and sliding vector field is consistent with that in \cite{bernardo2008} \S 5.2.2, except that the system of interest in this paper is only defined on one side of the sliding region.
		However, our main Theorem \ref{thm:main}, below, holds in either case.
		Hence our results also apply to Filippov systems with sliding regions bordered by vector fields on either side, as in  the example the stick-slip oscillator (\citet{LN2013} \S 6.5).
	\end{remark}
	
	The motion along a trajectory is specified differently depending on the location of a  point.
	For a point in the interior, the dynamics is determined by  $F^\text{interior}$.
	For a point on $\Sigma$, the velocity obeys either $F^\text{interior}$ or else $F^\text{slide}$ that is tangent to $\Sigma$, depending on whether $F^\text{interior}$ is directed inwardly or outwardly at a given boundary point.
	This dual definition of the vector field has the effect that points driven into the boundary do not exit through the hard boundary, but rather slide along the boundary until the trajectory crosses the liftoff boundary $\mathcal{L}$.

	Using the preceding notation, in the neighborhood of a hard boundary $\Sigma$, a system with a limit cycle component confined to the sliding region takes the following form
	\begin{eqnarray}\label{eq:1zoneFP}
	\frac{d\textbf{x}}{dt}=F(\mbx):=\left\{
	\Eqn{
		F^{\rm interior}(\mbx), & \mbx \in \RR^{\rm interior}\\
		F^{\rm slide}(\mbx), & \mbx \in\RR^{\rm slide}\subset \Sigma \\
	}
	\right.
	\end{eqnarray}
	where $F^{\rm interior},\RR^{\rm interior}, F^{\rm side}, \RR^{\rm slide}, \Sigma$ are given in Definitions \ref{def:sigma}-\ref{def:slidingvf}.
	
	\begin{definition}\label{def:lcsc}\rm
		In a general Filippov system which locally at a hard boundary $\Sigma$ has the form \eqref{eq:1zoneFP}, we call a closed, isolated periodic orbit that passes through a sliding region a \textit{limit cycle with sliding component}, denoted as \rm LCSC.\end{definition}

	\begin{figure}[!t]
		\centering
		\includegraphics[width=0.7\textwidth]{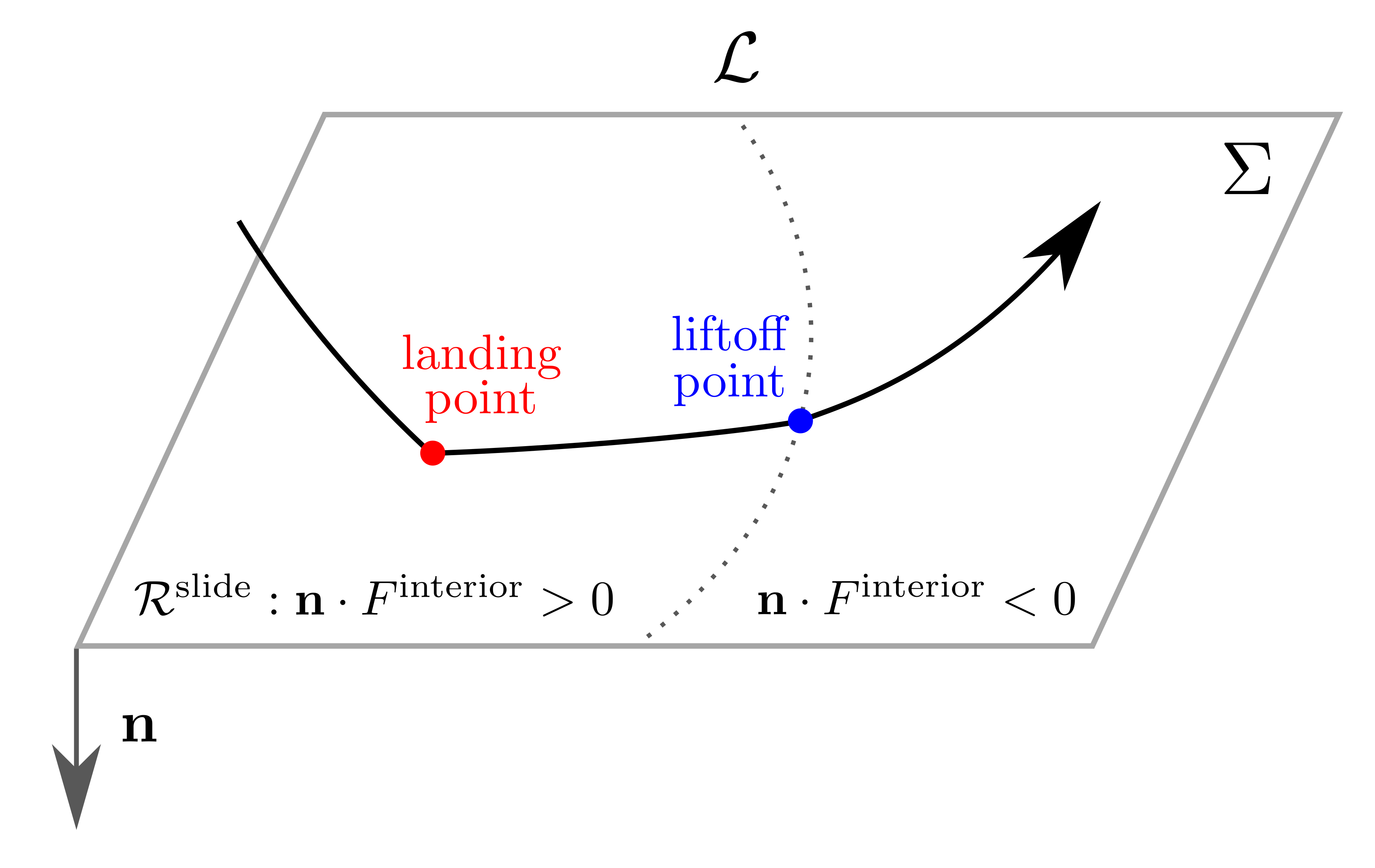}
		\caption{
			Trajectory from the interior with vector field $F^{\rm interior}$ making a transverse entry to a hard boundary $\Sigma$, followed by motion confined to the sliding region $\mathcal{R}^{\rm slide}$, then a smooth liftoff at $\mathcal{L}$ back into the interior of the domain.
			Red dot: landing point (point at which the trajectory exits the interior and enters the hard boundary surface).
			Blue dot: liftoff point (point at which the trajectory crosses the liftoff boundary $\mathcal{L}$ and reenters the interior).
			After a suitable change of coordinates, the geometry may be arranged as shown, with the hard boundary $\Sigma$ coinciding with one coordinate plane.
			Downward vertical arrow: $\mathbf{n}$, the outward normal vector for $\Sigma$.
			The region $\mathbf{n}\cdot F^\text{interior}>0$ defines the sliding region within $\Sigma$; the condition $\mathbf{n}\cdot F^\text{interior}=0$ defines the liftoff boundary $\mathcal{L}$. }
		\label{fig:sliding-motion-coordinates}
	\end{figure}
	
	\RED{We make the following assumptions about the vector field $F$, the hard boundary $\Sigma$ and the liftoff boundary $\LL\subset \Sigma$ throughout:
		\begin{assumption}\label{ass:assumption1}
			\begin{enumerate}
				\item $F^{\rm interior}: \overline{\RR}^{\rm interior}\to \mathbb{R}^n$ is at least $C^1$.
				\item Under an appropriate smooth change of coordinates, the hard boundary $\Sigma$ can be transformed into a lower dimensional manifold with a constant normal vector $n$  (cf.~Fig.~\ref{fig:sliding-motion-coordinates}).
				\item When the trajectory crosses the liftoff boundary $\LL$, the nondegeneracy condition \eqref{eq:nondege-1} holds so a liftoff point can be uniquely defined. 
			\end{enumerate}
	\end{assumption}}
	
	\subsection{Review of variational dynamics and iPRC in Filippov systems}\label{sec:nonsmooth-review}
	
	\RED{Below we review the behaviors of the variational dynamics $\mbu$ and the iPRC $\z$ in the case of transversal intersection, as well as the characterization of discontinuities in $\mbu$ in the case of sliding motion on a hard boundary. 
		As discussed before, the iSRC $\gamma_1$ and the lTRC $\eta$ experience the same discontinuities as $\mbu$ and $\z$, respectively.}
	
	\subsubsection{Variational dynamics: transversal crossing and sliding motion}
	
	For a sufficiently small instantaneous perturbation, the displacement $\mbu(t)$ evolves continuously over the domain in which \eqref{eq:2zoneFS} is smooth, and can be obtained to first order in the initial displacement by solving the variational equation \eqref{eq:var}.
	As $\gamma$ crosses $\Sigma$ at time $t_p$, $\mbu(t)$ exhibits discontinuities (or ``saltations'') since the Jacobian evaluated at $\mbx_p$ is not uniquely defined.
	The discontinuity in $\mbu$ at $\mbx_p$ can be expressed with the \textit{saltation matrix $S_p$}  as
	\begin{eqnarray}\label{eq:salt-def}
	{\mbu_p^+=S_p\,\mbu_p^-}
	\end{eqnarray}
	where $\mbu_p^-=\lim_{t\to t_p^-} \mbu(t)$ and $\mbu_p^+=\lim_{t\to t_p^+} \mbu(t)$ represent the displacements between perturbed and unperturbed solutions just before and just after the crossing, respectively.
	It is straightforward to show (cf.~\citet{LN2013} \S 7.2 or \citet{bernardo2008} \S 2.5) that
	$S_p$ can be constructed using the vector fields in the neighborhood of the crossing point and the vector $n_p$ normal to the switching boundary at $\mbx_p$ as
	\begin{equation}\label{eq:salt}
	S_p=I+\frac{(F_p^+-F_p^-)n_{p}^\intercal}{n_{p}^\intercal F_p^-}
	\end{equation} 
	where $F_p^-=\lim_{\mbx\to \mbx_p^-}F(\mbx),\,F_p^+=\lim_{\mbx\to \mbx_p^+}F(\mbx)$ are the vector fields of \eqref{eq:2zoneFS} just before and just after the crossing at $\mbx_p$.
	Throughout this paper, $I$ denotes the identity matrix with size $n\times n$.
	
	\begin{remark}\label{remark-salt}
		If the vector field $F$ evaluated along the limit cycle is continuous when crossing the boundary $\Sigma$ transversely, so that $F_p^-=F_p^+$ and $n_p^\intercal F_p^-\neq 0$, then the saltation matrix $S_p$ at such a boundary crossing point is the identity matrix, and there is no discontinuity in $\mbu$ or $\gamma_1$ at time $t_p$.
	\end{remark}   
	
	\RED{ 
		When the transversal flow condition $n_p^\intercal F_p^-\not=0$ is violated, the expression \eqref{eq:salt} cannot be used.
		Nevertheless,
		\citep{filippov1988,bernardo2008} showed how to adapt the definition of the saltation matrix to capture the discontinuity of $\mbu$ at non-transversal crossings including both landing and liftoff points. 
		Specifically, at a landing point, the saltation matrix is given by 
		\[
		S_p = I - n_{p}n_p^\intercal,
		\]
		whereas at a liftoff point, the saltation matrix is 
		\[
		S_p = I.
		\]
		The above results about $S_p$ at the non-transversal crossing points are also summarized as parts (a) and (b) in Theorem \ref{thm:main}, together with our new results about the iPRC in the case of sliding motion along a hard boundary.  For the derivation of $S_p$ at the landing and liftoff points, we refer to \citep{filippov1988,bernardo2008} or our proof of Theorem \ref{thm:main} in Appendix \ref{ap:proof}. }

	\subsubsection{iPRC: transversal crossing}
	
	Now we consider discontinuous jumps in the iPRC $\z$ for \eqref{eq:2zoneFS} with transversal intersection.
	This curve obeys the adjoint equation \eqref{eq:prc} and is continuous within the interior of each subdomain in which \eqref{eq:2zoneFS} is smooth.
	When the limit cycle path crosses the switching boundary at the point $\mbx_p$, $\z$ exhibits a discontinuous jump which can be characterized by the \textit{jump matrix} $J_p$ 
	\begin{eqnarray}\label{eq:jump-def}
	\z_p^+=J_p\z_p^-
	\end{eqnarray}
	where $\z_p^-=\lim_{t\to t_p^-} \z(t)$ and $\z_p^+=\lim_{t\to t_p^+} \z(t)$ are the iPRC just before and just after crossing the switching boundary at time $t_p$.
	As discussed in \citet{park2018}, the relation \eqref{eq:var-prc} between $\mbu$ and $\z$ for smooth systems remains valid at any transversal boundary crossing point.
	In other words,  $\mbu_p^{-\intercal}\z_p^-=\mbu_p^{+\intercal}\z_p^+$ holds at the transversal crossing point $\mbx_p$.
	This leads to a relation between the saltation and jump matrices at $\mbx_p$
	\begin{equation}\label{eq:jump-salt-relation}
	J_p^\intercal S_p=I.
	\end{equation}
	The saltation matrix $S_p$ given by \eqref{eq:salt} has full rank at any transverse crossing point.
	It follows that $J_p$ can be written as
	\begin{equation}\label{eq:jump-trans}
	J_p = ({S_p}^{-1})^\intercal.
	\end{equation}
	\RED{The definition of $J_p$ given by \eqref{eq:jump-trans}, however, does not hold at a landing or a liftoff point. See Remarks \ref{rem:uztable} and \ref{rem:sala-jump} for more details. This motivates us to characterize discontinuous jumps in $\z$ at non-transversal crossing points, which will then allow us to compute iPRC in the context of LCSC. }
	
	\subsection{Jumps in iPRC: sliding motion on a hard boundary}\label{sec:sliding}
	
	\RED{In this section,  we now present our new results on the iPRC for LCSC.} As discussed in the previous section, the existence of the jump matrix \eqref{eq:jump-trans} is guaranteed by the transversal flow condition \eqref{eq:FP-trans}, which, however, will no longer hold when part of the limit cycle slides along a boundary (e.g., Figure \ref{fig:fillippov}B).
	Next, we establish the conditions relating $\z$ at non-transversal crossings, including the \textit{landing point} and the \textit{liftoff point}.
	\RED{Before that, we need to impose the following assumption on the asymptotic phase function for LCSC:
		
		\begin{assumption}\label{ass:assumption2}
			Within the stable manifold of a LCSC there is a well defined asymptotic phase function $\phi(\mbx)$ satisfying $d\phi/dt=1$ along the trajectory, where $\phi$ is Lipschitz continuous. Moreover, on the hard boundary $\Sigma$, the directional derivatives of $\phi$ with respect to directions tangential to the surface are also Lipschitz continuous, except (possibly) at the liftoff and landing points.
		\end{assumption}
		
		Our main result, Theorem \ref{thm:main} as shown below, gathers together several conclusions about the variational and infinitesimal phase response curve dynamics of a LCSC local to a sliding boundary. For completeness, we review the established behavior of variational dynamics local to a sliding boundary (parts (a) and (b)) and present our new results regarding the iPRC in parts (c) through (e) (see Appendix \ref{ap:proof} for the proof).	
		
		\begin{theorem}\label{thm:main}\rm 
			Consider a general LCSC described locally by~\eqref{eq:1zoneFP} in the neighborhood of a hard boundary $\Sigma$, satisfying Assumption \ref{ass:assumption1} and Assumption \ref{ass:assumption2}. The following properties hold for the variational dynamics $\mbu$ and the iPRC $\z$ along $\Sigma$:
			\begin{enumerate}
				\item[(a)] At the landing point of $\Sigma$, the saltation matrix is $S=I-n n^\intercal$, where $I$ is the
				identity matrix.
				\item[(b)] At the liftoff point of $\Sigma$, the saltation matrix is $S=I$.
				\item[(c)] Along the sliding region within $\Sigma$, the component of $\z$ normal to $\Sigma$ is zero.
				\item[(d)] The normal component of $\z$ is continuous at the landing point.
				\item[(e)] The tangential components of $\z$ are continuous at both landing and liftoff points.
			\end{enumerate}
	\end{theorem} 	}
	
	We make the following additional observations about Theorem \ref{thm:main}:
	
	\begin{remark}\label{rem:Jland=I} The following two statements follow directly from Theorem \ref{thm:main}:
		\begin{itemize}
			\item It follows from (a) in Theorem \ref{thm:main} that the component of $\mbu$ normal to $\Sigma$ vanishes when the LCSC hits $\Sigma$.
			Once on the sliding region, the Jacobian used in the variational equation switches from $DF^{\rm interior}$ to $DF^{\rm slide}$ where $F^{\rm slide}$ has zero normal component by construction \eqref{eq:slidingVF}.
			Hence, the normal component of $\mbu$ is stationary over time, and remains zero on the sliding region.
			\item It follows from parts (d) and (e) in Theorem \ref{thm:main} that the jump matrix of $\z$ at a landing point is trivial (identity matrix).
		\end{itemize}
	\end{remark}
	
	\begin{remark}\label{rem:iprc-liftoff}
		Theorem \ref{thm:main} excludes discontinuities in $\z$ \emph{except} at the liftoff point, and then only in its normal component.
		Since the normal component of $\z$ along each sliding component of a LCSC is zero by Theorem \ref{thm:main}, a discontinuous jump occurring at a liftoff point must be a nonzero instantaneous jump, which cannot be specified directly in terms of the value of $\z$ prior to the jump.
		However, a time-reversed version of the jump matrix at the liftoff point, denoted as $\mathcal{J}$, is well defined as follows:
		\begin{eqnarray}
		\Eqn{
			{\z_{\text{lift}}^{-}}=\mathcal{J} {\z_{\text{lift}}^{+}}
		}
		\end{eqnarray}
		where ${\z_{\text{lift}}^{-}}$ and ${\z_{\text{lift}}^{+}}$ 
		are the iPRC just before and just after the trajectory crosses the liftoff point in forwards time, and $\mathcal{J}$ at the liftoff point has the same form as the saltation matrix $S$ at the corresponding landing point
		\begin{equation}
		\mathcal{J}=I-n {n}^\intercal .
		\end{equation} That is, the component of $\z$ normal to $\Sigma$ becomes $0$ as the trajectory enters $\Sigma$ in backwards time.
	\end{remark}

	\begin{remark}\label{rem:uztable}
		Combining Theorem \ref{thm:main}, Remarks \ref{rem:Jland=I} and \ref{rem:iprc-liftoff}, we summarize the behavior of the solutions of the variational and adjoint equations $\mbu$ and $\z$ in limit cycles with sliding components:
		\begin{center}
			\begin{tabular}{llll}
				& Landing            & Sliding        & Liftoff \\
				\midrule
				$\mbu$ & $S=I-nn^\intercal$ & $\mbu^\perp=0$ & $S=I$ \\
				$\z$   & $\mathcal{J}=I$              & $\z^\perp=0$   & $\mathcal{J}=I-nn^\intercal$\\
			\end{tabular}
		\end{center}
		where $S$ is the regular saltation matrix, $\mathcal{J}$ is the time-reversed jump matrix and $^\perp$ denotes the normal component.
	\end{remark}
	
	\begin{remark}\label{rem:sala-jump}
		It follows directly from Remark~\ref{rem:uztable} that the relation between the saltation and jump matrices at a transversal boundary crossing point $J^\intercal S=I$ (see \eqref{eq:jump-salt-relation}) is no longer true at a landing or a liftoff point.
		Instead, the following condition holds
		\[
		\mathcal{J}^\intercal_p S_p = I-n n^\intercal
		\]
		where $\mathcal{J}_p$ and $S_p$ denote the time-reversed jump matrix and the regular saltation matrix at a landing or a liftoff point.
	\end{remark}
	
	\RED{\begin{remark}
			Assumption \ref{ass:assumption2} is necessary for the proof of part (c) in Theorem \ref{thm:main}.
			A stable limit cycle arising in a $C^r$-smooth vector field, for $r\ge 1$, will have $C^r$ isochrons \citep{Wiggins1994book,JosicShea-BrownMoehlis2006}.
			In \citep{park2018} and \citep{wilson2019} the authors assume differentiability of the phase function with respect to a basis of vectors spanning a switching surface.
			The assumption we require here is similarly plausible; it appears to hold at least for the model systems we have considered.
	\end{remark}}
	
	Next we illustrate the behavior of a limit cycle with sliding component via an analytically tractable planar model in \S\ref{sec:toy-model} \RED{and a stick-slip oscillator in \S\ref{sec:stick-slip}.
		In these two examples,}  we will see that a nonzero instantaneous jump discussed in Remark \ref{rem:iprc-liftoff} can occur in the normal component of $\z$ at the liftoff point, reflecting a ``kink'' or nonsmooth feature in the isochrons (\textit{cf.}~Figure \ref{fig:LC-2d-isochrons}).
	In \RED{both systems,} the discontinuity in the iPRC reflects a curve of nondifferentiability in the asymptotic phase function propagating backwards along a trajectory from the liftoff point to the interior of the domain (see Figure \ref{fig:LC-2d-isochrons}A).
	The presence of a discontinuous jump from zero to a nonzero normal component in $\z$ in forward time implies that numerical evaluation of the iPRC (presented in Appendix~\ref{sec:algorithm}) should be accomplished via backward integration along the limit cycle. \RED{In \S\ref{sec:two-mass}, we further apply the iPRC developed for the LCSC system to study the synchronization of two weakly coupled stick-slip oscillators, which together form a four-dimensional nonsmooth system with two sliding components.}  
	In addition to the iPRC, we also provide numerical algorithms for calculating the lTRC, the variational dynamics and the iSRC for a LCSC in a general nonsmooth system with hard boundaries, in Appendix~\ref{sec:algorithm}.

	\section{Applications to a planar limit cycle model with sliding components
	}\label{sec:toy-model}
	
	In this section, we apply our methods to a two dimensional, analytically tractable model that has a single interior domain with purely linear flow and hard boundary constraints that create a limit cycle with sliding components (LCSC).
	We find the surprising result that the isochrons exhibit a nonsmooth ``kink'' propagating into the interior of the domain from the locations of the liftoff points, \textit{i.e.}~the points where the limit cycle smoothly departs the boundary.
	In addition, we show that using local timing response curve analysis gives significantly greater accuracy of the shape response than using a single, global, phase response curve.
	
	MATLAB source code for simulating the model and reproducing the figures is available: \url{https://github.com/yangyang-wang/LC_in_square}.
	
	In the interior of the domain $[-1,1]\times [-1,1]$, we take the vector field of a simple spiral source to define the interior dynamics of the planar model
	
	\begin{equation}\label{eq:toy-model}
	\frac{d\textbf{x}}{dt}=F(\textbf{x})=\begin{bmatrix} \alpha x - y\\x+\alpha y \end{bmatrix} 
	\end{equation}
	where $\textbf{x}=\begin{bmatrix} x\\y\end{bmatrix}$ and $\alpha$ is the expansion rate of the source at the origin.
	The rotation rate is fixed at a constant value $1$.
	The Jacobian matrix $DF$ evaluated along the limit cycle solution in the interior of the domain is
	\begin{equation}\label{eq:toy-model-2}
	DF=\begin{bmatrix} \alpha & -1 \\ 1 & \alpha \end{bmatrix} .
	\end{equation}
	In what follows we will require $0<\alpha < 1$, so we have a weakly expanding source.
	For illustration, $\alpha=0.2$ provides a convenient value.
	Every trajectory starting from the interior, except the origin, will eventually collide with one of the walls at $x=\pm 1$ or $y=\pm 1$ (in time not exceeding $\frac 1{2\alpha}\ln(2/(x(0)^2+y(0)^2)$).
	As in \S\ref{sec:sliding}, we set the sliding vector field when the trajectory is traveling along the wall to be equal to the continuation of the interior vector field in the component parallel to the wall, while the normal component is set to zero (except where it is oriented into the domain interior).
	
\begin{figure}[!t]
		\centering
		\setlength{\unitlength}{1mm}
		\begin{picture}(145,60)(0,0)
		\put(-5,0){\includegraphics[width=6.5in]{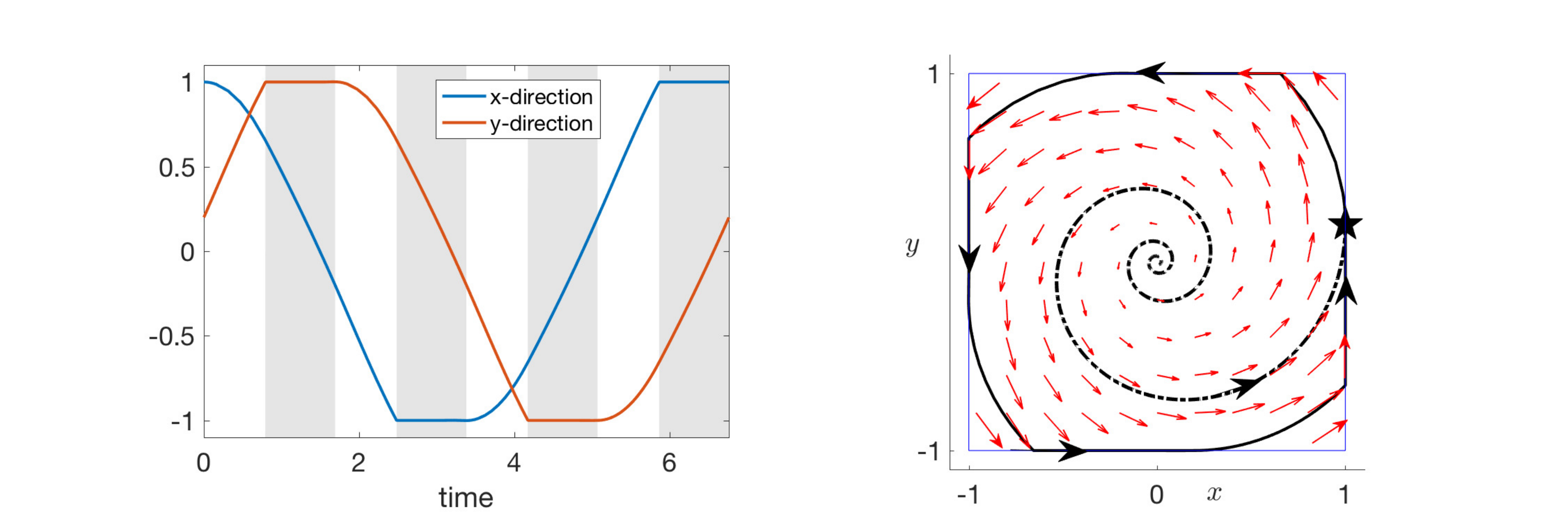}}
		\put(3,45){\bf{\large{(A)}}}
		\put(82,45){\bf{\large{(B)}}}
		\end{picture}  
		\caption{ Simulation result of the planar LCSC model with parameter $\alpha=0.2$.
			(A): Time series of the limit cycle $\gamma(t)$ generated by the planar LCSC model over one cycle with initial condition $\gamma(0)=[1,\alpha]^\intercal$.
			(B): Projection of $\gamma(t)$ onto $(x,y)$ phase space (solid black) and the osculating trajectory (dashed black), starting near the center that ends up running into the wall at the liftoff point $[1,\alpha]^\intercal$ (black star).
			Red arrows represent the vector field.
		} \label{fig:simu-toymodel} 
	\end{figure}

	The resulting vector fields of the planar LCSC model $F(\mbx)$ on the interior and along the walls are given in Table \ref{tab:natural}, and illustrated in Fig.~\ref{fig:simu-toymodel}B.
	
	\begin{table}[!htp]
		\centering
		\begin{tabular}{cccc}
			$x$ range & $y$ range & $dx/dt$ & $dy/dt$ \\ \hline
			$|x|<1$ & $|y|<1$ & $\alpha x-y$ & $x+\alpha y$\\ \hline
			$x=1$ & $-1\le y<\alpha$ & $0$ & $1+\alpha y$\\
			$x=1$ & $\alpha\le y < 1$ & $\alpha - y$ & $1+\alpha y$\\\hline
			$y=1$ & $1\ge x > - \alpha$ & $\alpha x-1$ & $0$\\
			$y=1$ & $-\alpha \ge x > -1$ & $\alpha x - 1$ & $x+\alpha$\\\hline
			$x=-1$ & $1\ge y>-\alpha$ & $0$ & $-1+\alpha y$\\
			$x=-1$ & $-\alpha\ge y > -1$ & $-\alpha - y$ & $-1+\alpha y$\\\hline
			$y=-1$ & $-1\le x <  \alpha$ & $\alpha x + 1$ & $0$\\
			$y=-1$ & $\alpha \le x < 1$ & $\alpha x + 1$ & $x-\alpha$\\\hline
		\end{tabular}
		\caption{Vector field of the planar LCSC model on the interior and along the boundaries.}
		\label{tab:natural} 
	\end{table}
	
	The trajectory will naturally lift off the wall and return to the interior when the normal component of the unconstrained vector field changes from outward to inward, i.e., $(F^{\text{interior}}\cdot n^{\text{wall}})|_{\text{wall}}=0$ (see \eqref{eq:cond-liftoff}).
	For instance, on the wall $x=1$ with a normal vector $n=[1,0]^\intercal$, we compute
	$$F^{\text{interior}}|_{\text{wall}}\cdot n^{\text{wall}}=(\alpha x-y)|_{\text{wall}}=\alpha-y=0.$$
	It follows that $y=\alpha$ defines the liftoff condition on the wall $x=1$.
	For this planar model, there are four lift-off points with coordinates $(1,\alpha), (-\alpha,1), (-1,-\alpha), (\alpha,-1)$ on the walls $x=1,\,y=1,\,x=-1, y=-1$, respectively.
	
	Denote the LCSC produced by the planar model by $\gamma(t)$, whose time series over $[0, T_0]$ is shown in Figure \ref{fig:simu-toymodel}A, where $T_0$ is the period.
	The projection of $\gamma(t)$ onto the $(x,y)$-plane is shown in the right panel, together with an osculating trajectory that starts near the center and ends up running into the wall $x=1$ at the lift-off point $(1,\alpha)$ (black star).

	Next we implement algorithms given in Appendix~\ref{sec:algorithm} to find the timing and shape responses of the LCSC to both instantaneous perturbations and sustained perturbation.
	We start by finding the iPRC for the LCSC to understand the timing response, and then solve the variational equation to find the linear shape response of the planar LCSC model to an instantaneous perturbation.
	Lastly, we compute the iSRC when the applied sustained perturbations are both uniform and nonuniform, to understand the shape response of the planar LCSC model to sustained perturbations.
	
	\subsection{Infinitesimal phase response analysis}\label{sec:iprc}
	
	In the case of weak coupling or small perturbations of a strongly stable limit cycle, a linearized analysis of the phase response curve -- the iPRC -- suffices to predict the behavior of the perturbed system.
	When trajectories slide along a hard boundary, however, the linearized analysis breaks down.
	For nonsmooth systems such as the LCSCs, the asymptotic phase function $\phi(\mbx)$ may itself be nonsmooth at certain locations, even when it remains well defined; its gradient (i.e., the iPRC) may therefore be discontinuous at those locations.
	Nevertheless, one may be able to derive a consistent first order approximation to the phase response curve notwithstanding that the directional derivative \eqref{eq:definition_of_iPRC} may not be well defined, as discussed in \S\ref{sec:nonsmooth-theory}.
	
	The dynamics of the planar LCSC model are smooth except for the discontinuities when crossing the switching boundaries, that is, entering or exiting the walls.
	The iPRC, $\z(t)$, will be continuous in the interior domain as well as in the interior of the four boundaries.
	As discussed in \S\ref{sec:sliding}, the discontinuity of iPRCs only occurs at the liftoff point.
	By Remark~\ref{rem:uztable}, the time-reversed jump matrix at a liftoff point, which takes the iPRC just after crossing the liftoff point to the iPRC just before crossing the liftoff point in backwards time, is given by 
	\begin{equation*}
	\mathcal{J}=\Matrix{0 & 0  \\
		0 & 1  }
	\end{equation*}
	when the trajectory leaves the walls $x=\pm 1$, and is given by
	\begin{equation*}
	\mathcal{J}=\Matrix{1 & 0  \\
		0 & 0  }
	\end{equation*}
	when the trajectory leaves the walls $y=\pm 1$.
	
	\begin{figure}[!t] 
		\centering
		\setlength{\unitlength}{1mm}
		\begin{picture}(145,60)(0,0)
		\put(0,0){\includegraphics[height=5.7cm]{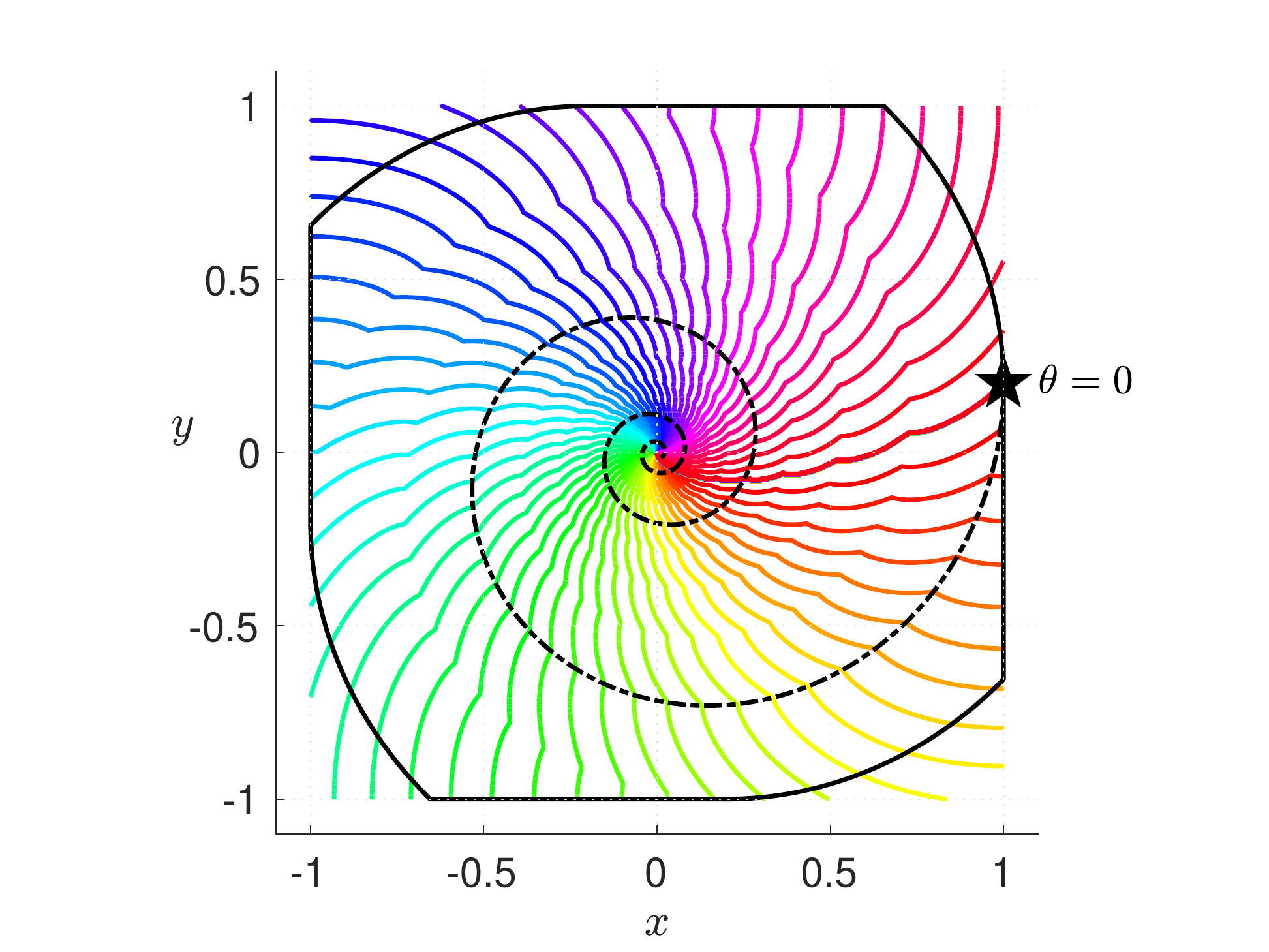}}
		\put(0,46){\bf{\large{(A)}}}
		\put(75,0){\includegraphics[height=5.5cm]{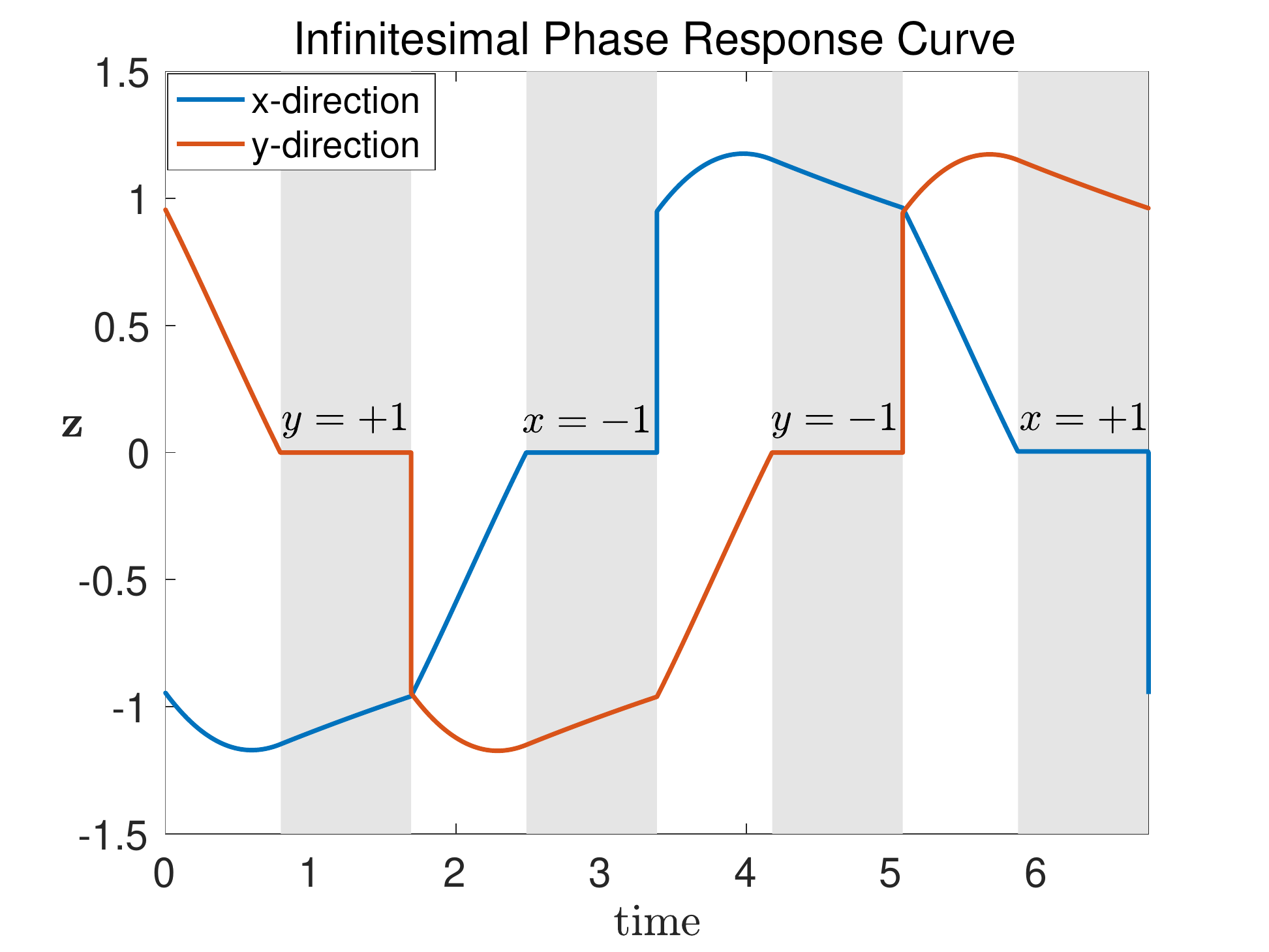}}
		\put(70,46){\bf{\large{(B)}}}
		\end{picture}  
		\caption{iPRCs for the planar LCSC model with parameter $\alpha=0.2$.
			(A): Trajectories and isochrons for the LCSC model.
			The solid black and dashed black curves are the same as in Figure \ref{fig:simu-toymodel}B.
			The colored scalloped curves are isochrons of the LCSC $\gamma(t)$ (black solid) corresponding to 50 evenly distributed phases $nT_0/50$, $n=1,\cdots, 50$.
			We define the phase at the liftoff point (black star) to be zero.
			(B): iPRCs for the planar LCSC model.
			The blue and red curves represent the iPRC for perturbations in the positive $x$ and $y$ directions, respectively.
			The intervals during which $\gamma(t)$ slides along a wall are indicated by the shaded regions.
			While $y=+1$, the iPRC vector is parallel to the wall ($\z_y\equiv 0$) and oriented opposite to the direction of flow ($\z_x < 0$).
			Similarly, on the remaining walls, the iPRC vector has zero normal component relative to the active constraint wall, and parallel component opposite the direction of motion. }
		\label{fig:LC-2d-isochrons} 
	\end{figure}

	Figure \ref{fig:LC-2d-isochrons}A shows the limit cycle (solid black curve), the osculating trajectory (dashed black curve) corresponding to the liftoff point (black star, $\theta=0$)  and the isochrons computed from a direct method, starting from a grid of initial conditions and tracking the phase of final locations (colored scalloped curves).
	There appears to be a ``kink" in the isochron function, propagating backwards in time along the trajectories that encounter the boundaries exactly at the liftoff points, such as the dashed curve.
	This apparent discontinuity in the gradient of the isochron function in the \emph{interior} of the domain exactly corresponds, at the boundary, with the point of discontinuity occurring in the iPRC along the limit cycle (cf.~Remark \ref{rem:iprc-liftoff}).
	According to Figure \ref{fig:LC-2d-isochrons}, the isochron curves are perpendicular to the sliding region of the wall at which the interior vector field is pointing outward.
	That is, the normal component of the iPRC when the trajectory slides along a wall is equal to 0.
	There is no jump in $\z$ when the trajectory enters the wall, but instead a discontinuous jump from zero to nonzero occurs in the normal component of $\z$ at the liftoff point.
	All of these observations are consistent with iPRC $\z$ (Figure \ref{fig:LC-2d-isochrons}, right) that is computed using \textbf{Algorithm for $\z$} in \S\ref{sec:algorithm-iprc} based on Theorem \ref{thm:main}.

	After the trajectory lifts off the east wall ($x=1$) at the point marked $\theta=0$ in Figure \ref{fig:LC-2d-isochrons}A (black star), a perturbation along the positive $x$-direction (resp., positive $y$-direction) causes a phase delay (resp., advance).
	While the timing sensitivity of the LCSC to small perturbations in the $x$-direction reaches a local maximum before reaching the next wall $y=1$, the phase advance caused by the $y$-direction perturbation decreases continuously to $0$ as the trajectory approaches $y=1$.
	As the trajectory is sliding along the wall ($y=1$), the positive $y$-direction perturbation that is normal to the wall has no effect on the LCSC and hence will not affect its phase.
	Moreover, we showed in Theorem \ref{thm:main} that a perturbation in the negative $y$-direction also has no effect on the phase, since the perturbed trajectory returns to the wall within time $O(\varepsilon)$, with a net  phase offset that is at most $O(\varepsilon^2)$, where $\varepsilon$ is the size of the perturbation.
	As the trajectory lifts off the wall, there is a discontinuous jump in $\z_y$, so that a negative $y$-direction perturbation applied immediately after the liftoff point leads to a phase advance.
	On the other hand, on the sliding region of the wall $y=1$, a perturbation along the positive $x$-direction, against the direction of the flow, results in a phase delay, which decreases in size as the phase increases, and becomes $0$ upon reaching the next wall, $x=-1$.
	The timing sensitivity of the LCSC to perturbations applied afterwards are similar to what are observed in the first quarter of the period due to the $\mathbf{Z}_4$-symmetry $\sigma (x,y) = (-y,x)$.
	
	The linear change in the oscillation period of the LCSC in response to a static perturbation can then be estimated by taking the integral of the iPRC multiplying the given perturbation, as shown in the last step in \textbf{Algorithm for $\z$} (\S~\ref{sec:algorithm-iprc}).
	As noted before, the change in period will be needed to solve \eqref{eq:src} for the iSRC to understand how this perturbation affects the shape of the LCSC.
	
	In this example, the interior vector field \eqref{eq:toy-model} is linear.
	Therefore its Jacobian is constant, and the iPRC may be obtained analytically \citep{park2018}.
	The resulting curves are indistinguishable from the numerically calculated curves shown in Fig.~\ref{fig:LC-2d-isochrons}B.
	
	\subsection{Variational analysis}
	
	Suppose a small instantaneous perturbation, applied at time $t=0$, leads to an initial displacement $\mbu(0)=\tilde{\gamma}(0)-\gamma(0)$, where $\gamma(0)=[1,\alpha]$ is the liftoff point  (black star in Figs.~\ref{fig:simu-toymodel}B and \ref{fig:LC-2d-isochrons}A) as in the previous section.
	We use the variational analysis to study how this perturbation evolves over time.
	
	Similar to the iPRC, $\mbu(t)$ will be continuous everywhere in the domain except when entering or exiting the walls.
	In contrast to the iPRC, $\mbu$ is continuous at all liftoff points, but exhibits discontinuous saltations when the trajectory enters a wall.
	According to Theorem \ref{thm:main}, the saltation matrix $S$, which takes $\mbu$ just before entering a wall to $\mbu$ just after entering the wall in forwards time, for the planar LCSC model is given by 
	\begin{equation*}
	S=\Matrix{0 & 0  \\
		0 & 1  }
	\end{equation*}
	when the trajectory enters the walls $x=\pm 1$, and is given by
	\begin{equation*}
	S=\Matrix{1 & 0  \\
		0 & 0  }
	\end{equation*}
	when the trajectory enters the walls $y=\pm 1$.

	Solutions to the variational equation of the planar LCSC model with the given initial condition $\mbu(0)$ can be computed using \textbf{Algorithm for $\mbu$} in \S\ref{sec:algorithm-u}.
	As discussed in Remark \ref{rem:fund-matrix}, an alternative way to find the displacement $\mbu(t)$ is to compute the fundamental solution matrix $\Phi(t, 0)$ by running \textbf{Algorithm for $\mbu$} twice and then to evaluate $\mbu(t) = \Phi(t,0) \mbu(0)$.
	The advantage of the latter approach is that once $\Phi(t,0)$ is obtained, it can be used to compute $\mbu(t)$ with any given initial value by evaluating a matrix multiplication instead of solving the variational equation.
	
\begin{figure}[!t] 
		\centering
		\includegraphics[width=6in]{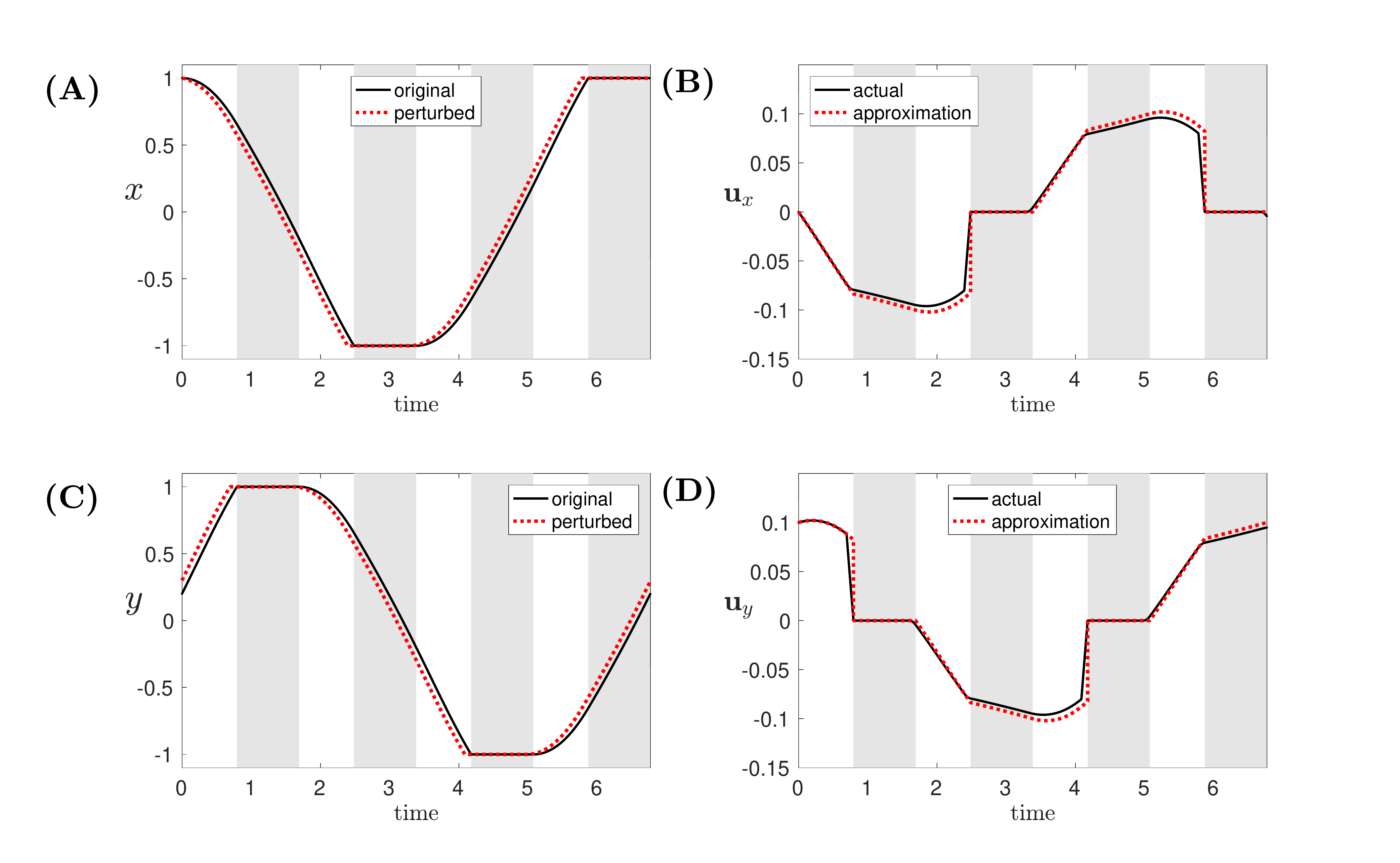}
		\caption{Linear shape response $\mbu(t)$ of the LCSC trajectory $\gamma(t)$ to an instantaneous perturbation applied at $\gamma(0)=[1,\alpha]$ when time $=0$ where $\alpha=0.2$.
			The initial displacement is $\mbu(0)=[0,\varepsilon]$ where $\varepsilon=0.1$.
			(A, C) Time series of $\gamma(t)$ (black solid) and $\tilde{\gamma}(t)$ (red dotted) with a perturbed initial condition $\tilde{\gamma}(0)=[1,\alpha+\varepsilon]$.
			(B, D) The difference between $\tilde{\gamma}(t)$ and $\gamma(t)$ obtained by direct calculation from the left panels (black) and the displacement solution $\mbu(t)$ obtained using the \textbf{Algorithm for $\mbu$} (red dotted).
			(A) and (B) show trajectories and the displacement along the $x$-direction, while (C) and (D) show trajectories and displacements along the $y$-direction.
			Shaded regions have the same meanings as in Figure \ref{fig:LC-2d-isochrons}.
		}
		\label{fig:var-toy} 
	\end{figure}

	Here, by taking $[1,0]$ and $[0, 1]$ as the initial conditions for $\mbu$ at the liftoff point A, we apply \textbf{Algorithm for $\mbu$} to compute the time evolution of the two columns for the fundamental matrix $\Phi(t,0)$.
	A simple calculation shows that the monodromy matrix $\Phi(T_0,0)$ has an eigenvalue $+1$, whose eigenvector $[0,1]$ is tangent to the limit cycle at the liftoff point, as expected (Remark \ref{rem:monodromy}).
	It follows that if the initial displacement at the liftoff point is along the limit cycle direction, then the displacement after a full period becomes the same as the initial one.
	To see this, we take the initial displacement $\mbu(0)=[0,\varepsilon]$ where $\varepsilon=0.1$ to be the tangent vector of the limit cycle at the liftoff point, and compute $\mbu(t)$, $x$ and $y$ components of which are shown in red dotted curves in Figure \ref{fig:var-toy}B,D.
	The saltations in $\mbu$ at time when the trajectory hits the walls can be clearly distinguished in the plot.
	Moreover, $\mbu(T_0)=\mbu(0)$ as we expect.
	
	To further validate the accuracy of $\mbu$, we solve and plot $\gamma(t)$ with $\gamma(0)=[1,\alpha]$ (black, Figure \ref{fig:var-toy}A,C) and the perturbed trajectory $\tilde{\gamma}(t)$ with $\tilde{\gamma}(0)=[1,\alpha]+\mbu(0)$ (red dotted, Figure \ref{fig:var-toy}A,C).
	The differences between the two trajectories along the $x$-direction and the $y$-direction are indicated by the black lines in Figure \ref{fig:var-toy}B and D, both showing good agreements with the approximated displacements computed from the variational equation, indicated by the red dotted lines in Figure \ref{fig:var-toy}B and D.
	Such an approximation becomes better as the perturbation size $\varepsilon$ gets smaller (simulation result not shown).

	Next, we study the effects of static perturbations on the timing and shape using the iPRC and iSRC.
	
	\subsection{Shape response analysis}\label{sec:planar-src}
	
	In this section, we illustrate how to compute the iSRC $\gamma_1$, the linear shape responses of the LCSC to small static perturbations.
	Recall that we use $\gamma_0(t)$ with period $T_0$ and $\gamma_{\varepsilon}(t)$ with period $T_{\varepsilon}$ to denote the original and the perturbed LCSC solutions.
	We write $\gamma_1$ for the linear shift in the limit cycle shape in response to the static perturbation as indicated by \eqref{eq:x_epsilon_of_tau}, which we also repeat here: 
	\[\gamma_\epsilon(\tau_\varepsilon(t))=\gamma_0(t)+\epsilon\gamma_1(t)+O(\epsilon^2),\] 
	where the time for the perturbed LCSC is rescaled to be $\tau_\varepsilon(t)$ to match the unperturbed time points.
	The iSRC $\gamma_1$ satisfies the nonhomogeneous variational equation \eqref{eq:src}.
	To solve this equation, an estimation of the timing scaling factor $\nu_1$, determined by the choice of time rescaling $\tau_\varepsilon(t)$, is needed.
	Here we consider two kinds of static perturbations on the planar LCSC model: global perturbation and piecewise perturbation.
	
	\paragraph{Global perturbation.}
	
	We apply a small static perturbation to the planar LCSC model by increasing the model parameter $\alpha$ by $\varepsilon$ globally: $\alpha\to \alpha+\varepsilon$.
	To compare the LCSCs before and after perturbation at corresponding time points, we rescale the perturbed trajectory uniformly in time so that $\tau_\varepsilon(t) = T_{\varepsilon}t/T_0$.
	It follows that $\nu_1 = T_1/T_0$, where the linear shift  $T_1:=\lim_{\varepsilon\to 0}(T_{\varepsilon}-T_0)/\varepsilon$
	can be estimated using the iPRC (see \eqref{eq:T1}).
	
\begin{figure}[!t]
		\begin{center}
			\setlength{\unitlength}{1mm}
			\begin{picture}(145,60)(0,0)
			\put(0,0){\includegraphics[height=5.5cm]{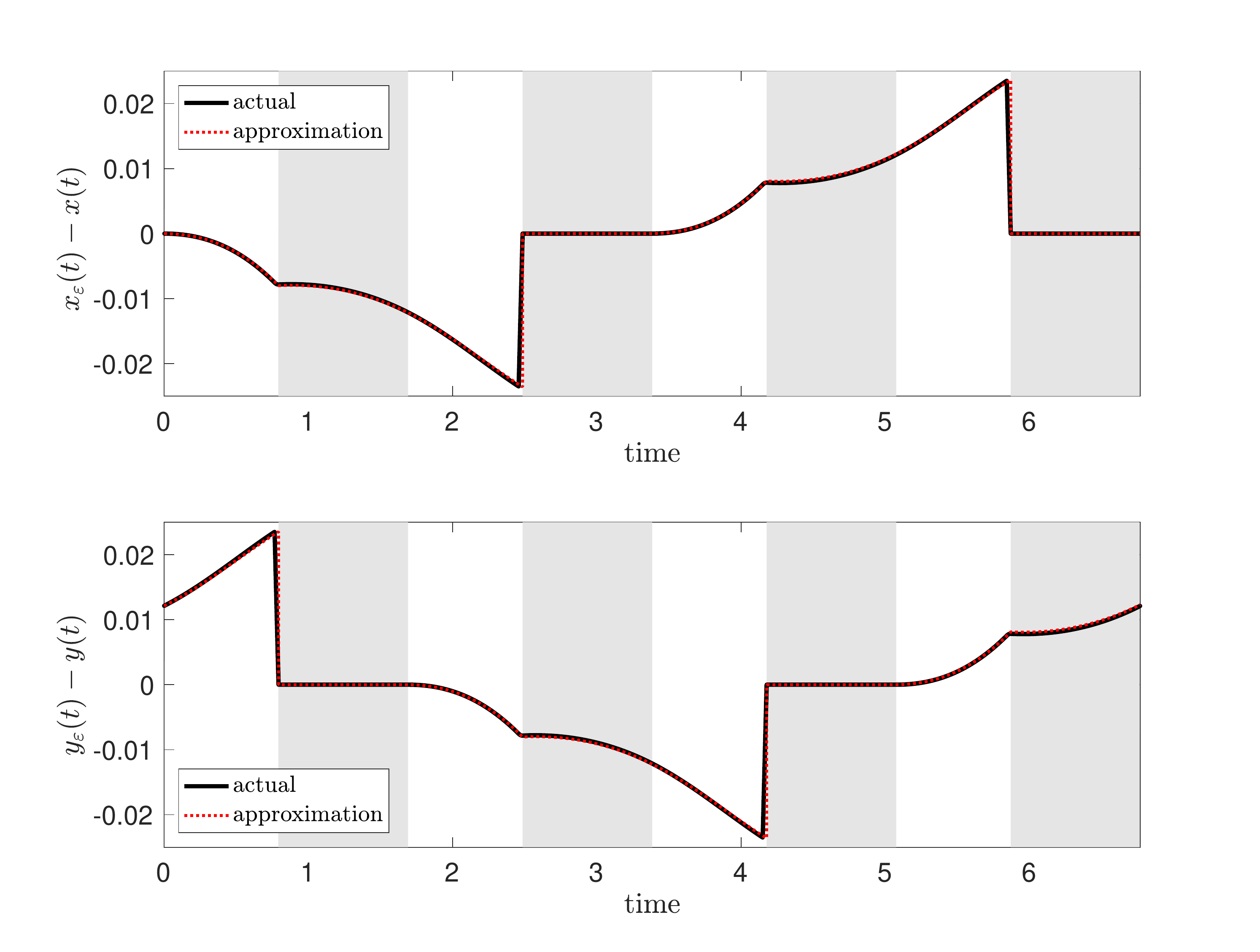}}
			\put(-8,46){\bf{\large{(A)}}}
			\put(75,0){\includegraphics[height=5.5cm]{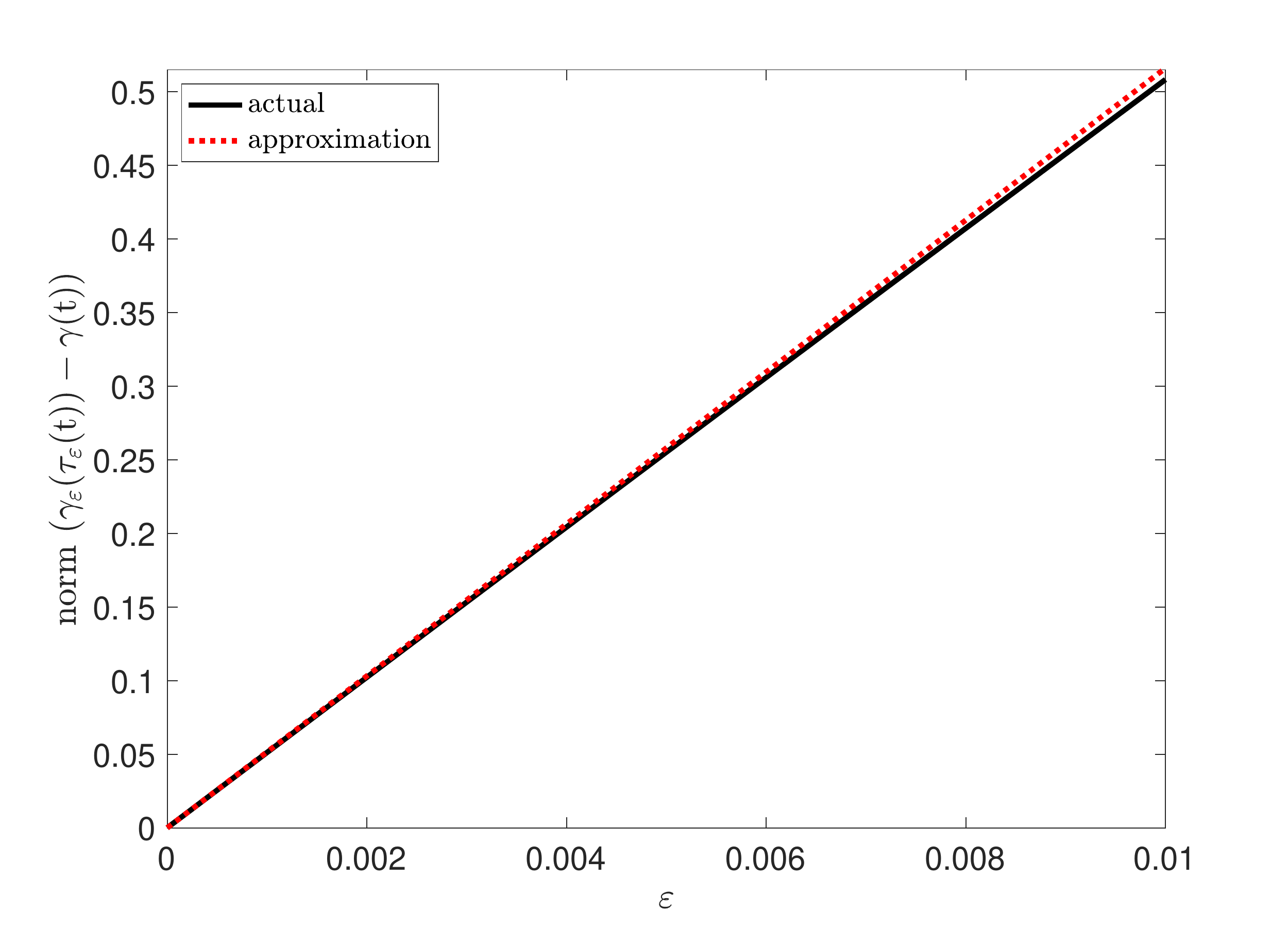}}
			\put(70,46){\bf{\large{(B)}}}
			\end{picture}
		\end{center}
		\caption{\label{fig:change-shape} 
			iSRC of the LCSC model to a small perturbation $\alpha\to \alpha+\epsilon$ with unperturbed parameter $\alpha=0.2$.
			(A) Time series of the difference between the perturbed and unperturbed solutions along the $x$-direction (top panel) and the $y$-direction (bottom panel) with $\epsilon=0.01$.
			The black curve denotes the numerical displacement computed by subtracting the unperturbed solution trajectory from the perturbed trajectory, after globally rescaling time.
			The red dashed curve denotes the product of $\varepsilon$ and the shape response curve solution.
			Shaded regions have the same meanings as in Figure \ref{fig:LC-2d-isochrons}.
			(B) The norm of the numerical difference (black) and the product of $\varepsilon$ and the iSRC (red dashed) grow linearly with respect to $\varepsilon$ with nearly identical slope, indicating that the iSRC is very good for approximating the numerical difference over a range of $\varepsilon$ and improves with smaller $\varepsilon$.
		}
	\end{figure}

	Using \textbf{Algorithm for $\gamma_1$ with uniform rescaling}, we numerically compute the iSRC $\gamma_1(t)$ for $\varepsilon=0.01$.
	The $x$ and $y$ components of $\varepsilon\gamma_1(t)$ are shown by the red curves in Figure \ref{fig:change-shape}A, both of which show good agreement with the numerical displacement $\gamma_\epsilon(\tau_\varepsilon(t))-\gamma_0(t)$ (black solid), as expected from our theory.
	
	For $\varepsilon$ over a range $[0, 0.01]$, we repeat the above procedure and compute the Euclidean norms of both the numerical displacement vector $\gamma_{\varepsilon}(\tau_\varepsilon(t))-\gamma(t)$ (Figure \ref{fig:change-shape}B, black solid) and the approximated displacement vector $\varepsilon\gamma_1(t)$ (Figure \ref{fig:change-shape}B, red dotted) over one cycle.
	From the plot, we can see that the iSRC with uniform rescaling of time gives a good first-order $\varepsilon$ approximation to the shape response of the planar LCSC model to a global static perturbation.
	
	\paragraph{Piecewise perturbation.}
	
	Uniform rescaling of time as used above is the simplest choice among many possible rescalings, and is shown to be adequate in the global perturbation case for computing an accurate iSRC.
	As discussed in \S\ref{sec:smooth-theory}, in certain cases we may instead need the technique of \emph{local} timing response curves (lTRCs) to obtain nonuniform choices of rescaling for greater accuracy.
	
\begin{figure}[!t]
		\begin{center}
			\setlength{\unitlength}{1mm}
			\begin{picture}(145,60)(0,0)
			\put(0,0){\includegraphics[height=5.5cm]{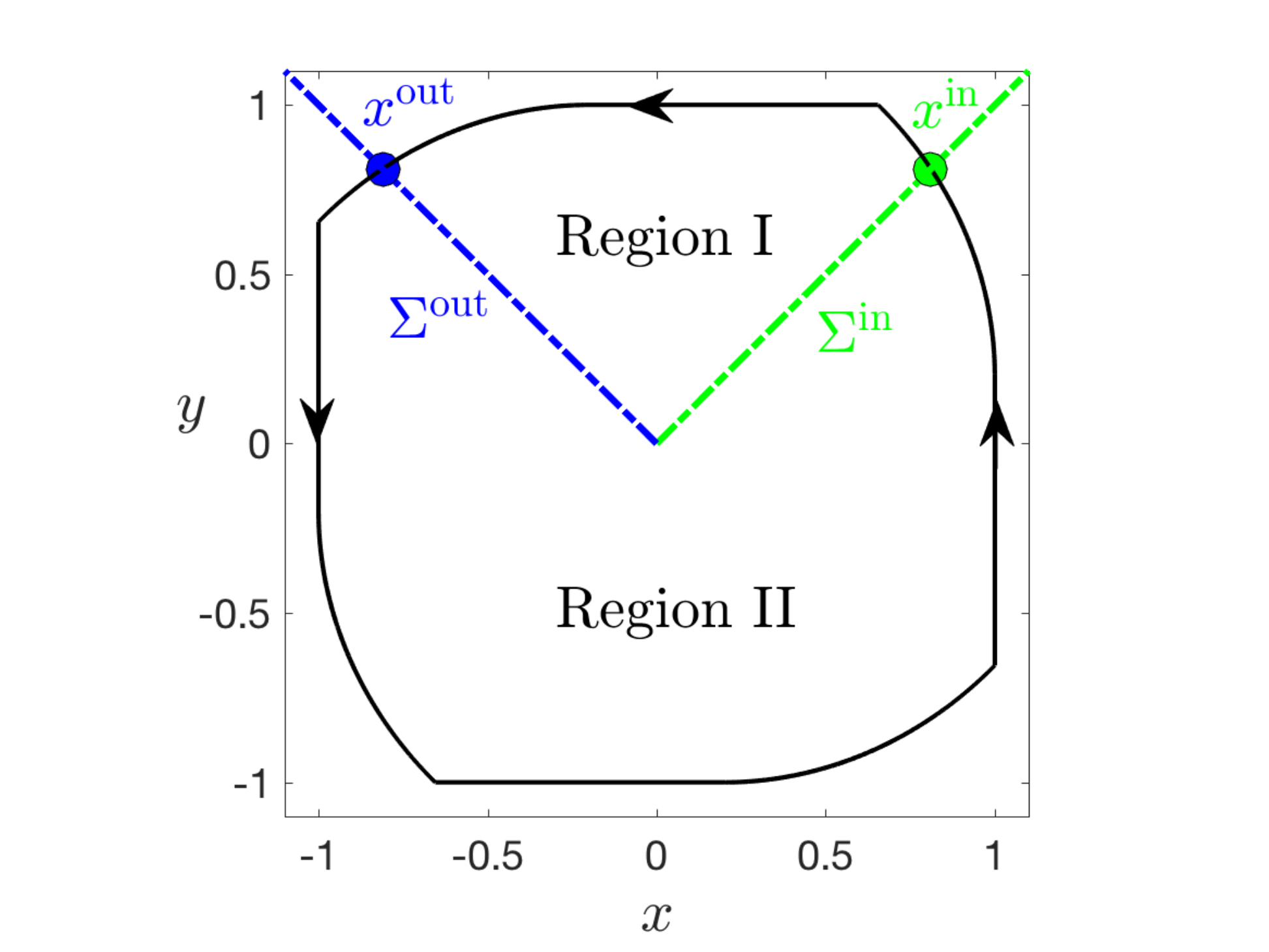}}
			\put(-5,46){\bf{\large{(A)}}}
			\put(75,0){\includegraphics[height=5.5cm]{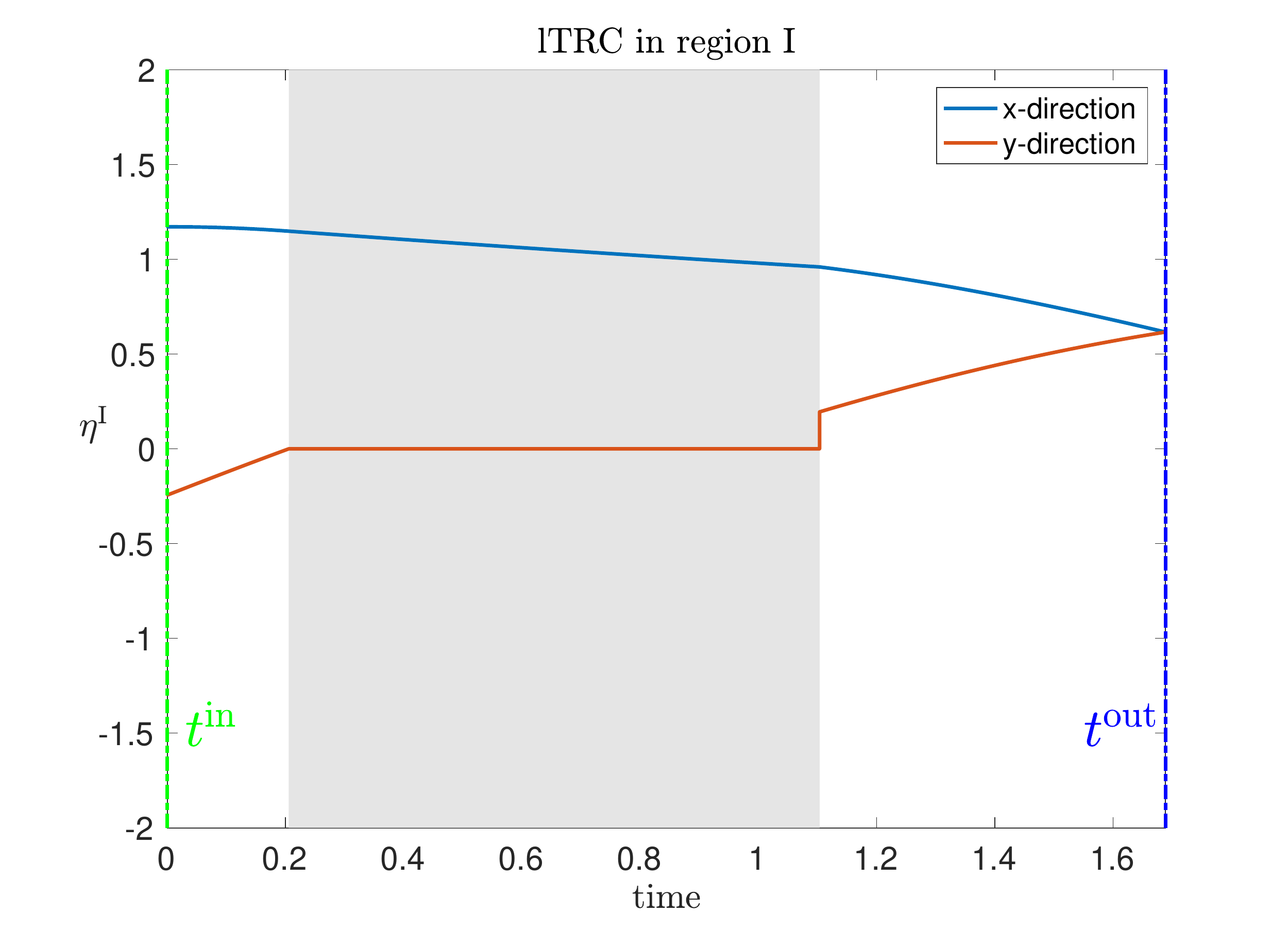}}
			\put(70,46){\bf{\large{(B)}}}
			\end{picture}
		\end{center}
		\caption{\label{fig:time-response} lTRC of the planar LCSC model under perturbation $(\alpha,\omega)\rightarrow (\alpha+\varepsilon,\omega-\varepsilon)$ over region I with unperturbed parameters $\alpha=0.2$ and $\omega=1$ held fixed in region II.
			(A) Projection of the limit cycle solution to the planar model with two new added switching surfaces $\Sigma^{\rm in}$ (green dashed line) and $\Sigma^{\rm out}$ (blue dashed line) onto its phase plane.
			(B) Time series of the lTRC $\eta^{\rm I}$ from $t^{\rm in}$ (the time of entry into region I at $\mbx^{\rm in}$) to $t^{\rm out}$ (the time of exiting region I at $\mbx^{\rm out}$).
			A discontinuous jump occurs when the trajectory exits the wall $y=1$ indicated by the right boundary of the shaded region, which has the same meaning as in Figure \ref{fig:LC-2d-isochrons}.
		}
	\end{figure}

	As an illustration, we add two local timing surfaces $\Sigma^{\rm in}$ and $\Sigma^{\rm out}$ to the planar LCSC model (see Figure \ref{fig:time-response}A).
	We denote the subdomain above $\Sigma^{\rm in}$ and $\Sigma^{\rm out}$ by region I ($\RR^{\rm I}$) and denote the remaining subdomain by region II ($\RR^{\rm II}$).
	Moreover, we introduce a new parameter $\omega$, the rotation rate of the source at the origin, that has previously been fixed at $1$, and rewrite the interior dynamics of the planar LCSC model as
	\begin{equation}\label{eq:toy-model-omega}
	\frac{d\textbf{x}}{dt}=F(\textbf{x})=\begin{bmatrix} \alpha x -\omega y\\\omega  x+\alpha y \end{bmatrix}.
	\end{equation}
	The vector fields on a given wall are obtained by replacing the coefficient of $y$ in $dx/dt$ (in Table \ref{tab:natural}) by $-\omega$ and replacing the coefficient of $x$ in $dy/dt$ (in Table \ref{tab:natural}) by $\omega$ on that wall.
	
\begin{figure}[!t]
		\begin{center}
			\setlength{\unitlength}{1mm}
			\begin{picture}(145,60)(0,0)
			\put(0,0){\includegraphics[height=5.5cm]{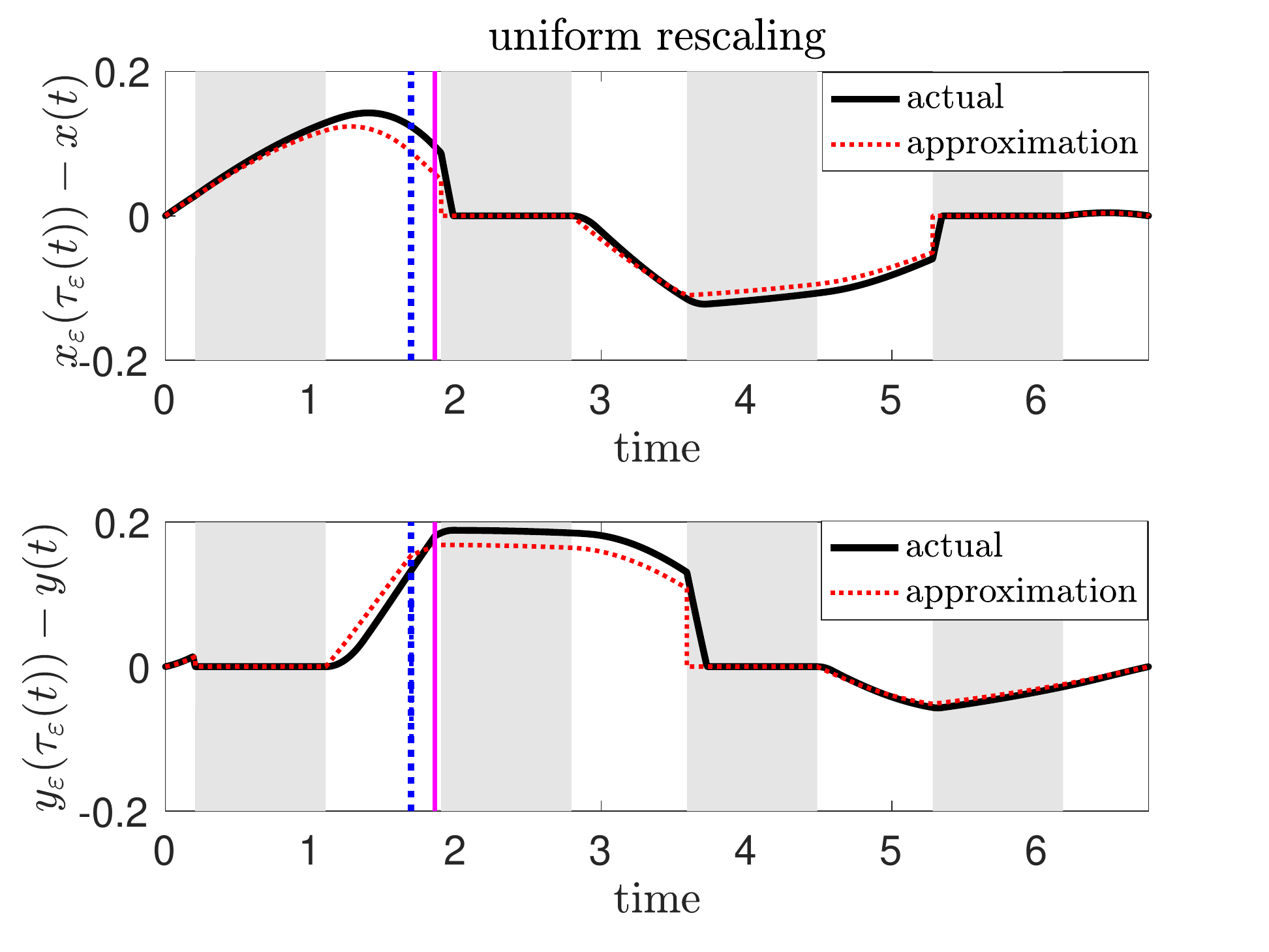}}
			\put(-8,46){\bf{\large{(A)}}}
			\put(75,0){\includegraphics[height=5.5cm]{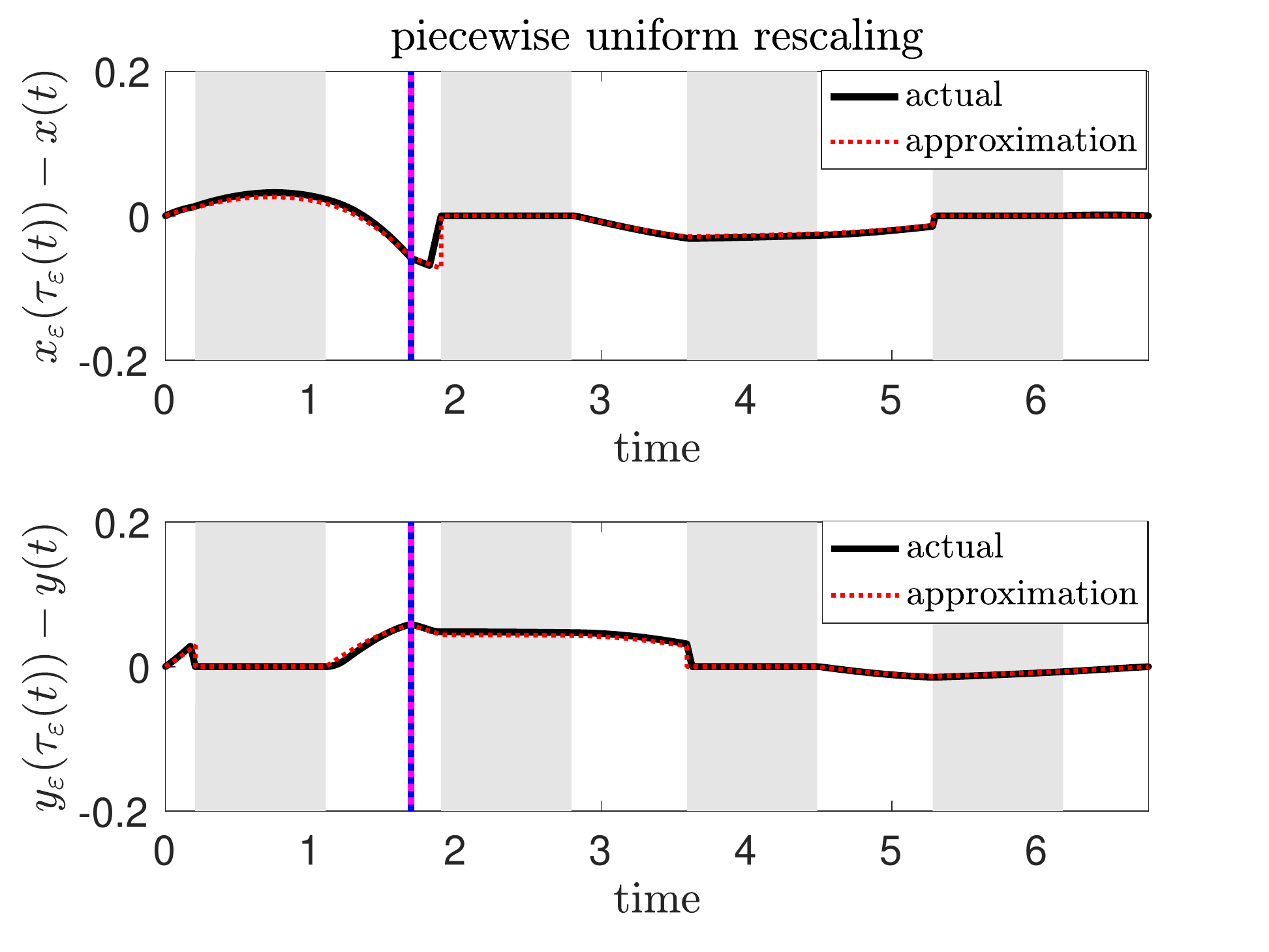}}
			\put(68,46){\bf{\large{(B)}}}
			\end{picture}
		\end{center}
		\caption{\label{fig:shape-response}
			A small perturbation is applied to the planar model over region I in which $(\alpha,\omega)\to(\alpha+\varepsilon,\omega-\varepsilon)$ with unperturbed parameters $\alpha=0.2$, $\omega=1$ and perturbation  $\varepsilon=0.1$.
			Time series of the difference between the perturbed and unperturbed solutions along the $x$-direction (top panel) and the $y$-direction (lower panel) using (A) the global rescaling factor and (B) two different rescaling factors within regions I and II.
			The vertical blue dashed line denotes the exit time of the unperturbed trajectory from Region I, while the vertical magenta solid line denotes the exit time of the rescaled perturbed trajectory from Region I.
			Other color codings of lines are the same as in Figure \ref{fig:change-shape}A.
			Shaded regions have the same meanings as in Figure \ref{fig:LC-2d-isochrons}.
		}
	\end{figure}

\begin{figure}[!t]
		\begin{center}
			\setlength{\unitlength}{1mm}
			\begin{picture}(145,60)(0,0)
			\put(0,0){\includegraphics[height=5.5cm]{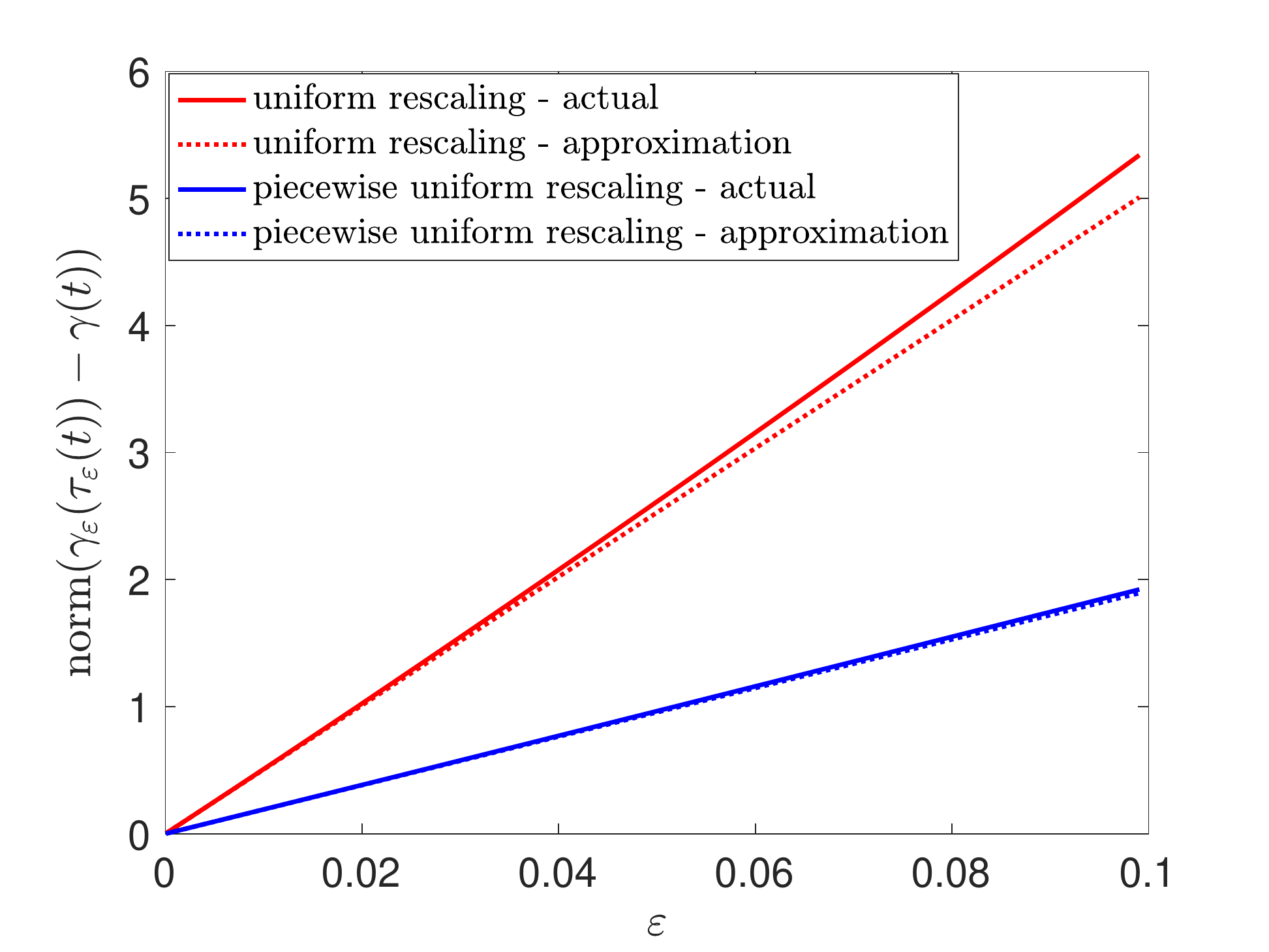}}
			\put(-5,46){\bf{\large{(A)}}}
			\put(75,0){\includegraphics[height=5.5cm]{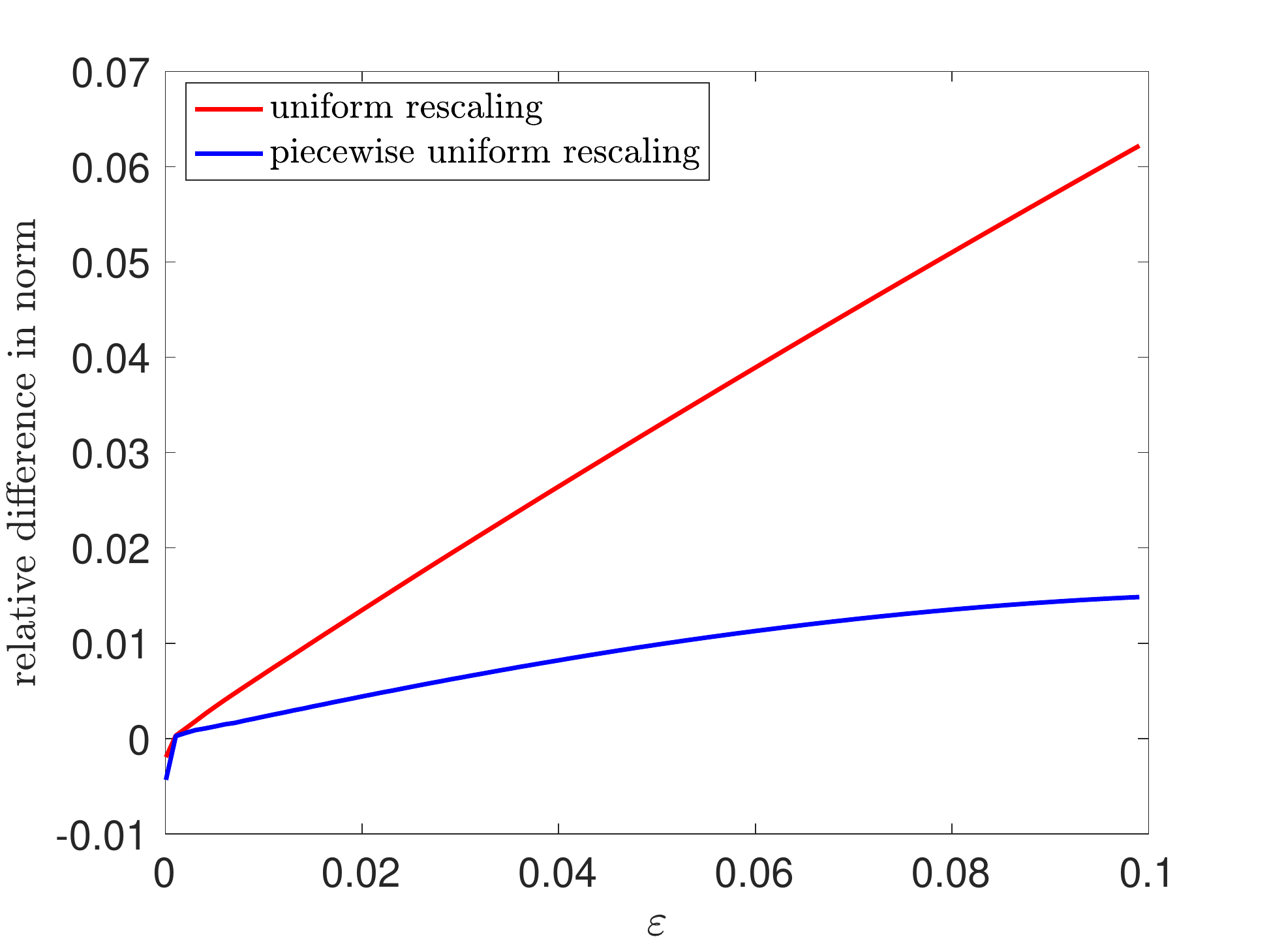}}
			\put(70,46){\bf{\large{(B)}}}
			\end{picture}  
		\end{center}
		\caption{\label{fig:norm-shape-response} A small parametric perturbation is applied to the planar model over region I in which $(\alpha,\omega)\to(\alpha+\varepsilon,\omega-\varepsilon)$ with unperturbed parameters $\alpha=0.2$, $\omega=1$.
			(A): Values of the Euclidean norm of $(\gamma_\varepsilon(\tau_\varepsilon(t))-\gamma(t))$ computed numerically (solid curve) versus those computed from the iSRC (dashed curve), as $\varepsilon$ varies.
			The norms grow approximately linearly with respect to $\varepsilon$. The approximation obtained by the iSRC when using piecewise uniform rescaling (blue) is closer to the actual simulation than using the uniform rescaling (red). (B): The relative difference between the actual and approximated norms with a uniform rescaling (red) is larger than that when piecewise uniform rescaling is used (blue).
			The difference between the two curves expands as $\varepsilon$ increases.}
	\end{figure}

	We apply a static piecewise perturbation to the system by letting $(\alpha,\omega)\rightarrow (\alpha+\varepsilon,\omega-\varepsilon)$ 
	over region I but not region II.
	Such a piecewise constant perturbation affects both the expansion and rotation rates of the source in region $\rm I$, and hence will lead to different timing sensitivities of $\gamma(t)$ in the two regions.
	It is therefore natural to use piecewise uniform rescaling when computing the shape response curve as opposed to using a uniform rescaling.
	In the following, we first compute the lTRC (see Figure \ref{fig:time-response}) and use it to estimate the two time rescaling factors for $\RR^{\rm I}$ and $\RR^{\rm II}$, which are denoted by $\nu_1^{\rm I}$ and $\nu_1^{\rm II}$, respectively.
	We then show the iSRC computed using the piecewise uniform rescaling factors provides a more accurate representation of the shape response to the piecewise static perturbation than using a uniform rescaling (see Figure \ref{fig:shape-response} and \ref{fig:norm-shape-response}).
	
	Although the lTRC $\eta$ is defined throughout the domain, estimating the effect of the perturbation localized to region I only requires evaluating the lTRC in this region.
	Figure \ref{fig:time-response}B shows the time series of $\eta^{\rm I}$ for the planar LCSC model in region I, obtained by numerically integrating the adjoint equation (\ref{eq:ltrc}) backward in time with the initial condition of $\eta^{\rm I}$ given by its value at the exit point of region I denoted by $\mbx^{\rm out}$ (see \textbf{Algorithm for $\eta^j$}).
	
	Similar to the iPRC, the $y$ component of the lTRC $\eta^\text{I}$ shown by the red curve in Figure \ref{fig:time-response}B is zero along the wall $y=1$, and the only discontinuous jump of $\eta^{\rm I}$ occurs at the liftoff point.
	Note that $\eta^{\rm I}$ is defined as the gradient of the time remaining in $\RR^{\rm I}$ until exiting through $\Sigma^{\rm out}$.
	If the $x$ or $y$ component of $\eta^{\rm I}$ is positive then the perturbation along the positive $x$-direction or $y$-direction increases the time remaining in $\RR^{\rm I}$, and the exit from $\RR^{\rm I}$ will occur later.
	On the other hand, if the $x$ or $y$ component $\eta^{\rm I}$ is negative then the perturbation along the positive $x$-direction or $y$-direction decreases the time remaining in $\RR^{\rm I}$, and the exit from $\RR^{\rm I}$ will occur sooner.
	The relative shift in time spent in $\RR^{\rm I}$ caused by a static perturbation can therefore be estimated using the lTRC (see \eqref{eq:local-time-shift}) as illustrated in the last step of \textbf{Algorithm for $\eta^j$}.
	Note that the first term in \eqref{eq:local-time-shift} implies that the timing change in a region generically depends on the shape change at the corresponding entry point, leading to the possibility of bidirectional coupling between timing and shape changes.
	However, in this planar system, a  perturbed trajectory with $\varepsilon\ll 1$ will converge  back to the original trajectory, within region II (where the perturbation is absent), in finite time.
	Under these circumstances, there is no shift between the perturbed and unperturbed trajectories in the entry location to region I.
	Hence, in this case, the local timing shift does not depend on the shape change.
	
	Let $T_{0}^{\rm I}$ denote the time spent in region I, and let $T_{0}^{\rm II}=T_{0}-T_{0}^{\rm I}$ denote the time spent in region II (recall $T_{0}$ is the total period).
	The linear shift in $T_0^{\rm I}$, denoted by $T_1^{\rm I}$, can be estimated using the lTRC $\eta^{\rm I}$ as discussed above.
	By definition the two time rescaling factors required to compute the iSRC are given by $\nu_1^{\rm I}=\frac{T_1^{\rm I}}{T_0^{\rm I}}$ and $\nu_1^{\rm II}=\frac{T_1-T_1^{\rm I}}{T_0^{\rm II}}$ where the global relative change in period, $T_1$, can be estimated using the iPRC as discussed before.
	With $\nu_1^{\rm I}$ and $\nu_1^{\rm II}$ known, we take $\mbx^{\rm in}$, the coordinate of the entry point into $\RR^{\rm I}$, as the initial condition for $\gamma(t)$ and apply \textbf{Algorithm for $\gamma_1$ with piecewise uniform rescaling} to compute the iSRC $\gamma_1$ for $\varepsilon=0.1$.
	The $x$ and $y$ components of $\varepsilon \gamma_1$ are shown by the red dashed curves in Figure \ref{fig:shape-response}B, both of which show good agreement with the numerical displacement $\gamma_{\varepsilon}(\tau_\varepsilon(t)) - \gamma(t)$ (black solid curves).
	Here the rescaling $\tau_\varepsilon(t)$ is piecewise uniform:
	\begin{eqnarray}
	\label{eq:piecewise-uniform-time-rescaling}
	\tau_\varepsilon(t)=\left\{
	\Eqn{
		t^{\rm in}+T_{\varepsilon}^{\rm I}(t-t^{\rm in})/T_0^{\rm I}, & \gamma(t) \in \RR^{\rm I}\\\\
		t^{\rm in}+T_{\varepsilon}^{\rm I}+T_{\varepsilon}^{\rm II}(t-t^{\rm out})/T_0^{\rm II}, & \gamma(t) \in \RR^{\rm II}\\
	}
	\right.
	\end{eqnarray}
	where $T_{\varepsilon}^i$ denotes the time $\gamma_{\varepsilon}$ spends in $\RR^i$ with $i\in\{\rm I, II\}$.
	It follows that the exit time of the trajectory from region I before (Figure \ref{fig:shape-response}B, vertical blue line) and after (Figure \ref{fig:shape-response}B, vertical magenta line) perturbation are the same.

	As a comparison, for $\varepsilon=0.1$, we also compute the iSRC and the numerical displacement using the uniform rescaling of time as we did in the global perturbation case (see Figure \ref{fig:shape-response}A).
	The difference between the vertical blue and magenta lines (the time when the unperturbed and perturbed trajectories leave region I) indicates region I and region II have different timing sensitivities.
	As expected, the resulting $\varepsilon \gamma_1$ no longer shows good agreement with the numerical displacement obtained from subtracting the unperturbed solution from the rescaled perturbed solution.

	Piecewise uniform rescaling, on the other hand, leads to a more accurate iSRC for the LCSC model \eqref{eq:toy-model-omega} than uniform rescaling, when the  LCSC $\gamma(t)$ experiences distinct timing sensitivities for $\varepsilon=0.1$.
	Fig.~\ref{fig:shape-response} contrasts the accuracy of the linearized shape response using global (A) versus local (B) timing response curves, for $\varepsilon=0.1$.
	We also show the same conclusion holds for other $\varepsilon$ values.
	To this end, for $\varepsilon$ over a range of $[0,0.1]$ we repeat the above procedure and compute the Euclidean norms of both the numerical displacement vector and the displacement vector approximated by the iSRC, as illustrated in Figure \ref{fig:norm-shape-response}A.
	The numerical and approximated norm curves using the uniform rescaling are shown in red solid and red dotted lines, while the numerical and approximated norm curves using the piecewise uniform rescaling are shown in blue solid and blue dotted lines.
	Unsurprisingly, the norms of the displacements between the perturbed and original trajectory grow approximately linearly with respect to $\varepsilon$, and the displacement norms with piecewise uniform rescaling are smaller than that with uniform rescaling.
	The fact that the difference between the lines in red is much bigger than the difference between the lines in blue suggests that the piecewise uniform rescaling gives a  more accurate iSRC than using the uniform rescaling for $\varepsilon\in [0,0.1]$, as we expect.
	This heightened accuracy is further demonstrated in Figure \ref{fig:norm-shape-response}B, where the relative difference between the numerical and approximated norms with uniform rescaling (red curve) is significantly larger than the relative difference when using piecewise uniform rescaling (blue curve).

	\RED{\section{Applications to stick-slip oscillators}\label{sec:stick-slip}
		
		In this section, we return to our motivating example --- the stick-slip oscillator --- presented in \S\ref{sec:intro}, which is a piecewise smooth system that exhibits a LCSC (see Figure \ref{fig:ss-motivating}). 
		In \S\ref{sec:one-mass}, we apply the analysis of lTRC and iSRC to compute the shape response of the single stick-slip oscillator to parametric perturbations as done in \S\ref{sec:toy-model}. In \S\ref{sec:two-mass} we use the phase reduction method to show that two weakly coupled identical stick-slip oscillators exhibit anti-phase synchronization. 
		
		\subsection{Shape response of a stick-slip oscillator to parametric perturbations}\label{sec:one-mass}
		
		We consider a stick-slip system consisting of a block of mass $m$ supported by a moving belt with constant velocity $u$. The block is connected to a fixed support by an elastic spring with stiffness $k$ and a linear dashpot with damping coefficient $c$. The surface between the block and the belt is rough so that the belt exerts a static friction force on the block which sticks to the belt during the \textit{stick phase} until the elastic force due to the spring and the damping force generated by the dashpot build up to exceed the maximum static friction force. At this point the \textit{slip phase} begins and the slipping motion is described by the following equation \citep{GB99,DL11},  
		\begin{equation}\label{eq:one-mass}
		mx''+cx'+kx = f(x'-u),
		\end{equation}
		where $x(t)$ is the displacement of the oscillator from the position at which the spring assumes its natural length; $'$ and $''$ indicate first and second order differentiation with respect to $t$.
		The \textit{kinetic friction force} of the block for $x'< u$ is given by 
		\begin{equation}\label{eq:stick-slip-force}
		f(x'-u) = \frac{1-\delta}{1-\gamma(x'-u)}+\delta+\eta(x'-u)^2 
		\end{equation}
		where $\delta\in[0, 1],\gamma>0$, and $\eta>0$. When $x'>u$, the kinetic friction force is $\frac{-(1-\delta)}{1+\gamma(x'-u)}-\delta-\eta(x'-u)^2$. The maximum static friction force is assumed to be $f_s=1$ when there is zero relative velocity (i.e., $x'=u$) and therefore the friction force is continuous in the stick-slip transition, whereas in general such a transition may be characterized by a finite jump in the friction force.  The slipping motion governed by \eqref{eq:one-mass} will continue to the point where there is no relative motion between the block and belt so that $x'=u$, and the elastic and damping forces are balanced by the static friction force so that $x''=0$.  
		
		\begin{figure}[!t]
			\centering
			\includegraphics[width=0.8\textwidth]{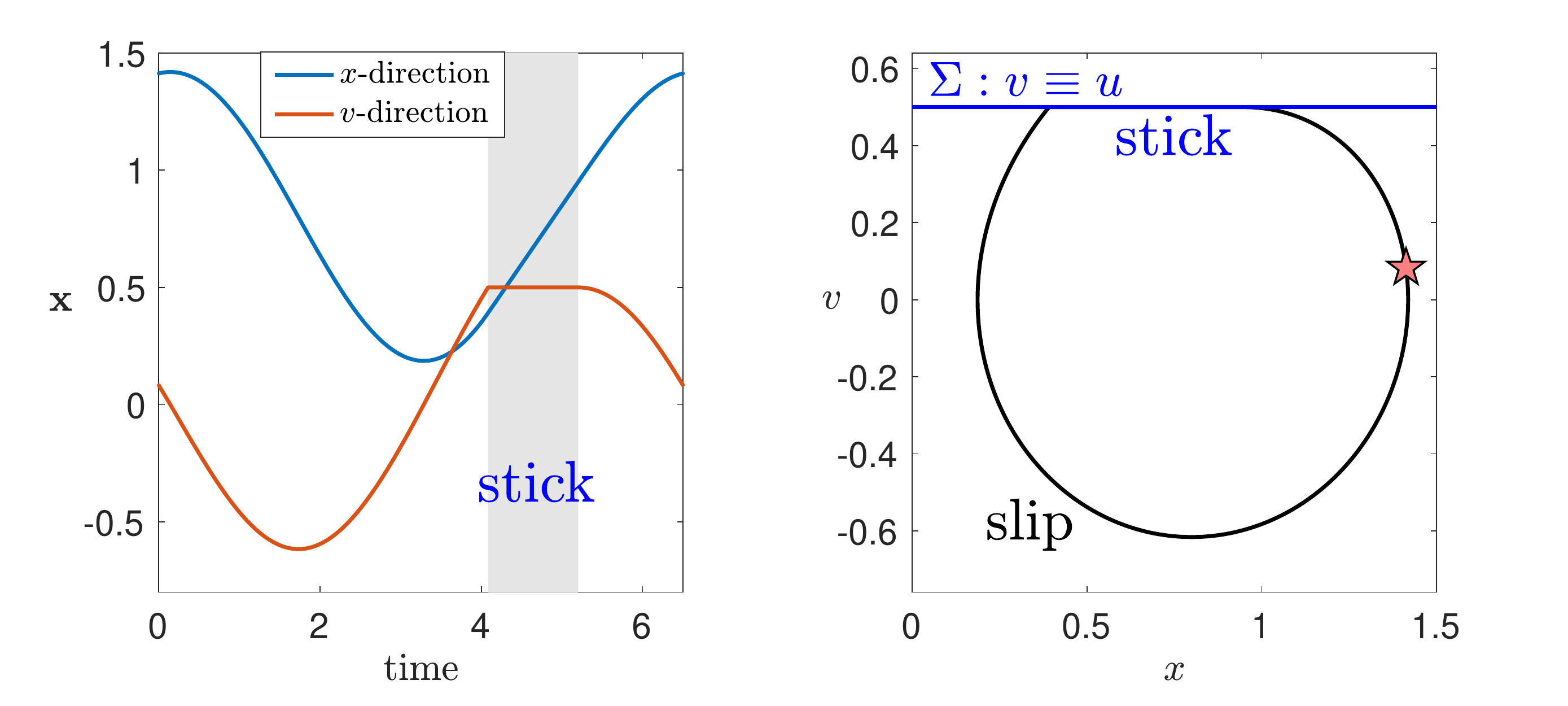}
			\caption{\RED{Simulation result of the stick-slip system \eqref{eq:one-mass-ps} when $m=1,k=1,c=0.1,\delta=0.5,\gamma=1,\eta=0.001,u=0.5$. Left: Time series of the LCSC with initial condition $[1.4127, 0.0829]^\intercal$. The shaded region indicates the stick phase. Right: Projection of the LCSC from the left onto $(x,v)$ phase space (solid black) and the hard boundary $\Sigma: v\equiv u$ (solid blue). The red star symbol denotes the starting point.}}
			\label{fig:stick-slip-simu}
		\end{figure}
		
		The continuous repetition of sticking and slipping motions with appropriate parameters can lead to a stick-slip limit cycle with sliding components (LCSC) within one zone $x'\leq u$. See Figure \ref{fig:stick-slip-simu}. The sliding region $\RR^{\rm slide}$ when the system is constrained to one less degree of freedom by a hard boundary, corresponds physically to the stick phase of the stick-slip system, when the block is captured by the moving belt and is carried along with the velocity $x'\equiv u$ until it escapes. It follows that the hard boundary is $\Sigma=\{x'=u\}$ with a unit normal vector $n=[0, 1]^\intercal$.  Letting $v = x'$ and $\mbx =[x,v]^\intercal$, we rewrite the stick-slip system in the form of \eqref{eq:1zoneFP} where 
		\begin{equation}\label{eq:one-mass-ps}
		\frac{d\mbx}{dt}=F(\mbx):=\begin{cases}
		\begin{array}{lclr}
		F^{\rm interior}(\mbx) &=& \Matrixc{v \\ -\frac{k}{m}x-\frac{c}{m}v+\frac{f(v-u)}{m}},& \mbx\in \RR^{\rm interior} \\
		F^{\rm slide}(\mbx) &=& \Matrixc{u \\ 0}, &\mbx\in \RR^{\rm slide}\subset \Sigma.
		\end{array}
		\end{cases}
		\end{equation}
		where $F^{\rm interior}$ is the vector field during the slip phase and $F^{\rm slide}$ is the vector field during the stick phase. It follows from Definition \ref{def:sliding} and Definition \ref{def:liftoff} that  $\RR^{\rm slide} = \{\mbx\in \Sigma\, |\, x < \frac{1-cu}{k} \}$ and the flow exits the sliding region (i.e., the stick phase ends) when $x=\frac{1-cu}{k}$. In other words, the stick phase terminates when the maximum static friction force $f_s=1$ acting on the block is exceeded by the other two forces; that is, $kx+cx'=f_s$ where $x'=u$ during the stick phase. The interior domain is therefore $\RR^{\rm interior} = \{\mbx\in \Sigma\, |\, x \geq \frac{1-cu}{k} \} \cup \{\mbx\in \mathbb{R}^2\,|\, v< u\}.$ 
		
		\begin{remark}\rm
			In contrast to our setup in which the stick-slip system \eqref{eq:one-mass-ps} is restricted to the domain $v\leq u$, \citep{GB99,DL11} used \emph{Filippov's Convex Method} to write the stick-slip system in the full space $\mathbb{R}^2$
			\begin{eqnarray}\label{eq:stick-slip-full}
			\frac{d\textbf{x}}{dt}=F(\mbx):=\left\{
			\Eqn{
				F^{\rm I}(\mbx),  & v<u \\
				F^{\Sigma}(\mbx), & v=u \\
				F^{\rm II}(\mbx),  & v>u
			}
			\right.
			\end{eqnarray}
			where $F^{\Sigma}(\mbx)=(1-\alpha(\mbx))F^{\rm I}+\alpha(\mbx) F^{\rm II}$ and $\alpha(\mbx)=\frac{n^\intercal F^{\rm I}}{n^\intercal (F^{\rm I}-F^{\rm II})}$. Nonetheless, the two systems \eqref{eq:one-mass-ps} and \eqref{eq:stick-slip-full} exhibit the same LCSC that has the same iPRC and iSRC in response to parametric perturbations. This is because $F^{\rm I}$ is identical to $F^{\rm interior}$ in \eqref{eq:one-mass-ps} and, moreover, it follows from direct calculation that $F^{\Sigma}$ is identical to $F^{\rm slide}$ in \eqref{eq:one-mass-ps}.  Since we are only interested in understanding the timing and shape response of the stick-slip LCSC to small parametric perturbations, it is enough for us to work with \eqref{eq:one-mass-ps} that is defined in the domain $v\leq u$ where the LCSC exists. 
		\end{remark}
		\begin{figure}[!t]
			\centering
			\includegraphics[width=0.6\textwidth]{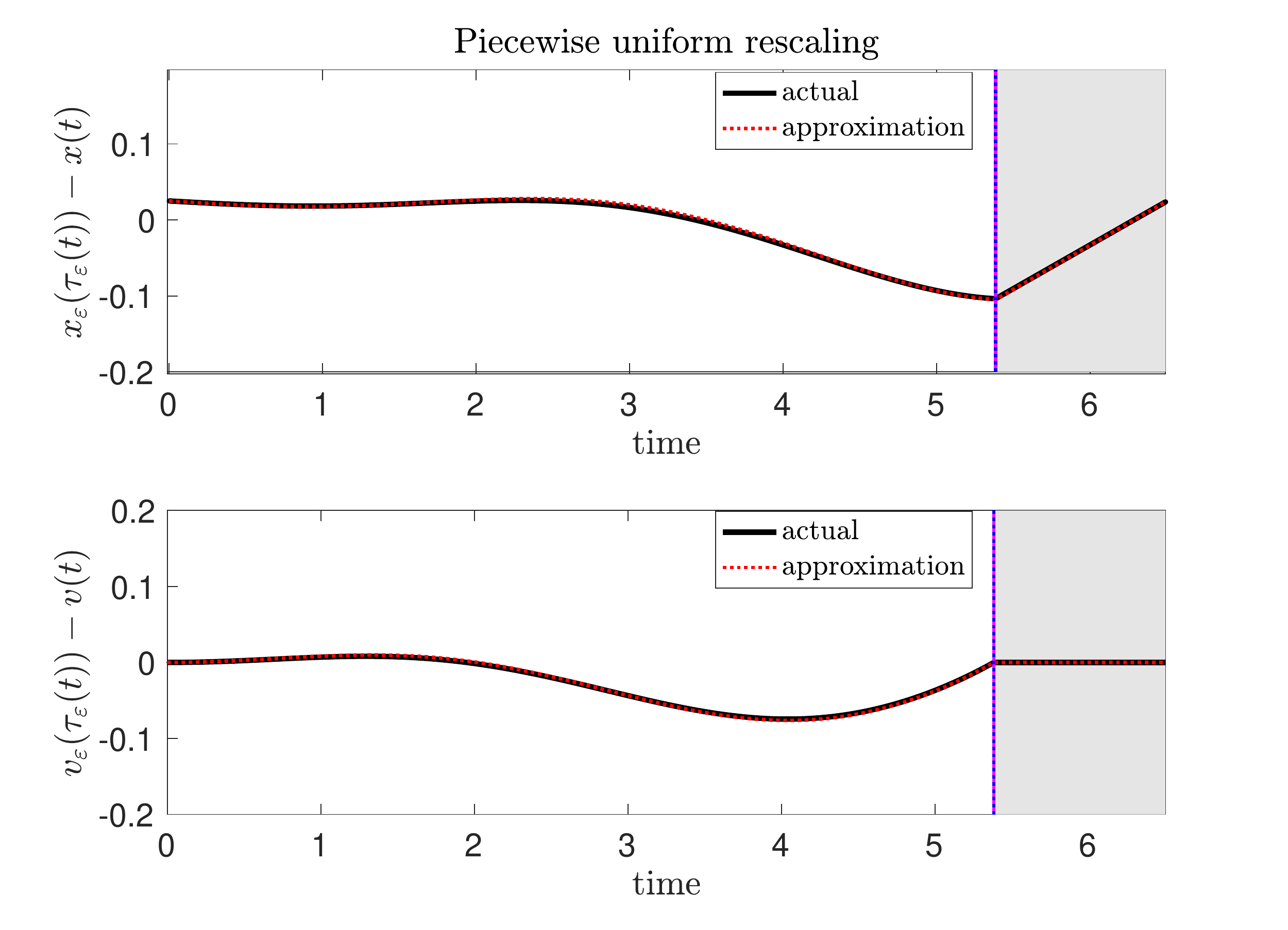}
			\caption{\RED{Shape response of the stick-slip system \eqref{eq:one-mass-ps} to a small parametric perturbation of the damping coefficient $c\to c+\varepsilon$ with $\varepsilon=-0.05$ and other unperturbed parameters the same as in Figure \ref{fig:stick-slip-simu}. Time series of the displacement between the perturbed and unperturbed solutions along the $x$-direction (resp., the $v$-direction) are shown in the top panel (resp., the bottom panel). The black curve denotes the numerical displacement computed by subtracting the unperturbed solution trajectory from the time-rescaled perturbed trajectory. The red dashed curve denotes the product of the iSRC solution and the perturbation size $\varepsilon$. The white and shaded regions indicate the slip and stick phases, respectively, and the blue/magenta dashed line indicates the transition.}}
			\label{fig:stick-slip-isrc}
		\end{figure}
		
		When $m=1,k=1,c=0.1,\delta=0.5,\gamma=1,\eta=0.001,u=0.5$, the system \eqref{eq:one-mass-ps} exhibits a LCSC, denoted by $\gamma(t)$, that slides along the hard boundary $\Sigma$ (see Figure \ref{fig:stick-slip-simu}). 
		In the following, we study the shape response of $\gamma(t)$ to a small parametric perturbation using the iSRC and lTRC, as discussed in \S\ref{sec:intro}.  
		Similar to the constructed planar model \eqref{eq:toy-model-omega} with two local timing surfaces, variations of the model parameters in \eqref{eq:one-mass-ps} (e.g., the damping coefficient $c$) only affect the vector field in the interior domain (i.e., during the slip phase), which naturally leads to a piecewise perturbation on the system. The lTRC is therefore needed to compute rescaling factors in different regions (interior/boundary) or phases (slip/stick), which are required to compute the iSRC. 
		The iSRC obtained by using piecewise uniform rescaling factors is shown by red dashed lines in Figure \ref{fig:stick-slip-isrc}, agreeing with the actual displacements (black solid line) between rescaled perturbed and unperturbed trajectories.

		\subsection{Anti-synchrony of two weakly coupled stick-slip oscillators}\label{sec:two-mass}
		
		The theory of weakly coupled oscillators has been used to predict the synchronization properties in networks of oscillators in smooth systems \citep{SL2012,Park2017} and nonsmooth systems with transversal crossing boundaries \citep{park2018}. In this section, we consider a nonsmooth system with hard boundaries composed of two blocks connected by a spring on a moving belt and apply the weakly coupled oscillator theory to study the synchrony between the two coupled stick-slip oscillators. Our results predict a non-intuitive result, namely, that the anti-synchrony  solution is stable whereas in-phase synchronization is unstable for an identical pair of stick-slip oscillators. Moreover, the anti-synchrony solution has an extremely slow rate of convergence. 
		
		A slip phase of the coupled stick-slip oscillators begins according to the following equations of motion \citep{Galvanetto2001},
		\begin{equation}\label{eq:ss-coupled}
		\begin{array}{lcl}
		m_1 x_1'' &=&-k_1x_1-k_{3}(x_1-x_2)+f_{1}(x_1'-u), \\
		m_2 x_2'' &=&-k_2x_2-k_{3}(x_2-x_1)+f_{2}(x_2'-u),
		\end{array}
		\end{equation}
		where $x_i$ is the displacement, $m_i$ is the mass, $k_i$ for $i=1,2$ is the stiffness of the spring connecting block $i$ to the fixed support and $k_3$ is the stiffness of the coupling spring, $f_{i}(x_i'-u)$ is the kinetic friction force of the $i$-th block and $u$ is the velocity of the moving belt. The damping coefficient of the spring is assumed to be $0$ for simplicity.    
		Each mass can undergo a stick phase, which leads to two hard boundaries:
		\[
		\Sigma^1=\{x_1'=u\}, \quad \Sigma^2=\{x_2'=u\}.
		\]
		
		Letting $X_1=[x_1,v_1]^\intercal$ and $X_2 = [x_2,v_2]^\intercal$ where $v_1 = x_1'$ and $v_2 = x_2'$ and assuming $m_1=m_2=m$, $k_1=k_2=k$ and $f_1=f_2=f$ so the two uncoupled oscillators are identical, we rewrite the coupled stick-slip systems in the following form 
		\begin{equation}\label{eq:ss-twomass-vec}
		\begin{array}{lcl}
		X_1'&=& F(X_1)+k_3 G(X_2,X_1)\\
		X_2' &=& F(X_2)+k_3G(X_1,X_2)
		\end{array}
		\end{equation}
		where 
		\begin{eqnarray*}
			F(X_i):=\left\{
			\Eqn{
				\Matrixc{v_i \\ -\frac{k}{m}x_i +\frac{f(v_i-u)}{m}}, & X_i \in \RR^{{\rm interior_i}}\\
				\Matrixc{u \\ 0}, & X_i \in\RR^{{\rm slide_i}}\subset \Sigma^i}
			\right.
		\end{eqnarray*}
		\begin{eqnarray*}
			G(X_j,X_i):=\left\{
			\Eqn{
				\Matrixc{0\\-(x_i-x_j)/m}, & X_i \in \RR^{{\rm interior_i}}\\
				\Matrixc{0 \\ 0}, & X_i \in\RR^{{\rm slide_i}}\subset \Sigma^i}
			\right.
		\end{eqnarray*}
		and the kinetic force $f$ is given by \eqref{eq:stick-slip-force}.
		The sliding regions ($\RR^{{\rm slide_1}}, \RR^{{\rm slide_2}}$), confined to $\Sigma^1$ and $\Sigma^2$, and the interior domains ($\RR^{\rm interior_1}$, $\RR^{\rm interior_2}$) can be found according to Definition \ref{def:sliding} and Definition \ref{def:interior}, as we did for the one-mass stick-slip system \eqref{eq:one-mass-ps}. 
		
		When there is no coupling with $k_3=0$, the two oscillators are identical and exhibit a $T_0$-periodic LCSC solution denoted as $\gamma(t)$. As discussed before, the periodic solutions of an asymptotically stable oscillator can be represented by a single variable phase model (see \eqref{eq:phase-reduction}) and the coupled oscillators can then be converted to the following phase model
		\begin{equation}\label{eq:ss-phase-1}
		\begin{array}{lcl}
		\theta_1'&=& 1+k_3 \z(\theta_1)\cdot G(X_2,X_1)\\
		\theta_2' &=& 1+k_3\z(\theta_2)\cdot G(X_1,X_2)
		\end{array}
		\end{equation}
		where $\theta_1\in [0, T_0]$ and $\theta_2\in [0, T_0]$ are the phases of the two uncoupled oscillators. $\z$ is the iPRC curve for the uncoupled stick-slip oscillator. If the coupling strength $k_3$ is sufficiently small, the uncoupled oscillator is almost identical to the periodic solutions $X_j(t)$. The system \eqref{eq:ss-phase-1} can then be approximated by
		\begin{equation}
		\begin{array}{lcl}
		\theta_1'&=& 1+k_3 \z(\theta_1)\cdot G(\gamma(\theta_2),\gamma(\theta_1))\\
		\theta_2' &=& 1+k_3\z(\theta_2)\cdot G(\gamma(\theta_1),\gamma(\theta_2)).
		\end{array}
		\end{equation}
		Averaging the right hand sides over one cycle $[0, T_0]$ and defining $\psi=\theta_2-\theta_1$, we can obtain the following scalar equation of the relative phase:
		\begin{equation}\label{eq:H}
		\begin{array}{lcl}
		\psi' &=& k_3 (H(-\psi) - H(\psi)) \equiv k_3 \HH(\psi),
		\end{array}
		\end{equation}
		where 
		\begin{equation}
		H(\psi) = \frac{1}{T_0}\int_{0}^{T_0} \z(t)\cdot G(\gamma(t+\psi),\gamma(t)) \,d t.
		\end{equation}
		The autonomous and scalar ODE \eqref{eq:H} can then be used to predict the synchronization rates and stability using a standard stability analysis on the phase line. 
		
		\begin{figure}[!htp]
			\centering
			\includegraphics[width=0.45\textwidth]{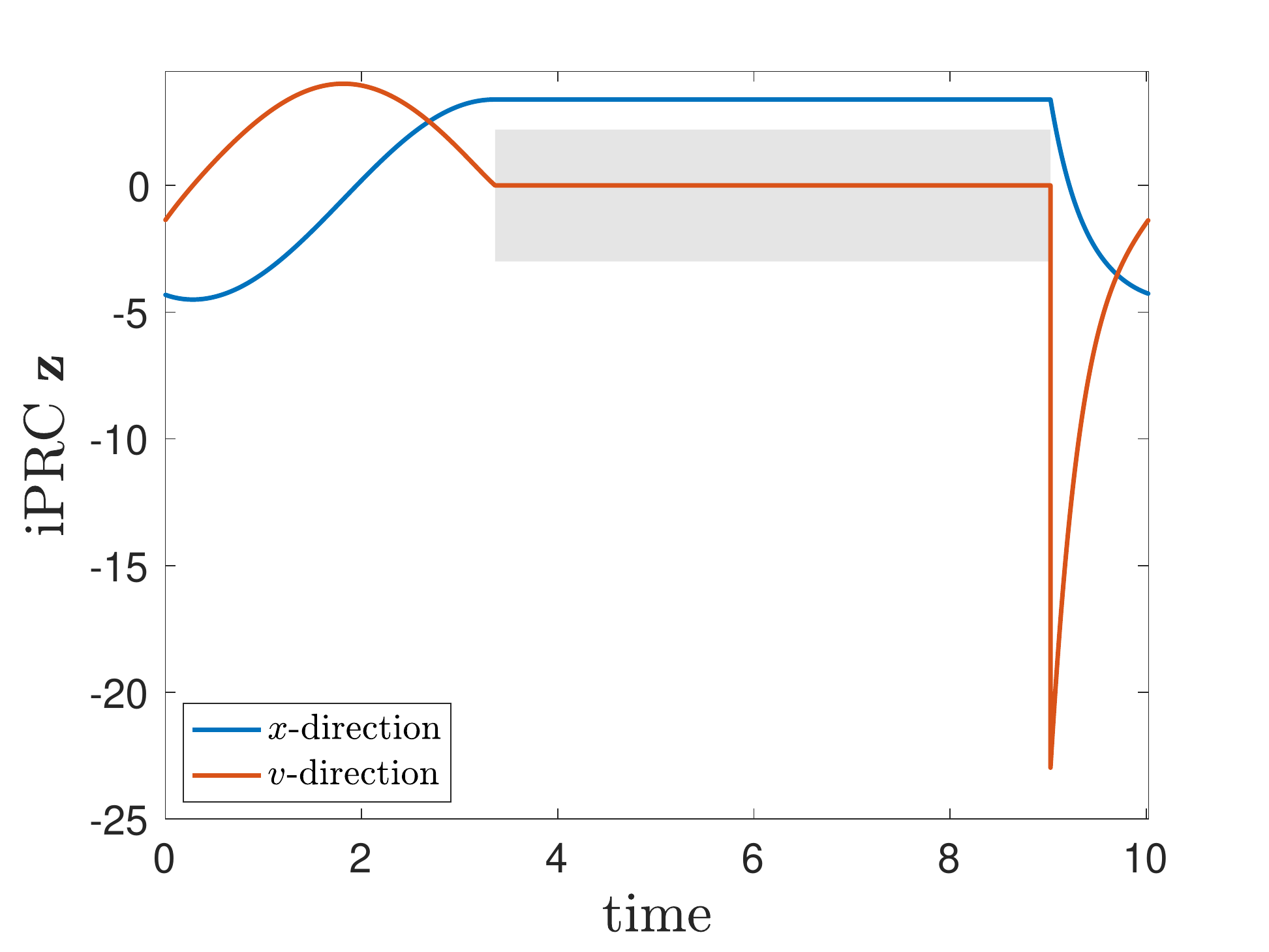}
			\includegraphics[width=0.45\textwidth]{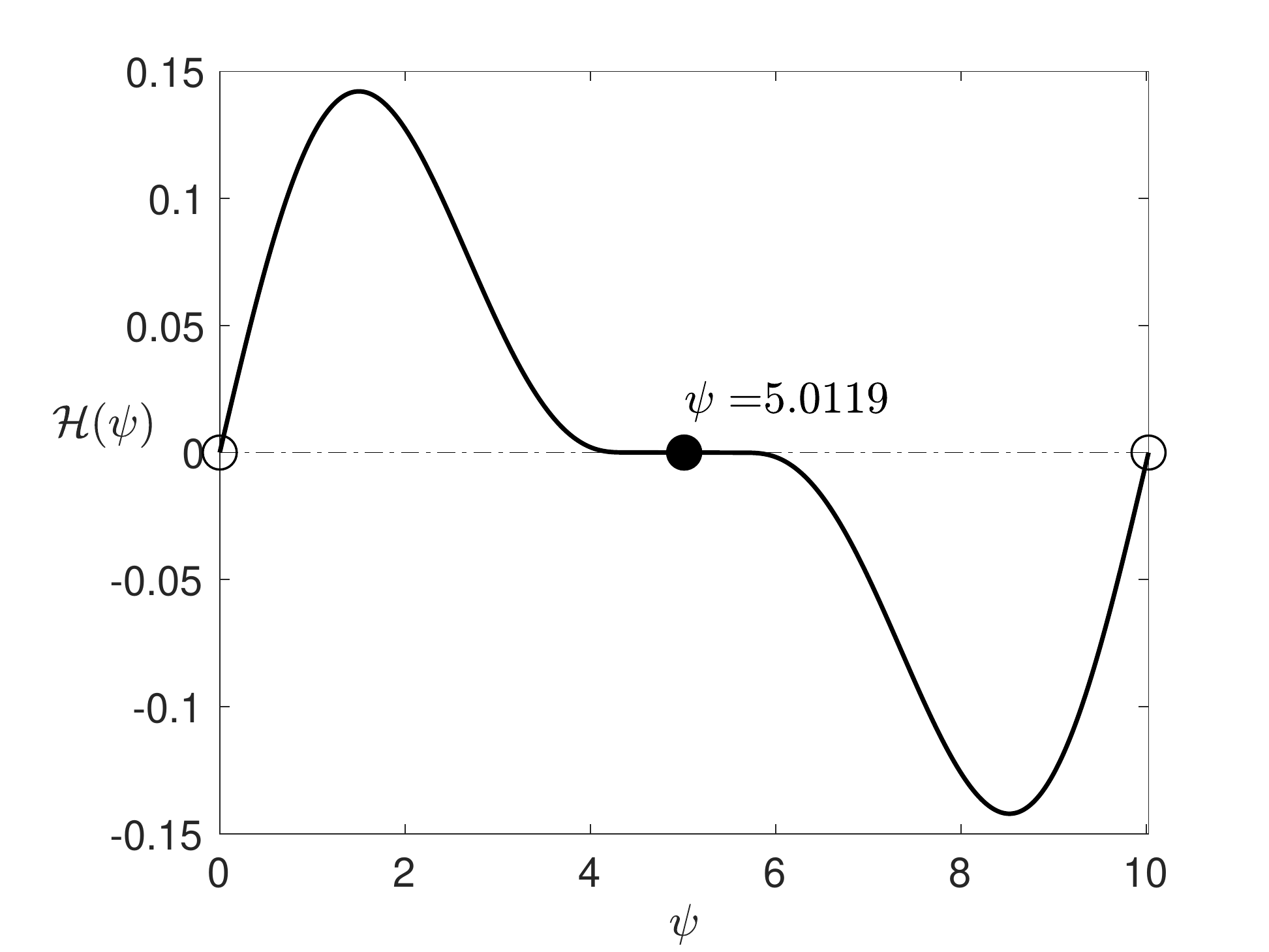}
			\caption{\RED{iPRC for the uncoupled stick-slip system and the right hand side $\HH$ of \eqref{eq:H} when $m=1, k=1, \delta=0, \gamma=3, \eta=0$ and $u=0.295$.  Left: iPRC for the unperturbed oscillator. The shaded region indicates the stick phase. Right: Stability analysis of the right hand side function $\HH$ of \eqref{eq:H}. Filled (resp., hollow) circle denotes an asymptotically stable (resp., unstable) phase locked solution. }}
			\label{fig:coupled-H_psi}
		\end{figure}
		
		When $m=1, k=1, \delta=0, \gamma=3, \eta=0$, and $u=0.295$, the iPRC for the uncoupled stick-slip oscillator is shown in Figure \ref{fig:coupled-H_psi}, left panel. As discussed before, the discontinuous jump in the iPRC occurs at the liftoff boundary, that is, when the stick phase ends. The right panel shows the right hand side function $\HH$ of \eqref{eq:H}. On the phase line, a filled circle at $\psi=T_0/2$ (resp., open circle at $\psi=0, T_0$) corresponds to an asymptotically stable (resp., unstable) phase locked solution. Hence, the phase model predicts that coupled stick-slip oscillators will diverge from an in-phase synchrony and asymptotically converge to an anti-phase synchrony. This prediction is supported by the numerical simulation in Figure \ref{fig:coupled-simu}, left panel, showing that oscillations beginning in-phase eventually converge to an anti-phase synchrony solution.

		The accuracy of the synchronization rates is demonstrated in Figure \ref{fig:coupled-simu}, right panel, where the red dashed curve is the solution to \eqref{eq:H} and the black curve is the numerical phase difference in the full model \eqref{eq:ss-twomass-vec} for $k_3=0.001$. The two curves agree relatively well until $\psi$ is close to the plateau region in $\HH(\psi)$ shown in Figure \ref{fig:coupled-H_psi}. During this region, $\HH(\psi)$ is nearly zero so the convergence to the anti-phase synchrony is very slow. Hence in the right panel of Figure \ref{fig:coupled-simu}, we only include the time evolution of $\psi$ for $t\in [0, 80000]$ which is not long enough for the phase difference to converge to the anti-synchrony state. From the plot, we can see that $\psi$ increases rapidly to about $4$ over the first $2000$ unit of time. After that, the converging rate significantly slows down because $\psi$ enters the plateau region of $\HH(\psi)$.  On the other hand, the shape of the $\HH(\psi)$ curve also suggests that the anti-synchronous point is near-neutrally stable, which may explain why  the prediction of the convergence rates near the plateau region is less accurate than the prediction in non-plateau regions. Improving the accuracy needs further work, such as accounting for higher order effects that are neglected in our first-order phase reduction of coupled oscillators, which is beyond the scope of this paper.

		\begin{figure}[!t]
			\centering
			\includegraphics[width=0.45\textwidth]{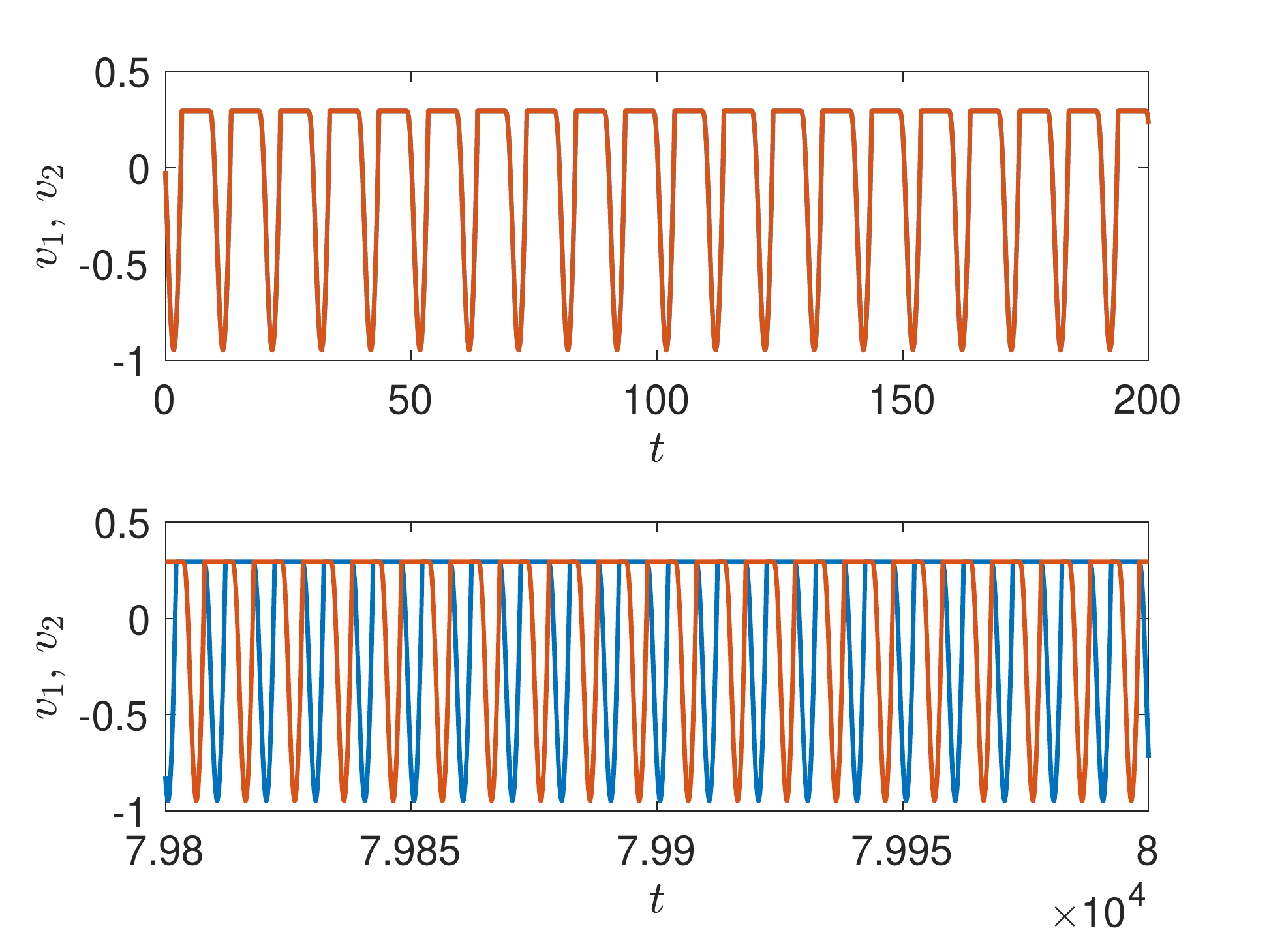}
			\includegraphics[width=0.45\textwidth]{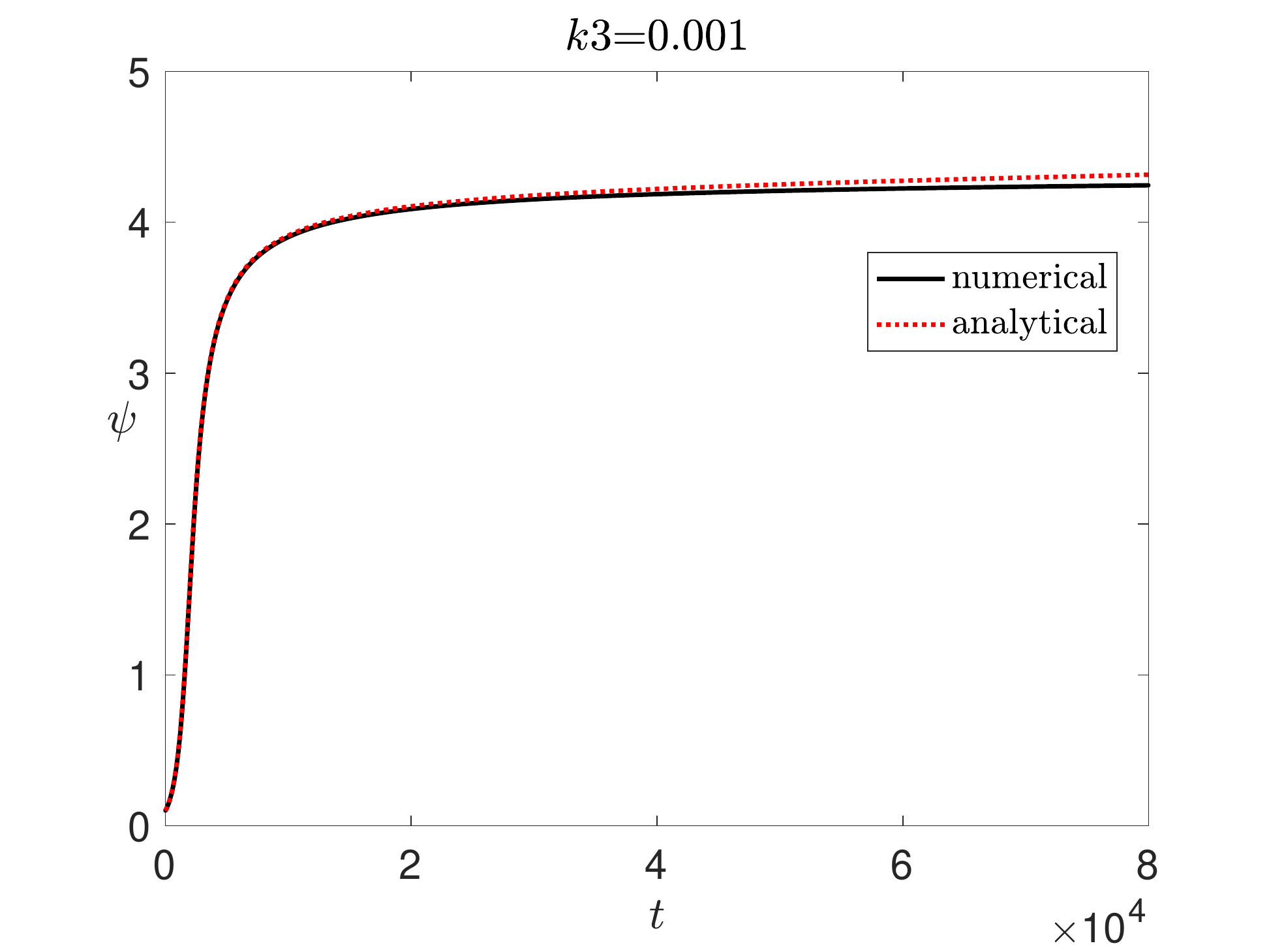}
			\caption{\RED{Simulation result of the coupled stick-slip system \eqref{eq:ss-twomass-vec} when the coupling strength $k_3=0.001$ and other parameters are the same as in Figure \ref{fig:coupled-H_psi}. Left: Time series of velocities of the two oscillators: initial conditions at a phase difference of 0 ($\psi=0$) (top panel) approach a phase difference of $4.2449$ after duration time $80000$ (bottom panel) and eventually converge to a phase-locked solution with a phase difference of $T_0/2$ where $T_0=10.02$ is the period of the uncoupled stick-slip oscillator (not shown here since the convergence rate of $\psi$ near the plateau region of $\HH(\psi)$ as shown in Figure \ref{fig:coupled-H_psi} is very slow).  Right: Time series of the predicted phase difference (red) and the actual phase difference (black) for time over $[0, 80000]$.}}
			\label{fig:coupled-simu}
		\end{figure}
	}

	\section{Discussion}\label{sec:discussion}
	
	Rhythmic motions making and breaking contact with a constraining boundary, and subject to external perturbations, arise in motor control systems such as walking, running, scratching, biting and swallowing, as well as other natural and engineered hybrid systems \RED{\citep{branicky1998,BRS2015}}. Dynamical systems describing such rhythmic motions are therefore nonsmooth and often exhibit limit cycle trajectories with sliding components.
	In smooth dynamical systems, classical analysis for understanding the change in periodic limit cycle orbits under weak perturbation relies on the Jacobian linearization of the flow near the limit cycle.
	These methods do not apply directly to nonsmooth systems, for which the Jacobian matrices are not well defined.
	In this work, we describe for the first time the infinitesimal phase response curves (iPRC) for limit cycles with sliding components (LCSC).
	Moreover, we give a rigorous derivation of the jump matrix for the iPRC at the hard boundary crossing point.
	We also report, for the first time, how the presence of a liftoff point, where a limit cycle leaves a constraint surface, can create a nondifferentiable ``kink'' in the asymptotic phase function, propagating backwards in time along an osculating trajectory (see Figure~\ref{fig:LC-2d-isochrons}A).
	Most significantly, we have developed the infinitesimal shape response curve (iSRC) to analyze the joint variation of both shape and timing of limit cycles with sliding components, under parametric perturbations.
	We show that taking into account local timing sensitivity \emph{within} a switching region improves the accuracy of the iSRC over global timing analysis alone.
	This improvement in accuracy is facilitated by our introduction of a novel \textit{local timing response curve} (lTRC) measuring the timing sensitivity of an oscillator within a given local region.

	Our results clarify an important distinction between the effects of the boundary encounter on the timing and shape changes in limit cycles with sliding components.
	We have extended the iPRC developed for smooth limit cycle systems to the LCSC case, presented here as Theorem \ref{thm:main}.
	In addition, our analysis yields an explicit expression for the iPRC jump matrix that characterizes the behavior of the iPRC at the landing and liftoff points.
	Surprisingly, we find that the iPRC experiences no discontinuity when the trajectory first contacts a hard boundary, while the variational equation suffers a discontinuity, captured by the saltation matrix.
	Even more interesting, at the liftoff point -- where the saltation matrix for the variational problem is trivial -- the iPRC \emph{does} show a discontinuous change, captured by a nontrivial jump matrix.
	Specifically, there is a discontinuous jump from zero to a nonzero normal component in the iPRC.
	Consequently, numerical evaluation of the iPRC must be obtained by backward integration along the limit cycle, as discussed in \S\ref{sec:algorithm-iprc}.
	Finally, we find that both the iPRC and the variational dynamics have zero normal components during the sliding component of the limit cycle, due to dimensional compression at the hard boundary.
	
	\RED{Limit cycles with sliding components can sometimes arise as the singular limits of smooth singularly perturbed systems \citep{Jeffrey2018,JW2020}. Specifically, \citep{Jeffrey2018} shows that a piecewise smooth system can be understood as a singular perturbation problem in the limiting situation by blowing up the discontinuity into a switching layer.    
		To the best of our knowledge, this literature does not address phase response curve and variational dynamics.
		In principle, one might obtain results analogous to those we present here by first analyzing a smooth system and subsequently taking the singular limit.
		While such an undertaking would be both interesting and challenging, our methods avoid the associated technical challenges by calculating the iSRC, iPRC and lTRC for the nonsmooth system directly.
		Moreover, our results may provide some insights into phase response curve and variational dynamics of singularly perturbed systems that exhibit LCSC in the singular limit. For example, \citep{SK93,izhikevich2000} analyze and predict synchronization properties of relaxation oscillators using the results from iPRC analysis in the singular relaxation limit. Whether a similar relationship holds between singularly perturbed systems and LCSCs in singular limits has yet to be understood. }	
	
	Standard variational and phase response curve analysis typically neglects changes in timing or shape, focusing instead on only one of the two aspects \citep{kuramoto1975}.
	However, in many applications such as motor control systems, both the shape and timing of the trajectory are often affected under slow or parametric perturbations.
	In this paper, we consider both timing and shape aspects using the iSRC\RED{, a first-order approximation to the change in shape of the limit cycle under a parameteric perturbation}.
	We have discussed two ways of incorporating timing changes into the iSRC: uniform timing rescaling based on the global timing analysis (iPRC) and piecewise uniform timing rescaling based on the local timing analysis (lTRC).
	As demonstrated in the planar system example in \S\ref{sec:toy-model}, when the trajectory exhibits approximately constant timing sensitivities, the iSRC with global timing rescaling is good enough for approximating the shape change (see Figure \ref{fig:change-shape}); otherwise, we need take into account local timing changes to increase the accuracy of the iSRC (see Figure \ref{fig:shape-response}).
	LCSC with piecewise timing sensitivities naturally arise in many motor control systems due to nonuniform perturbations \RED{as well as the stick-slip mechanical system as studied in \S\ref{sec:stick-slip}.}
	Local timing analysis (lTRC) will then provide a better understanding of such systems compared with the global timing analysis (iPRC).
	\RED{Ours is not the first work to characterize linear responses of limit cycles to parameteric perturbations. \cite{Taylor2008} defines the ``parametric impulse phase response curve (pIPRC)" to capture the timing sensitivity of limit cycle systems to parametric perturbations, which is estimated in our paper through the iPRC as described at the end  of \S\ref{sec:smooth-inst}.}

	Other investigators have also considered variational \citep{bernardo2008,LN2013} and phase response analysis in nonsmooth systems  \citep{shirasaka2017,park2018,CGL18,wilson2019}, but the studies on the iPRC were subject to transverse flow conditions.
	Our work extends the iPRC analysis to the LCSC case in which the transversal crossing condition fails.
	Combined timing and shape responses of limit cycles to perturbations have also been explored in other works.
	\citet{monga2018} examined energy-optimal control of the timing of limit cycle systems including spiking neuron models and models of cardiac arrhythmia.
	They showed that when one of the nontrivial Floquet multipliers of an unperturbed limit cycle system has magnitude close to unity, control inputs based solely on standard phase reduction, which neglects the effect on the shape of the controlled trajectory, can dramatically fail to achieve control objectives.
	They and other authors have introduced  augmented phase reduction techniques that use a system of coordinates (related to the Floquet coordinates) transverse to the limit cycle to improve the accuracy of phase reduction and control \RED{\citep{castejon2013,WM2015,WM2016,WE2018,monga2018b,wilson2019,wilson2020a,wilson2020b,perez-cervera2020}}.
	These methods require the underlying dynamics be smoothly differentiable, and rely on calculation of the Jacobian (first derivative) and in some cases the Hessian (second derivative) matrices \citep{WE2018}. 
	For nonsmooth limit cycle systems with sliding components, our analysis is the first to address the combined effects of shape and timing, an essential element of improved control in biomedical applications as well as for understanding mechanisms of control in naturally occurring motor control systems.
	
	For trajectories with different timing sensitivities in different regions, we rely on the local timing response curve (lTRC) to estimate the relative shift in time in each sub-region, in order to compute the full infinitesimal shape response curve (iSRC).
	Conversely, solving for the lTRC in a given region may also require an understanding of the impact of the perturbation on the entry point associated with that region (see \eqref{eq:local-time-shift}).
	Thus, in general, the iSRC and the lTRC are interdependent.
	While we have not derived a closed-form expression for the shape and timing response in the most general case, we have provided effective algorithms for solving each of them separately, which requires preliminary numerical work to find the trajectory shape shift at the entry point.
	In the future, it may be possible to derive general closed-form expressions for the iSRC and lTRC in systems with distinct timing sensitivities.
	
	While our methods are illustrated using a planar limit cycle system with hard boundaries \RED{and coupled stick-slip systems}, they apply to \RED{higher dimensional and more realistic} systems as well.
	For instance,  preliminary investigations suggest that the  methods developed in this paper are applicable to analyzing the nonsmooth dynamics arising in the control system of feeding movements in the sea slug \textit{Aplysia} \citep{shaw2012,shaw2015,lyttle2017}. More generally, limit cycles with discontinuous trajectories arise in neuroscience (e.g., integrate and fire neurons) and mechanics (e.g., ricochet dynamics). If such systems manifest limit cycles with sliding components, our methods could be combined with variational methods adapted for piecewise continuous trajectories \citep{coombes2012,shirasaka2017}.
	
	It was observed heuristically by \citet{lyttle2017} that sensory feedback could in some circumstances lead to significant robustness against an increase in applied load, in the sense that although modest relative increases in external load (\textit{c.}~20\%) led to comparable changes in both the timing and shape of trajectories, the net effect on the performance (rate of intake of food) was an order of magnitude smaller (\textit{c.}~1\%).
	Similarly, \citet{diekman2017} showed that in a model for control of a central pattern generator regulating the breathing rhythm, mean arterial partial pressure of oxygen (PPO$_2$) remained approximately constant under changing metabolic loads when chemosensory feedback from the arterial PPO$_2$ to the central pattern generator was present, but varied widely otherwise \citep{diekman2017}.
	Understanding how rhythmic biological control systems respond to such perturbations and maintain robust, adaptive performance is one of the fundamental problems within theoretical biology.
	Solving these problems will then require variational analysis along the lines we develop here.
	Nonsmooth dynamics arise naturally in many biological systems \citep{AiharaSuzuki2010,coombes2012}, and thus, the approach in this paper is likely to have broad applicability to many other problems in biology.

	\section*{Acknowledgement}
	This work was made possible in part by grants from the National Science Foundation (DMS-1413770, DEB-1654989, IOS-174869 \RED{and IOS-1754869 to H.J.C}). P.J.T thanks the Oberlin College Department of Mathematics for research support.  This research has been supported in part by the National Science Foundation Grant DMS-1440386 to the Mathematical Biosciences Institute.

	\newpage
	\appendix

	\section{Table of Common Symbols}\label{ap:symbols}
	\begin{center}
		\begin{tabular}{@{}cl@{}}
			
			\toprule
			
			\textbf{Symbol} & \textbf{Meaning} \\
			
			\midrule
			
			$\mbx$                                       & state variables \\
			$t$                                          & time \\
			$\theta(t)$                                  & phase of a limit cycle \\
			$\phi(\mbx)$                                 & asymptotic phase of a stable limit cycle \\
			
			\midrule
			
			$F(\mbx)$                                    & unperturbed velocity vector field \\
			$\gamma(t)$                                  & unperturbed limit cycle solution \\
			$T$                                          & period of the unperturbed limit cycle \\
			
			\midrule
			
			$\varepsilon P$                              & small instantaneous perturbation vector \\
			$\tilde{\gamma}(t)$                          & trajectory near limit cycle after instantaneous perturbation \\
			$\mbu(t) \simeq \tilde{\gamma}(t) - \gamma(t)$    & displacement from limit cycle after instantaneous perturbation \\
			
			\midrule
			
			$\varepsilon$                                & sustained (parametric) perturbation \\
			$F_\varepsilon(\mbx)$                        & perturbed velocity vector field \\
			$\gamma_\varepsilon(t)$                      & perturbed limit cycles solution \\
			$T_\varepsilon$                              & period of the perturbed limit cycle \\
			
			\midrule
			
			$F_0 = F$                                    & zeroth-order term of Taylor expansion of $F_\varepsilon$ around $\varepsilon=0$ \\
			$\gamma_0 = \gamma$                          & zeroth-order term of Taylor expansion of $\gamma_\varepsilon$ around $\varepsilon=0$ \\
			$T_0 = T$                                    & zeroth-order term of Taylor expansion of $T_\varepsilon$ around $\varepsilon=0$ \\
			
			\midrule
			
			$F_1 =      \partial F_\varepsilon /
			\partial \varepsilon
			\big|_{\varepsilon=0}$           & first-order term of Taylor expansion of $F_\varepsilon$ around $\varepsilon=0$ \\
			$\gamma_1 = \partial \gamma_\varepsilon /
			\partial \varepsilon
			\big|_{\varepsilon=0}$           & first-order term of Taylor expansion of $\gamma_\varepsilon$ around $\varepsilon=0$, \\
			& also called the infinitesimal shape response curve (iSRC) \\
			$T_1 =      \partial T_\varepsilon /
			\partial \varepsilon
			\big|_{\varepsilon=0}$           & first-order term of Taylor expansion of $T_\varepsilon$ around $\varepsilon=0$ \\
			
			\midrule
			
			$DF$; $D_\mathbf{w}\phi$                                         & Jacobian matrix; directional derivative of $\phi$ in $\mathbf{w}$ direction \\
			$I$                                          & identity matrix \\
			$S$                                          & saltation matrix (for variation equation) \\
			$J$                                          & jump matrix (for adjoint equation) \\
			$\mathcal{J}$                                & time-reversed jump matrix (for adjoint equation) \\
			
			\midrule
			
			$\Sigma^i$                                   & boundary $i$ \\
			$\RR^j$                                      & region $j$ \\
			$F^j(\mbx)$                                  & velocity vector field in region $j$ \\
			
			\midrule
			
			$\nu_\varepsilon = T_0/T_\varepsilon$        & relative frequency of perturbed limit cycle \\
			$\nu_1 = T_1/T_0$                            & first-order term of Taylor expansion of $\nu_\varepsilon$ around $\varepsilon=0$, \\
			& also called the relative change in frequency \\
			$\mathcal{T}^j$                           & time remaining in region $j$ along a trajectory \\
			
			\midrule
			
			$\mbu(t)$                                    & variational dynamics governed by \eqref{eq:var} \\
			$\z(t)=\nabla_\mbx\phi(\gamma(t))$           & infinitesimal phase response curve (iPRC) governed by \eqref{eq:prc} \\
			$\gamma_1(t)$                                & infinitesimal shape response curve (iSRC) governed by \eqref{eq:src} \\
			$\eta^j(t) = \nabla_\mbx \mathcal{T}^j(\gamma(t))$   & local timing response curve (lTRC) governed by \eqref{eq:ltrc} \\
			
			\bottomrule
		\end{tabular}
	\end{center}

	\RED{\section{Proof of Lemma \ref{lem:isrc}}
		\label{ap:isrc-ini}

		In this section we prove Lemma \ref{lem:isrc}, which we restate for the reader's convenience.
		
		\textbf{Lemma} {\em 
			Let $\gamma^\textbf{a}_1(t)$ and $\gamma^\textbf{b}_1(t)$ be
			two $T_0$-periodic solutions to the iSRC equation \eqref{eq:src} for a smooth vector field $F_0$ with a hyperbolically stable limit cycle $\gamma_0(t)$.  
			Then, their difference satisfies $\gamma^\textbf{b}_1(t)-\gamma^\textbf{a}_1(t)=\varphi F_0(\gamma_0(t))$, where $\varphi$ is a constant representing a fixed phase offset.}	
		
		\begin{proof}
			Consider
			\begin{equation}\label{eq:F0}
			\mbx'=F_0(\mbx)
			\end{equation} 
			with $\mbx(0)=\mbx_0\in \mathbb{R}^n$. Let $\Phi(t,0)$ be the fundamental matrix solution. 
			Then $\Phi(t,0)$ satisfies $\Phi'(t,0)= DF_0(\mbx(t))\Phi(t,0)$ and $\Phi(0,0)=I$, where $\mbx=\gamma_0(t)$ is the unperturbed limit cycle solution.
			Suppose the monodromy matrix $M=\Phi(T_0,0)$ is diagonalizable with eigenvalues $\{\mu_i, i=1,\cdots,n\}$ associated with linearly independent eigenvectors $\{\mathbf{v}_i, i=1,\cdots,n\}$. The eigenvalues $\mu_i$ are often
			referred to as Floquet multipliers of the periodic orbit solution $\gamma_0(t)$ of \eqref{eq:F0} \citep{meiss07}. Since $\gamma_0(t)$ is hyperbolically stable, $M$ has a single trivial Floquet multiplier. Without loss of generality, we assume $\mu_1=1$ and hence $\mathbf{v}_1=F_0(\gamma_0(0))$.
			
			Let the vector $\eta_i(t)$ be the solution to the variational equation
			\[
			\eta_i' = DF_0(\gamma_0(t))\eta_i
			\]
			that starts along the $i$-th Floquet eigenvector direction $\eta_i(0)=\mathbf{v}_i$. It follows that $\eta_i(T_0) = \Phi(T_0,0)\eta_i(0)=M \mathbf{v}_i = \mu_i \mathbf{v}_i$. For simplicity, we denote $A(t) = DF_0(\gamma_0(t))$ hereafter.
			
			Let $\rho_i = \ln(\mu_i)/T_0$ and let $\mathbf{q}_i(t)=e^{-t\rho_i}\eta_i(t)$ be the rescaled version of the trajectory $\eta_i(t)$, which is the $i$-th Floquet coordinate \citep{meiss07}. It follows from direct calculations that $\mathbf{q}_i(t)$ is periodic so that $\mathbf{q}_i(t) = \mathbf{q}_i(t+T_0)$ and satisfies the initial value problem
			\begin{equation}\label{eq:qi}
			\mathbf{q}_i(t)' = A(t)\mathbf{q}_i(t) - \rho_i \mathbf{q}_i(t)
			\end{equation}
			with $\mathbf{q}_i(0)=\mathbf{v}_i$.  Since $\rho_1=0$, the first Floquet coordinate $\mathbf{q}_1(t)$ satisfies the initial value problem
			\begin{equation}\label{eq:q1}
			\mathbf{q}_1'=A(t)\mathbf{q}_1,\quad \mathbf{q}_1(0)=F_0(\gamma_0(t)). 
			\end{equation} 
			Using the chain rule, upon differentiating $F_0(\gamma_0(t))$, one can show that setting $\mathbf{q}_1(t)=F_0(\gamma_0(t))$ for any $t$, solves the initial value problem \eqref{eq:q1}. 
			Note that \eqref{eq:qi} is similar to the equation satisfied by $\mathcal{I}_i$, the gradient of isostable coordinates that are related to Floquet coordinates $\mathbf{q}_i$ \citep{WM2016,perez-cervera2020}. In fact, direct calculation implies that the relationship between the iPRC and the variational dynamics (see Remark \ref{rem:uTz=const}) also holds for $\mathbf{q}_i$ and $\mathcal{I}_i$; that is, $\frac{d(\mathbf{q}_i^\intercal\mathcal{I}_i)}{dt}=0$ for $i=1,\cdots,n$.
			
			Note that at each $t$, the vector set $\{\mathbf{q}_i(t), i=1,\cdots,n\}$ spans $\mathbb{R}^n$. 
			Therefore we can write $\gamma_1(t)$, the general solution to the iSRC equation \eqref{eq:src}, as a linear combination of $\{\mathbf{q}_i(t), i=1,\cdots,n\}$ with coefficients $a_i(t)$:
			\[
			\gamma_1(t) = \sum_{i=1}^{n} a_i(t) \mathbf{q}_i(t).
			\]
			
			Let $Q(t)$ be the $n\times n$ matrix $Q(t)=\left(\mathbf{q}_1(t)|\cdots| \mathbf{q}_n(t)\right)$, and let  $R=\mathrm{diag}(\rho_1,\cdots,\rho_n)$ be the diagonal matrix with $\{\rho_1,\cdots, \rho_n\}$ as the diagonal entries, and $\mathbf{a}(t) =[a_1(t),\cdots,a_n(t)]^\intercal$. 
			Then 
			\begin{equation}\label{eq:gamma_1}
			\gamma_1(t)=Q(t)\mathbf{a}(t)
			\end{equation}
			and \eqref{eq:qi} can be rewritten as
			\begin{equation}\label{eq:qi_matrix}
			Q(t)' = A(t)Q(t)-Q(t)R.    
			\end{equation}
			
			Differentiating both sides of \eqref{eq:gamma_1} and substituting in \eqref{eq:qi_matrix} leads to 
			\begin{align*}
			\gamma_1'(t)&= Q'(t)\mathbf{a}(t) +Q(t)\mathbf{a}'(t)\\
			&=(A(t)Q(t)-Q(t)R )\mathbf{a}(t) +Q(t)\mathbf{a}'(t)\\
			&=(Q(t)\mathbf{a}'(t)-Q(t)R \mathbf{a}(t)) + A(t)Q(t)\mathbf{a}(t)\\
			&=(Q(t)\mathbf{a}'(t)-Q(t)R \mathbf{a}(t)) + A(t)\gamma_1(t).
			\end{align*}
			On the other hand, by \eqref{eq:src} we have
			\begin{eqnarray}
			\gamma_1'(t)
			&=& A(t) \gamma_1(t) +\mathbf{c}(t)\nonumber,
			\end{eqnarray}
			where $\mathbf{c}(t)=\nu_1 F_0(\gamma_0(t)) +\frac{\partial F_\varepsilon(\gamma_0(t))}{\partial \varepsilon}\Big|_{\varepsilon=0}$.
			It follows that 
			\begin{equation}\label{eq:ai'}
			Q(t)\mathbf{a}'(t)= Q(t)R \mathbf{a}(t)+ \mathbf{c}(t).  
			\end{equation}
			Since $\{\mathbf{q}_i(t)\}_{i=1}^n$ spans $\mathbb{R}^b$ for each time $t\in [0, T_0]$, the matrix $Q(t)$ is invertible at each $t$. 
			Thus multiplying both sides of \eqref{eq:ai'} by $Q(t)^{-1}$ gives
			\begin{equation}\label{eq:ai'-2}
			\mathbf{a}'(t)= R \mathbf{a}(t) + Q(t)^{-1}\mathbf{c}(t).  
			\end{equation}
			
			Suppose two different iSRC curves are given by $\gamma_1^a(t)=\sum_{i=1}^n a_i(t)\mathbf{q}_i(t)=Q(t)\mathbf{a}(t)$ and $\gamma_1^b(t)=\sum_{i=1}^n b_i(t)\mathbf{q}_i(t)=Q(t)\mathbf{b}(t)$. Then
			\[
			\mathbf{b}'(t)-\mathbf{a}'(t) = R\mathbf{b}(t)-R\mathbf{a}(t)= \Matrixc{\rho_1&&\\ &\ddots & \\ && \rho_n}(\mathbf{b}(t)-\mathbf{a}(t)).
			\]
			It follows that for $i=1,\cdots,n$ we have $b_i(t)-a_i(t)=C_ie^{\rho_i t}$ for some constant $C_i$. 
			Note that $\mathbf{a}(t)$ and $\mathbf{b}(t)$ are both $T_0$-periodic.
			Therefore $C_i=C_i e^{\rho_i T_0}$. 
			So either $C_i=0$, and hence $a_i(t)\equiv b_i(t)$, or else $\rho_i=0$. 
			However, recall there is only one trivial multiplier $\mu_1=1$, so that only $\rho_1=0$ and $\rho_i\neq 0$ for $i=2,\cdots,n$. Hence, $a_i(t)\equiv b_i(t)$ for $i=2,\cdots,n$; that is, there exists some constant $\phi$ such that $\mathbf{b}(t)-\mathbf{a}(t)=\phi \mathbf{e}_1$.
			
			Thus, 
			\[
			\gamma_1^b(t)-\gamma_1^a(t) = Q(t) \phi \mathbf{e}_1 =\phi \mathbf{q}_1(t).
			\]
			Consequently, the two iSRC curves only differ along the first Floquet coordinate direction $\mathbf{q}_1(t)=F_0(\gamma_0(t))$ and hence only differ by a shift in phase along the direction of the limit cycle $\gamma_0(t)$:
			\[
			\gamma_1^b(t)-\gamma_1^a(t) =\phi F_0(\gamma_0(t))
			\]
			where $\phi$ is the constant phase shift introduced by the initial conditions $\gamma_1^b(0)-\gamma_1^a(0) =\phi F_0(\gamma_0(0))$. 
			
	\end{proof}}
	
	\newpage 
	\section{Derivation of Equation \ref{eq:local-time-shift}}\label{ap:T1}
	
	This section establishes equation \eqref{eq:local-time-shift}, which specifies the first-order change in the transit time through region I, or $T_1^\text{I}$:
	$$T^{\rm I}_{1} = \eta^{\rm I}(\mbx^{\rm in})\cdot \frac{\partial \mbx_{\varepsilon}^{\rm in}}{\partial \varepsilon}\Big|_{\varepsilon=0}+\int_{t^{\rm in}}^{t^{\rm out}}\eta^{\rm I}(\gamma(t))\cdot \frac{\partial F_\varepsilon(\gamma(t))}{\partial \varepsilon}\Big|_{\varepsilon=0}dt,
	$$ 
	Recall
	$\mathcal{T}^{\rm I}(\mbx)$ 
	is the time remaining until exiting region I through $\Sigma^{\rm out}$, under the unperturbed vector field, starting from location $\mbx$;   $\eta^{\rm I}:= \nabla \mathcal{T}^{\rm I}(\mbx)$ is the local timing response curve (lTRC) for region I, defined for the component of the trajectory lying within region I, i.e.~for times $t\in [t^{\rm in}, t^{\rm out}]$; and
	$\mbx_\varepsilon^{\rm in}$ is the coordinate of the perturbed entry point into region I.

	We consider a single region $\mathcal{R}$ with entry surface $\Sigma^{\rm in}$ and exist surface $\Sigma^{\rm out}$.
	We assume that these two surfaces are fixed, independent of static perturbation with size $\epsilon$.
	The limit cycle solution $\mbx=\gamma_\epsilon(\tau)$ satisfies
	\[
	\frac{d\mbx}{d\tau}=F_{\varepsilon}(\mbx)
	\]
	where $\tau$ is the time coordinate of the perturbed trajectory.
	Moreover,
	$\gamma_\epsilon(\tau)$ enters $\mathcal{R}$ at $\mbx^\text{in}_\epsilon\in\Sigma^{\rm in}$ when $\tau=t^\text{in}_\epsilon$ and exits at $\mbx^\text{out}_\epsilon\in\Sigma^{\rm out}$ when $\tau=t^\text{out}_\epsilon$.
	Since the system is autonomous, we are free to choose the reference time along the limit cycle orbit.
	For convenience of calculation, we set $t^\text{out}_\epsilon\equiv 0$ for all $\epsilon$.
	
	Denote the transit time that $\gamma_\epsilon$ spends in $\mathcal{R}$ by $T^\RR_\varepsilon$.
	It follows that $t^\text{in}_\epsilon = -T^\RR_{\epsilon}$, where $\epsilon$ can be $0$.
	Assuming that  the transit time has a well behaved expansion in $\epsilon,$ we write 
	\begin{align}\label{eq:Teps}
	T^\RR_\epsilon&=T^\RR_0+\epsilon T^\RR_1+O(\epsilon^2)
	\end{align}
	where $T^\RR_0$ is the transit time for the unperturbed trajectory and $T^\RR_1$ is  the linear shift in the transit time.
	In the rest of this section, we drop the superscript $\RR$ on $T^\RR_\varepsilon,\,T^\RR_0$ and $T^\RR_1$ for simplicity.
	
	Our goal is to prove that $T_1$ is given by \eqref{eq:local-time-shift}.
	We do this in two steps.
	First, we show that the transit time $T_{\epsilon}$ can be expressed in terms of the perturbed vector field and perturbed local timing response curve (see \eqref{eq:integral-for-transit-time}).
	Second, we expand the expression for $T_{\epsilon}$ to first order in $\varepsilon$ to obtain the expression for $T_1$.
	
	Since the time remaining to exit, denoted as $\mathcal{T}_\varepsilon$, decreases at a constant rate along trajectories, for arbitrary $\epsilon$ we have 
	\begin{equation}\label{eq:dTdt=-1}
	-1=\frac{d\mathcal{T}_\varepsilon}{d\tau}=F_\varepsilon(\gamma_\varepsilon(\tau))\cdot\eta_\varepsilon(\gamma_\varepsilon(\tau)),
	\end{equation}
	where $\eta_\varepsilon(\mbx)=\nabla\mathcal{T}_\varepsilon(\mbx)$
	is defined as the local timing response curve under perturbation.
	By \eqref{eq:dTdt=-1},  the transit time $T_\varepsilon$ is therefore given by \begin{equation}\label{eq:integral-for-transit-time}
	T_\varepsilon=\int_{\tau=t^\text{out}_\epsilon}^{t^\text{in}_\epsilon}F_\varepsilon(\gamma_\epsilon(\tau))\cdot \eta_\epsilon(\gamma_\epsilon(\tau))\,d\tau.
	\end{equation}
	In this expression, we integrate
	\emph{backwards in time} along the limit cycle trajectory, from the egress point $\mbx^\text{out}_\varepsilon$ at time $t^\text{out}_\varepsilon$, to the ingress point $\mbx^\text{in}_\varepsilon$ at time $t^\text{in}_\varepsilon$:

	For $\epsilon=0$, and taking into account \eqref{eq:dTdt=-1}, this integral reduces to 
	\begin{equation}
	T_0=\int_{\tau=t^\text{out}_0}^{t^\text{in}_0} F_0(\gamma_0(\tau))\cdot\eta_0(\gamma_0(\tau))\,d\tau=\int_{\tau=t^\text{out}_0}^{t^\text{in}_0}(-1)\,d\tau = t^\text{out}_0-t^\text{in}_0=0-(-T_0),
	\end{equation}
	since $t^\text{in}_0= -T_0$ and $t_\epsilon^\text{out}\equiv 0$.
	
	In order to derive an expression for $T_1$, the first order shift in the transit time, we need to expand \eqref{eq:integral-for-transit-time} to first order in $\epsilon$.
	To this end, we need to know the Taylor expansions for all terms in \eqref{eq:integral-for-transit-time}.
	
	Suppose we can expand $F_\epsilon$,  $\mathcal{T}_{\varepsilon}$, and $\eta_\varepsilon$ as follows:
	\begin{eqnarray}\label{eq:eps}
	\Eqn{F_{\varepsilon}(\mbx)&=&F_0(\mbx)+\epsilon F_1(\mbx)+O(\epsilon^2),&\text{ as }\epsilon\to 0,\\
		\mathcal{T}_\varepsilon(\mbx)&=&\mathcal{T}_0(\mbx)+\epsilon\mathcal{T}_1(\mbx)+O(\epsilon^2),&\text{ as }\epsilon\to 0,\\
		\eta_\varepsilon(\mbx)&=&\eta_0(\mbx)+\epsilon\eta_1(\mbx)+O(\epsilon^2),&\text{ as }\epsilon\to 0, }
	\end{eqnarray}
	where $\eta_0(\mbx) = \nabla\mathcal{T}_0(\mbx)$ is the unperturbed local timing response curve.
	
	Following the idea of deriving the infinitesimal shape response curve in \S\ref{sec:smooth-sust}, we write the portion of the perturbed limit cycle trajectory within region $\mathcal{R}$ in terms of the unperturbed limit cycle, plus a small correction, 
	\begin{align}\label{eq:gammaeps}
	\gamma_\varepsilon(\tau)&=\gamma\left(\nu_\varepsilon\tau\right)+\epsilon \gamma_1\left(\nu_\varepsilon\tau\right)+O(\epsilon^2)
	\end{align}
	where $-T_\varepsilon\le \tau \le 0$ and $\nu_\varepsilon= \frac{T_0}{T_\varepsilon}$.
	
	Now we expand  \eqref{eq:integral-for-transit-time} to first order
	\begin{align}
	T_\epsilon 
	&=\int_{\tau=0}^{-T_{\varepsilon}}\Big[ F_0(\gamma_0(\nu_\varepsilon\tau))+\epsilon DF_0(\gamma_0(\nu_\varepsilon\tau))\cdot \gamma_1(\nu_\varepsilon\tau)+\epsilon F_1(\gamma_0(\nu_\varepsilon\tau)) \Big]\cdot\\ \nonumber
	&\quad\quad\quad\Big[\eta_0(\gamma_0(\nu_\varepsilon\tau))+\epsilon D\eta_0(\gamma_0(\nu_\varepsilon\tau)) \cdot\gamma_1(\nu_\varepsilon\tau)+\epsilon\eta_1(\gamma_0(\nu_\varepsilon\tau))  \Big]d\tau
	+O(\epsilon^2)\\ \nonumber
	&=\int_{\tau=0}^{-T_{\varepsilon}}  F_0(\gamma_0(\nu_\varepsilon\tau))\cdot \eta_0(\gamma_0(\nu_\varepsilon\tau)) d\tau + \epsilon \Big[ F_0(\gamma_0(\nu_\varepsilon\tau)) \cdot\eta_1(\gamma_0(\nu_\varepsilon\tau))+ F_1(\gamma_0(\nu_\varepsilon\tau))\cdot  \eta_0(\gamma_0(\nu_\varepsilon\tau)) \Big] d\tau + \\\nonumber  
	& \quad\quad\quad \epsilon\Big[ F_0(\gamma_0(\nu_\varepsilon\tau)) \cdot D\eta_0(\gamma_0(\nu_\varepsilon\tau)) \cdot\gamma_1(\nu_\varepsilon\tau) +  DF_0(\gamma_0(\nu_\varepsilon\tau))\cdot \gamma_1(\nu_\varepsilon\tau) \cdot \eta_0(\gamma_0(\nu_\varepsilon\tau)) \Big]  d\tau +O(\epsilon^2)
	\\ \nonumber
	&=\frac{1}{\nu_\varepsilon}\int_{t=0}^{-T_{0}}  F_0(\gamma_0(t))\cdot \eta_0(\gamma_0(t)) dt + \epsilon \Big[ F_0(\gamma_0(t)) \cdot\eta_1(\gamma_0(t))+ F_1(\gamma_0(t))\cdot  \eta_0(\gamma_0(t)) \Big] dt + \\\nonumber  
	& \quad\quad\quad \epsilon\Big[ F_0(\gamma_0(t)) \cdot D\eta_0(\gamma_0(t)) \cdot\gamma_1(t) +  DF_0(\gamma_0(t))\cdot \gamma_1(t) \cdot \eta_0(\gamma_0(t)) \Big]  dt +O(\epsilon^2)
	\\ \nonumber
	\end{align}
	
	To order $O(1)$, we recover
	\begin{equation}
	T_0=\int_{t=0}^{-T_0} F_0(\gamma_0(t))\cdot\eta_0(\gamma_0(t))\,dt.
	\end{equation}
	This leads to 
	$
	T_{\varepsilon} = \frac{1}{\nu_{\varepsilon}}T_0,
	$
	as required for consistency.
	We are therefore left with
	\begin{align}\label{eq:nu1_on_one_side} 
	0&=\int_{t=0}^{-T_{0}}    \Big[ F_0(\gamma_0(t)) \cdot\eta_1(\gamma_0(t))+ F_1(\gamma_0(t))\cdot  \eta_0(\gamma_0(t)) \Big] dt  \\\nonumber  
	& + \int_{t=0}^{-T_{0}}   \Big[ F_0(\gamma_0(t)) \cdot D\eta_0(\gamma_0(t)) \cdot\gamma_1(t) +  DF_0(\gamma_0(t))\cdot \gamma_1(t) \cdot \eta_0(\gamma_0(t)) \Big]  dt +O(\epsilon)
	\\ \nonumber
	&=\int_{t=0}^{-T_{0}}    \Big[ F_0(\gamma_0(t)) \cdot\eta_1(\gamma_0(t))+ F_1(\gamma_0(t))\cdot  \eta_0(\gamma_0(t)) \Big] dt  \\\nonumber  
	& + \int_{t=0}^{-T_{0}}   \Big[F_0(\gamma_0(t))^\intercal D\eta_0(\gamma_0(t)) + \eta_0(\gamma_0(t))^\intercal DF_0(\gamma_0(t))\Big] \cdot\gamma_1(t)  dt +O(\epsilon)
	\\ \nonumber
	\end{align}
	where the second equality follows from rearranging orders of factors in the second integral.
	
	Note that since $F_0\cdot \eta_0\equiv -1$ everywhere, we have the identity
	\begin{align}
	0=\frac\partial{\partial \mbx_j}\left(\sum_i \eta^i F^i\right)=\sum_i\frac{\partial \eta^i}{\partial \mbx_j} F^i+\sum_i\eta^i\frac{\partial F^i}{\partial \mbx_j}
	\end{align}
	where $F^i$ and $\eta^i$ are the $i$-th components for $F_0$ and $\eta_0$; $\mbx_j$ denotes the $j$th component of $\mbx$ for $j\in \{1,\cdots, n\}$.
	It follows that $F_0^\intercal(D\eta_0)+\eta_0^\intercal(D F_0)=0$ in \eqref{eq:nu1_on_one_side}, leaving only 
	\begin{equation}
	0=\int_{t=0}^{-T_0}  \Big[F_0(\gamma_0(t)) \cdot\eta_1(\gamma_0(t))+ F_1(\gamma_0(t))\cdot  \eta_0(\gamma_0(t))\Big] \,dt.
	\end{equation}
	Since $F_0(\gamma_0(t))=d\gamma_0/dt$ and $\eta_1(\mbx)=\partial \eta_\epsilon(\mbx)/\partial\epsilon |_{\epsilon=0} = \partial \nabla\mathcal{T}_\epsilon(\mbx)/\partial \epsilon|_{\epsilon=0}$, 
	\begin{align}
	\int_{t=0}^{-T_0} 
	F_0(\gamma_0(t))\cdot \eta_1(\gamma_0(t)) \,dt \nonumber
	&= \int_{t=0}^{-T_0}\left(\frac{d\gamma_0}{d t}\right)\cdot\left.\left(\frac\partial{\partial\epsilon}\left[\nabla \mathcal{T}_\epsilon(\gamma_0(t)) \right]\right)\right|_{\epsilon=0}\,dt\\\nonumber
	&=\int_{t=0}^{-T_0}\left(\frac{d\gamma_0}{dt}\right)\cdot\nabla\left.\left(\frac\partial{\partial\epsilon}\left[ \mathcal{T}_\epsilon(\gamma_0(t)) \right]\right)\right|_{\epsilon=0}\,dt\\\nonumber
	&=\int_{t=0}^{-T_0}\frac{d}{dt}\left.\left(\frac\partial{\partial\epsilon}\left[ \mathcal{T}_\epsilon(\gamma_0(t)) \right]\right)\right|_{\epsilon=0}\,dt\\\nonumber
	&=
	\left.\left(\frac\partial{\partial\epsilon}\left[ \mathcal{T}_\epsilon(\mbx^\text{in}_0) \right]\right)\right|_{\epsilon=0} -\left.\left(\frac\partial{\partial\epsilon}\left[ \mathcal{T}_\epsilon(\mbx^\text{out}_0) \right]\right)\right|_{\epsilon=0} \\\nonumber
	&=\left.\left(\frac\partial{\partial\epsilon}\left[ \mathcal{T}_\epsilon(\mbx^\text{in}_0) \right]\right)\right|_{\epsilon=0} -0\\\nonumber
	&=\mathcal{T}_1(\mbx^\text{in}_0).
	\end{align}
	Therefore
	\begin{equation}\label{eq:T1v1}
	\mathcal{T}_1(\mbx^\text{in}_0) = 
	\int_{t=-T_0}^{0} 
	F_1(\gamma_0(t))\cdot \eta_0(\gamma_0(t)) dt =\int_{t=t_0^{\rm in}}^{t_0^{\rm out}} 
	F_1(\gamma_0(t))\cdot \eta_0(\gamma_0(t)) dt .
	\end{equation}
	The second equality follows from our convention that $t^\text{in}_0= -T_0$ and $t_\epsilon^\text{out}\equiv 0$.
	
	We notice that
	\begin{align}\label{eq:Teps2}
	T_\epsilon&=\mathcal{T}_\epsilon(\mbx_\epsilon^\text{in})=T_0+ \epsilon\left( \mathcal{T}_1(\mbx^\text{in}_0)+\nabla\mathcal{T}_0(\mbx^\text{in}_0) \cdot\mbx^\text{in}_1 \right), 
	\end{align}
	where we have made use of the Taylor expansion $\mbx^\text{in}_\varepsilon=\mbx^\text{in}_0+\epsilon\mbx^\text{in}_1+O(\epsilon^2),\text{ as }\epsilon\to 0$.
	Equating the first order terms in \eqref{eq:Teps} and \eqref{eq:Teps2} leads to
	\begin{equation}\label{eq:T1v2}
	T_1=\mathcal{T}_1(\mbx^\text{in}_0)+\eta_0(\mbx^\text{in}_0)\cdot\mbx^\text{in}_1.
	\end{equation}
	Substituting \eqref{eq:T1v1} into \eqref{eq:T1v2}, we finally obtain
	\begin{equation}\label{eq:T1proof}
	T_1= \eta_0(\mbx^\text{in}_0)\cdot\mbx^\text{in}_1+ 
	\int_{t=t_0^{\rm in}}^{t_0^{\rm out}} 
	F_1(\gamma_0(t))\cdot \eta_0(\gamma_0(t)) dt
	\end{equation}
	which is \eqref{eq:local-time-shift}, as desired.
	
	\section{Proof of Theorem \ref{thm:main}}\label{ap:proof}
	In this section we present a proof of Theorem \ref{thm:main}, which we restate for the reader's convenience. \RED{As discussed before, parts (a) and (b) are already covered in \citep{filippov1988,bernardo2008}, whereas parts (c) through (d) are our new results. For completeness, we still include parts (a) and (b) as well as our versions of proofs.}
	
	\RED{
		\textbf{Theorem.}
		Consider a general LCSC described locally by~\eqref{eq:1zoneFP} in the neighborhood of a hard boundary $\Sigma$, satisfying Assumption \ref{ass:assumption1} and Assumption \ref{ass:assumption2}. The following properties hold for the variational dynamics $\mbu$ and the iPRC $\z$ along $\Sigma$:
		\begin{enumerate}
			\item[(a)] At the landing point of $\Sigma$, the saltation matrix is $S=I-n n^\intercal$, where $I$ is the
			identity matrix.
			\item[(b)] At the liftoff point of $\Sigma$, the saltation matrix is $S=I$.
			\item[(c)] Along the sliding region within $\Sigma$, the component of $\z$ normal to $\Sigma$ is zero.
			\item[(d)] The normal component of $\z$ is continuous at the landing point.
			\item[(e)] The tangential components of $\z$ are continuous at both landing and liftoff points.
	\end{enumerate}}
	
	\proof
	We choose coordinates $\mbx=(\mbw,v)=(w_1,w_2,\ldots,w_{n-1},v)$
	so that within a neighborhood containing both the landing and liftoff points, the hard boundary corresponds to $v=0$, the interior of the domain coincides with $v>0$, and the unit normal vector for the hard boundary is $\mbn=(0,\ldots,0,1)$.
	Writing the velocity vector $\mbF=(f_1,f_2,\ldots,f_{n-1},g)$ in these coordinates.
	In addition, we use $\mbF^{\rm slide}$ to denote the vector field for points on the sliding region, whereas the dynamics of other points is governed by $\mbF^{\rm int}$.
	The transversal intersection condition for the trajectory entering the hard boundary 
	is $g^\text{int}(\mbx_\text{land},0)<0$
	(cf.~eq.~\eqref{eq:slidingregion}; note that $\mbn$ defined here points in the opposite direction from the \emph{outward} normal vector in \eqref{eq:slidingregion}).
	At points $\mbx\in\mathcal{L}$ on the liftoff boundary, $\mbF^\text{slide}$ and $\mbF^\text{int}$ coincide and we will use whichever notation seems clearer in a given instance.
	Under the nondegeneracy condition at the liftoff point \eqref{eq:nondege-1}, we can further arrange the coordinates $(w_1,\ldots,w_{n-1})$ so that the unit vector  normal to the liftoff boundary $\mathcal{L}$ at the liftoff point is $\ell=(0,\ldots,0,1,0)$, and $g^\text{int}\gtrless 0 \iff w_{n-1}\gtrless 0$.
	With these coordinates, the nondegeneracy condition \eqref{eq:nondege-1} is  $\mbF^{\rm slide}(\mbx_{\rm lift})\cdot \ell=f^\text{slide}_{n-1}(\mbx_\text{lift})>0$.

	\paragraph{(a) At the landing point, the saltation matrix is $S=I-\mbn \mbn^\intercal$, where $I$ is the identity matrix.}
	The saltation matrix at a transition from the interior to a sliding motion along a hard boundary is given in (\citet{bernardo2008}, Example 2.14, p.~111) as 
	\begin{eqnarray}\label{eq:salt-filippov}
	S=I+\frac{(\mbF^{\rm slide}-\mbF^{\rm int})\mbn^\intercal}{\mbn^\intercal \mbF^{\rm int}},
	\end{eqnarray}
	provided the trajectory approaches the hard boundary transversally.
	
	It follows 
	from the definition of the sliding vector field $F^{\rm slide}$ given by \eqref{eq:slidingVF} that 
	\[
	S=I - \mbn \mbn^\intercal,
	\]
	as claimed.

	\paragraph{(b) At the liftoff point, the saltation matrix is $S=I$.}
	We adapt the argument in (\cite{bernardo2008}, \S 2.5) to our hard boundary/liftoff construction.
	The essential difference is that the trajectory is not transverse to the hard boundary at the liftoff point, indeed $\mbn^\intercal \mbF =0$ at $\mbx_\text{lift}$, so eq.~\eqref{eq:salt-filippov} does not give a well defined saltation matrix.
	However, by replacing the  vector $\mbn$ normal to the hard boundary with the vector $\ell$ normal to the liftoff boundary, we recover an equation analogous to \eqref{eq:salt-filippov}, as we will show.
	Since $\mbF^\text{slide}=\mbF^\text{int}$ at the liftoff point, we conclude that the saltation matrix at the liftoff point reduces to the identity matrix.
	
	Let $\Phi_\text{I}$ and $\Phi_\text{II}$ denote the flow operators on the sliding region and in the domain complementary to the sliding region, respectively.
	That is, $\Phi_\text{I}(\mbx,t)$ takes initial point $\mbx\in\mathcal{R}^\text{slide}$ at time zero to $\Phi_\text{I}(\mbx,t)$ at time $0\le t \le \mathcal{T}(\mbx)$.
	So $\Phi_\text{I}$ is restricted to act for  times up to the time $\mathcal{T}(\mbx)$ at which the trajectory starting at $\mbx$ reaches the liftoff point, $\Phi_\text{I}(\mbx,\mathcal{T}(\mbx))\in\mathcal{L}$.
	Such a trajectory necessarily has initial condition $\mbx=(w_1,\ldots,w_{n-1},0)$ satisfying $w_{n-1}<0$, by our coordinatization.
	Let $\mbx_\text{a}\in \mathcal{R}^\text{slide}$ 
	be a point on the periodic limit cycle solution, 
	so that $\Phi_\text{I}(\mbx_\text{a},\mathcal{T}(\mbx_\text{a}))=\mbx_\text{lift}$.
	Write $\tau=\mathcal{T}(\mbx_\text{a})$ for the time it takes for the trajectory to reach the liftoff point after passing location $\mbx_\text{a}$.
	We require a first-order accurate estimate of the effect of the boundary on the displacement between the unperturbed trajectory and a nearby trajectory.
	If we make a small (size $\epsilon$) perturbation into the domain interior, away from the constraint surface, the normal component of the perturbed trajectory will return to zero within a time interval of $O(\epsilon)$ duration, before the two trajectories reach the liftoff boundary.
	Therefore we need only consider perturbations tangent to the constraint surface.
	
\begin{figure}[!t]
		\begin{center}
			\setlength{\unitlength}{1mm}
			\begin{picture}(160,86)(0,0)
			
			\put(0,43){\includegraphics[width=80mm]{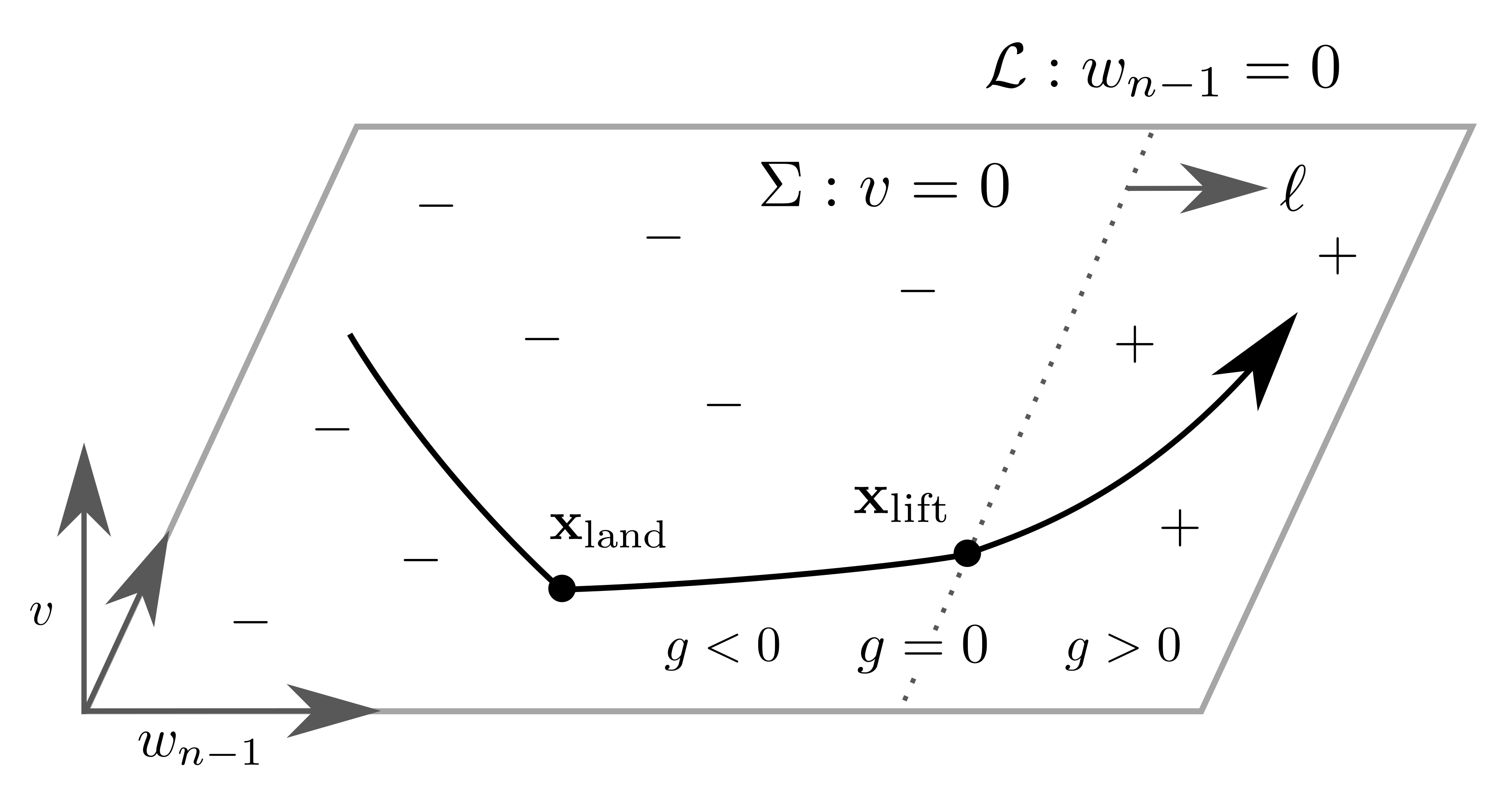}}
			\put(80,43){\includegraphics[width=80mm]{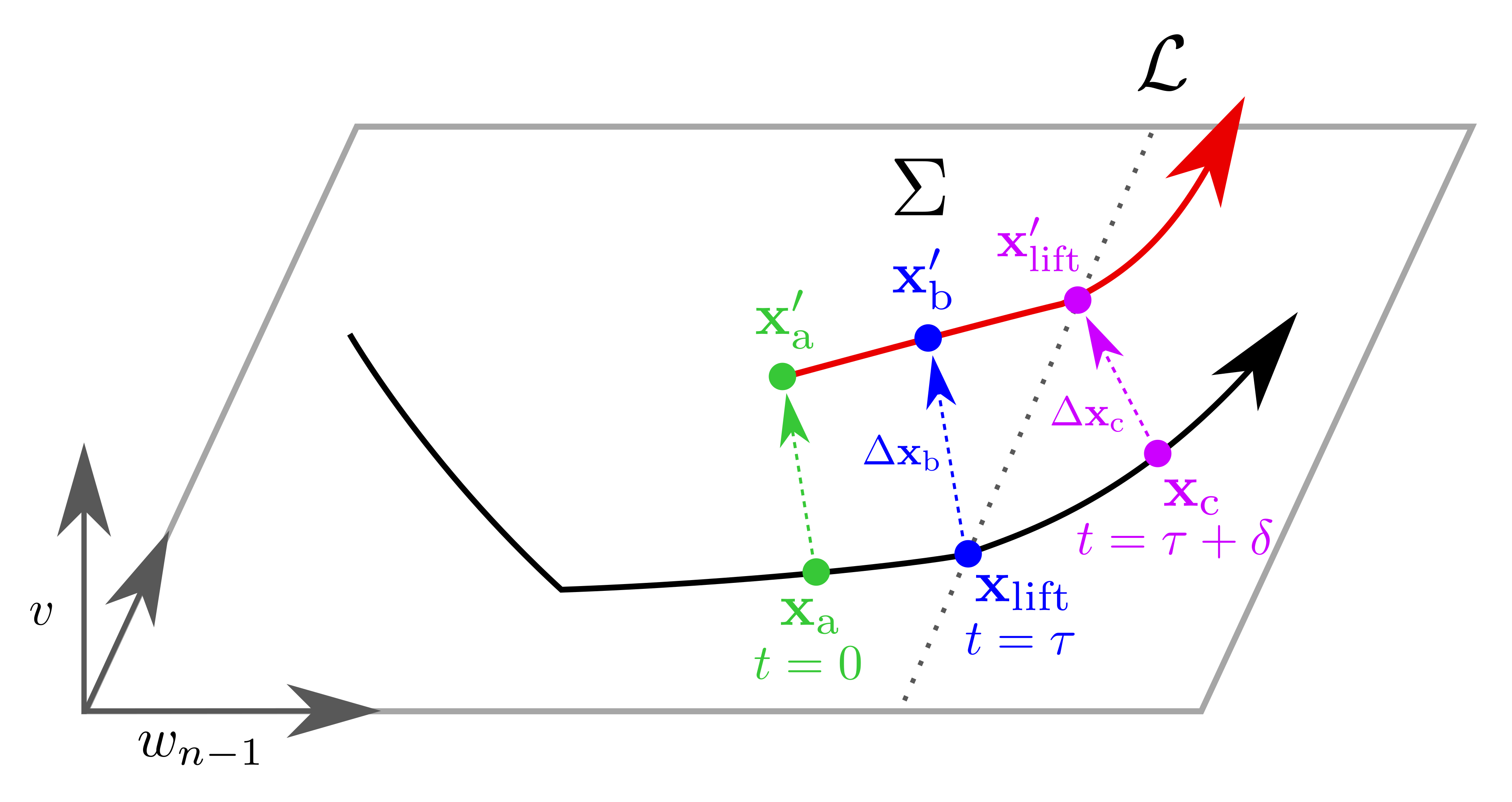}}
			\put(0,0){\includegraphics[width=80mm]{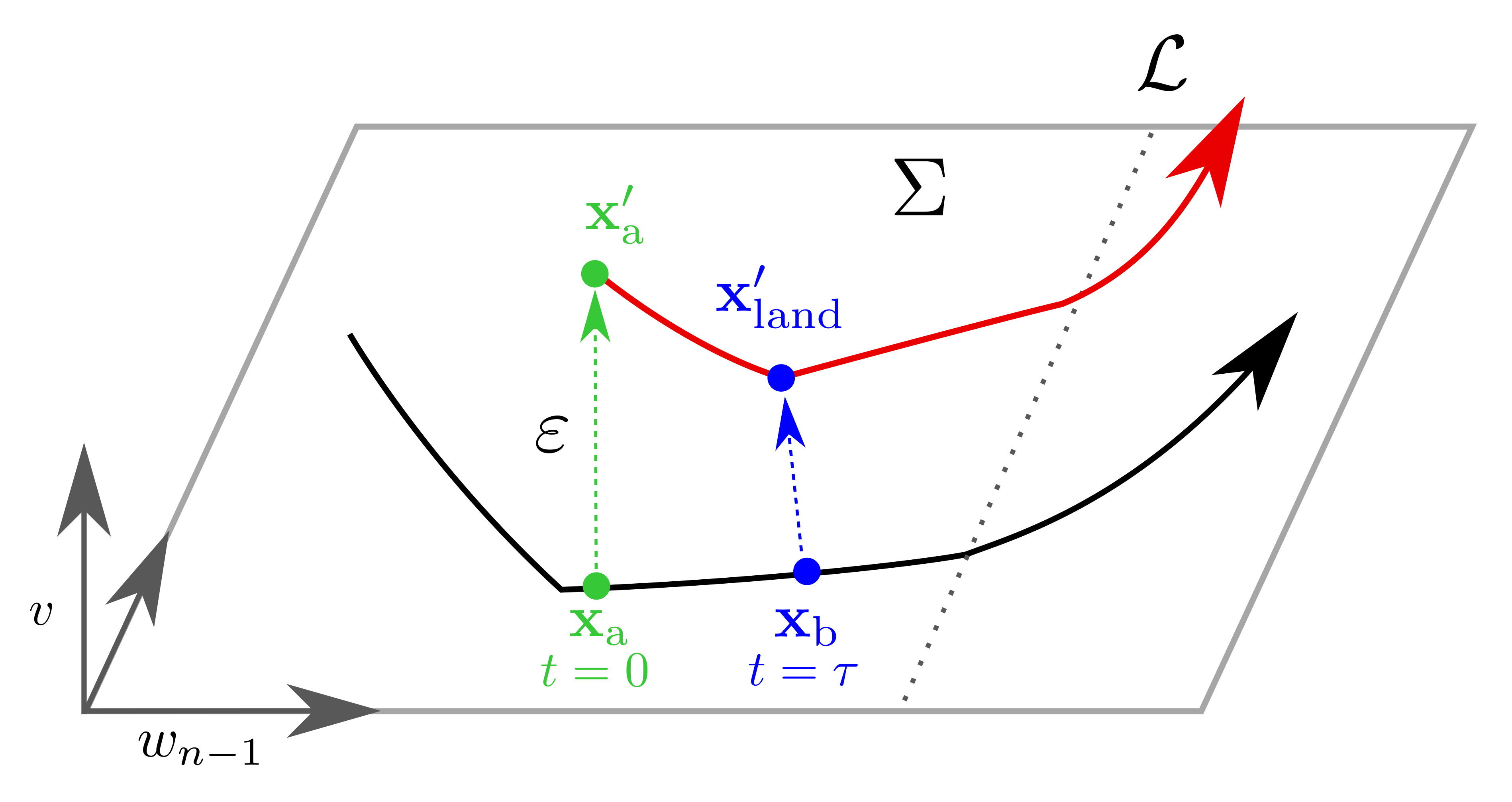}}
			\put(80,0){\includegraphics[width=80mm]{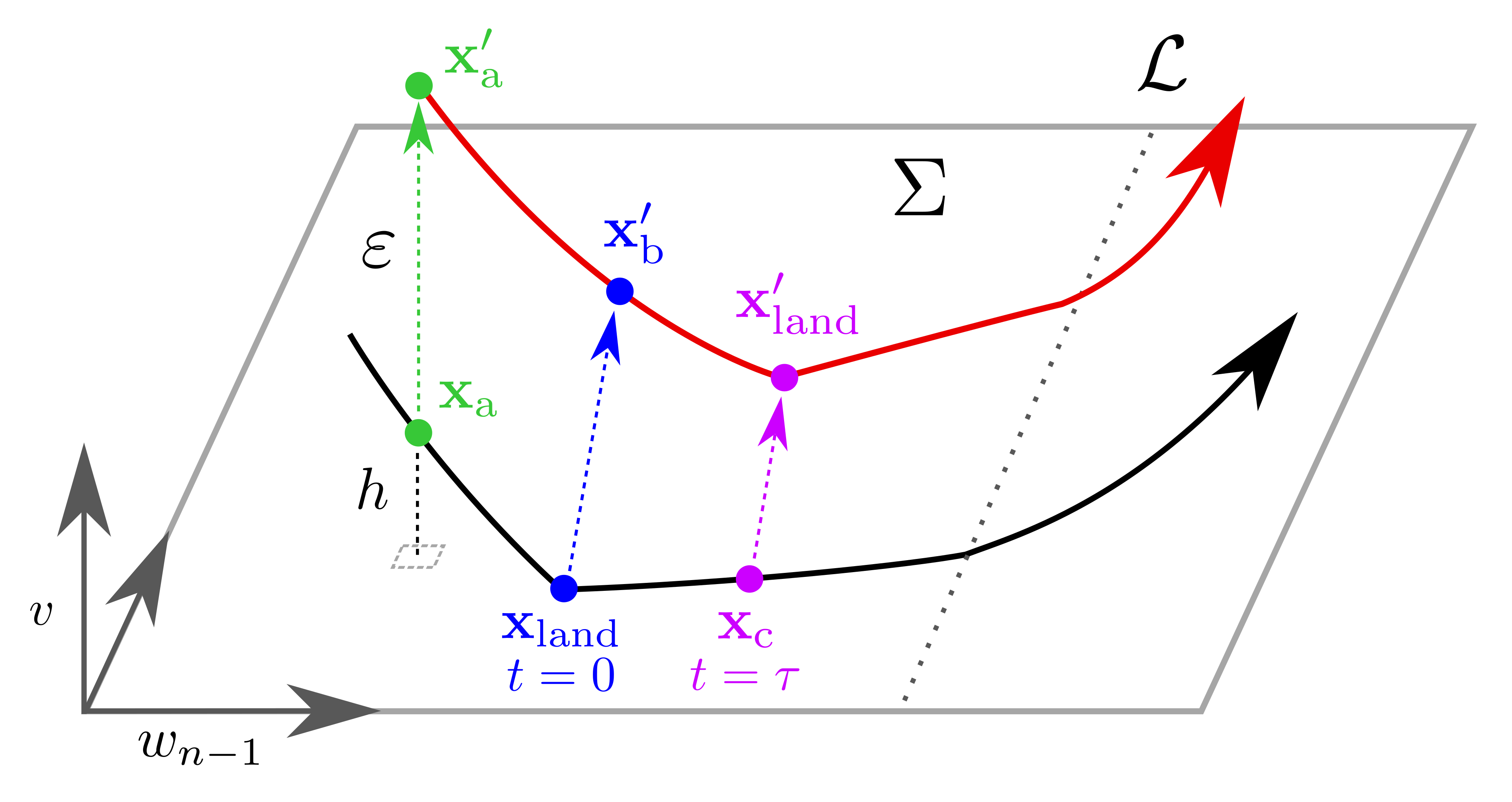}}
			
			\put(2,81){\bf{\large{(A)}}}
			\put(82,81){\bf{\large{(B)}}}
			\put(2,38){\bf{\large{(C)}}}
			\put(82,38){\bf{\large{(D)}}}
			\end{picture}
		\end{center}
		\caption{\label{fig:proof-parts-b-c-d} Unperturbed trajectory (black curve) and a perturbed trajectory (red curve) near the hard boundary $\Sigma$ (horizontal plane) in the $(\mbw, v)$ phase space.
			Dashed line: intersection of liftoff boundary $\mathcal{L}$ and $\Sigma$.
			(A) Trajectory moves downward towards the sliding region (the area in $\Sigma$ where $g<0$), hits $\Sigma$ at the landing point $\mbx_{\rm land}$, and exits $\Sigma$ at the liftoff point $\mbx_{\rm lift}$.
			(B) Construction for the proof of part (b).
			An instantaneous perturbation tangent to $\Sigma$ is made to the point $\mbx_a$ at $t=0$, pushing it to a point $\mbx_a'\in\Sigma$.
			The trajectory starting at $\mbx_a$ (resp., $\mbx_a'$) reaches the liftoff point $\mbx_{\rm lift}$ (resp., $\mbx_b'$) after time $\tau$, and reaches $\mbx_c$ (resp., $\mbx_{\rm lift}'$) after additional time $\delta$.
			The displacements $\Delta\mbx_b=\mbx_b'-\mbx_{\rm lift}$ and $\Delta\mbx_c=\mbx_{\rm lift}'-\mbx_{c}$ differ by an amount captured, to linear order, by the saltation matrix.
			(C) Construction for the proof of part (c).
			An instantaneous perturbation with size $\varepsilon$ in the positive $v$-direction (green arrow) is made to the point $\mbx_a\in \Sigma$, pushing it off the boundary to an interior point $\mbx_a'$.
			After time $\tau$, the trajectory starting at $\mbx_a'$ (resp., $\mbx_a$) reaches a landing point $\mbx_{\rm land}'$ (resp., $\mbx_b$).
			(D) The same perturbation (green arrow) as in panel (C) is applied to the point $\mbx_a$ located at a distance of $h$ above $\Sigma$, pushing it to a point $\mbx_a'$.
			The trajectory starting at $\mbx_a$ lands on $\Sigma$ at $\mbx_{\rm land}$.
			After the same amount of time, the perturbed trajectory starting at $\mbx_a'$ reaches $\mbx_b'$.
			After additional time $\tau$, the two trajectories reach $\mbx_c$ and $\mbx_{\rm land}'$, respectively.
		}
	\end{figure}
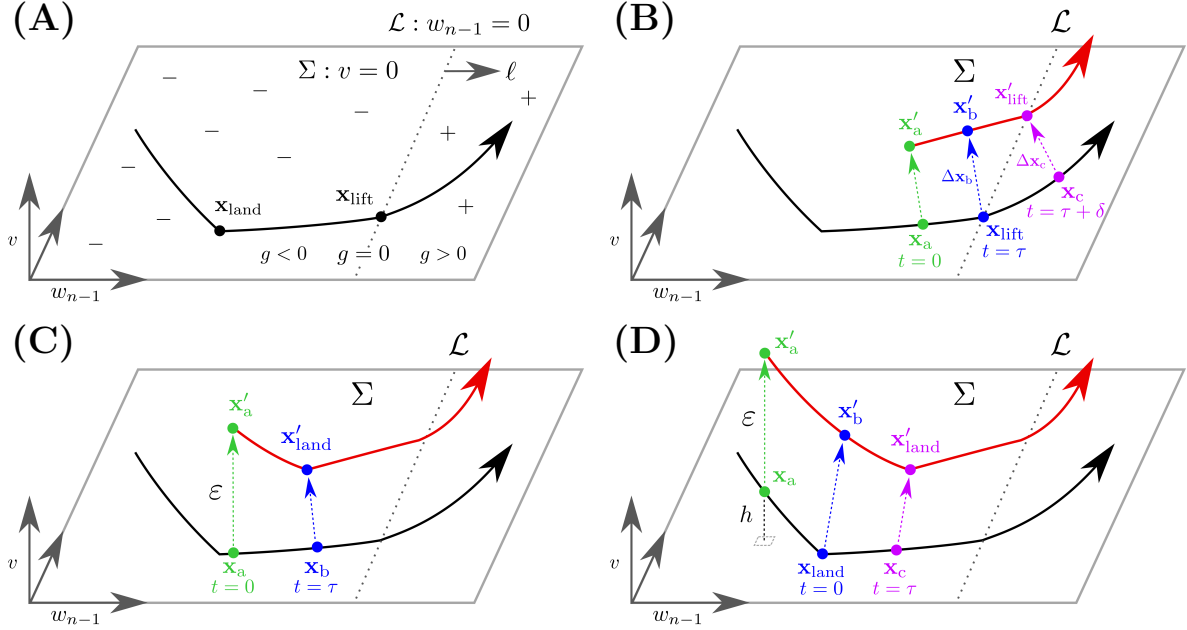

	Let $\mbx_\text{a}'\in\mathcal{R}^\text{slide}$ denote a point near $\mbx_\text{a}$, and suppose it takes time $\mathcal{T}(\mbx_\text{a}')=\tau+\delta$ for the trajectory through $\mbx_\text{a}'$ to liftoff, at some point $\mbx_\text{lift}'\in\mathcal{L}$.
	There are two cases to consider: either $\delta\ge 0$ or else $\delta \le 0$.
	The two cases are handled similarly; we focus on the first for brevity.
	In case $\delta>0$, the original trajectory arrives at the liftoff boundary before the perturbed trajectory, and the point $\mbx_\text{b}'=\Phi_\text{I}(\mbx_\text{a}',\tau)\in\mathcal{R}^\text{slide}$.
	We write $\mbx_\text{b}'=\mbx_\text{lift}+\Delta \mbx_\text{b}$ (see Fig.~\ref{fig:proof-parts-b-c-d}B) 
	and expand the flow operator as follows:
	\begin{align}\nonumber 
	\Phi_\text{I}(\mbx_\text{b}',\delta)=&\mbx_\text{b}'+\delta\, \mbF^\text{slide}(\mbx_\text{b}')+\frac{\delta^2}{2}\left(\nabla^\text{slide}\mbF^\text{slide}(\mbx_\text{b}') \right)\cdot \mbF^\text{slide}(\mbx_\text{b}')+O(\delta^3)\\
	\label{eq:PhiI_expand1}
	=&\mbx_\text{lift}+\Delta\mbx_\text{b} +\delta\, \mbF^\text{slide}(\mbx_\text{lift})+\delta\,\left( \nabla^\text{slide}\mbF^\text{slide}(\mbx_\text{lift} )\right)\cdot\Delta\mbx_\text{b}\\ \nonumber 
	&+\frac{\delta^2}{2}\left(\nabla^\text{slide}\mbF^\text{slide}(\mbx_\text{b}') \right)\cdot \mbF^\text{slide}(\mbx_\text{b}')+O(3),
	\end{align}
	where $\nabla^\text{slide}$ is the gradient operator restricted to $\mbx=(x_1,\ldots,x_{n-1})$.
	The Taylor expansion in \eqref{eq:PhiI_expand1} is justified in a neighborhood of $\mbx_\text{b}'$ contained in the sliding region of the hard boundary.
	The transversality of the intersection of the reference trajectory with $\mathcal{L}$ (that is, $\mbF_{n-1}(\mbx_\text{lift})>0$) means that $\delta$ and $|\Delta\mbx_\text{b}|$ will be of the same order.
	We write $O(n)$ to denote terms of order $\left(|\Delta\mbx_\text{b}|^p \delta^{n-p}\right)$ for $0\le p \le n$.
	
	Next we estimate $\delta$ and the location $\mbx_\text{lift}'$ at which the perturbed trajectory crosses $\mathcal{L}$.
	To first order,
	\begin{align} 
	\ell^\intercal\mbx_\text{b}'&=\ell^\intercal\mbF^\text{slide}(\mbx_\text{b}')\,\delta\\
	\ell^\intercal(\mbx_\text{lift}+\Delta\mbx_\text{b})&=\ell^\intercal\left(\mbF^\text{slide}(\mbx_\text{lift}+\Delta\mbx_\text{b})\right)\delta \label{eq:ineedabettername}\\
	\ell^\intercal\Delta\mbx_\text{b} &=\ell^\intercal\left( \mbF^\text{slide}(\mbx_\text{lift})+\left(\nabla^\text{slide}\mbF^\text{slide}(\mbx_\text{lift})\right)\cdot\Delta\mbx_\text{b} \right)\delta\\
	\nonumber
	&=\ell^\intercal\mbF^\text{slide}(\mbx_\text{lift})\delta+O(2)\\
	\delta &= \frac{\ell^\intercal\Delta\mbx_\text{b}}{\ell^\intercal\mbF^\text{slide}(\mbx_\text{lift})}+O(2).
	\end{align}
	Combining this result with \eqref{eq:PhiI_expand1}, the perturbed trajectory's liftoff location is \begin{equation}
	\label{eq:perturbed_liftoff_location}
	\mbx_\text{lift}'=\mbx_\text{lift}+\Delta\mbx_\text{b}+\mbF^\text{slide}(\mbx_\text{lift})\delta + O(2).
	\end{equation}
	
	Meanwhile, as the perturbed trajectory proceeds to $\mathcal{L}$, during a time interval of duration $\delta$, the unperturbed trajectory has reentered the interior and evolves according to $\Phi_\text{II}$, the flow defined for all initial conditions \emph{not} within the sliding region.
	At a time $\delta$  after reaching $\mathcal{L}$, the unperturbed trajectory is located, to first order, at a point 
	\begin{equation}
	\label{eq:delta-after-liftoff}
	\mbx_\text{c}=\mbx_\text{lift}+\mbF^\text{int}(\mbx_\text{lift})\,\delta+O(2).
	\end{equation}
	Thus,  combining \eqref{eq:perturbed_liftoff_location} and \eqref{eq:delta-after-liftoff} the displacement between the two trajectories immediately following liftoff  of  the perturbed trajectory, $\Delta\mbx_\text{c}=\mbx_\text{lift}'-\mbx_\text{c}$, is given (to first order) by
	\begin{align*}
	\Delta\mbx_\text{c}&=\mbx_\text{lift}'-\mbx_\text{c}\\
	&=\mbx_\text{lift}+\Delta\mbx_\text{b}+\mbF^\text{slide}(\mbx_\text{lift})\delta-\left( \mbx_\text{lift}+\mbF^\text{int}(\mbx_\text{lift})\,\delta \right)\\
	&=\Delta\mbx_\text{b} + \left( \mbF^\text{slide}(\mbx_\text{lift})-\mbF^\text{int}(\mbx_\text{lift}) \right)\delta\\
	&=\Delta\mbx_\text{b} + \frac{\left( \mbF^\text{slide}(\mbx_\text{lift})-\mbF^\text{int}(\mbx_\text{lift}) \right)\ell^\intercal\Delta\mbx_\text{b}}{\ell^\intercal\mbF^\text{slide}(\mbx_\text{lift})}\\
	&=S_\text{lift}\Delta\mbx_\text{b} + O(2).
	\end{align*}
	Therefore, the saltation matrix at the liftoff point is
	\begin{equation}
	S_\text{lift}=I+\frac{\left( \mbF^\text{slide}(\mbx_\text{lift})-\mbF^\text{int}(\mbx_\text{lift}) \right)\ell^\intercal}{\ell^\intercal\mbF^\text{slide}(\mbx_\text{lift})}.
	\label{eq:Slift}        
	\end{equation}
	
	We take the vector field on the sliding region to be the projection of the vector field defined for the interior onto the boundary surface (cf.~\eqref{eq:slidingVF}).
	Therefore for our construction $\mbF^\text{slide}(\mbx_\text{lift})=\mbF^\text{int}(\mbx_\text{lift})$, and hence $S_\text{lift}=I,$ as claimed.
	We note that equation \eqref{eq:Slift} will hold for more general constructions as well.
	This concludes the proof of part (b).

	In parts (c) and (d) of the proof, our goal is to show the normal component of the iPRC is zero along the sliding region on $\Sigma$ and is continuous at the landing point.
	To this end, we compute the normal component of the iPRC using its definition \eqref{eq:definition_of_iPRC}, which in $(\mbw,v)$ coordinates takes the form
	\begin{equation}\label{eq:z-norm}\z_v := \z\cdot \mbn = \lim_{\varepsilon\to 0}\frac{\phi(\mbx+\varepsilon \mbn) -\phi(\mbx)}{\varepsilon},
	\end{equation}
	where $\phi(\mbx)$ denotes the asymptotic phase at point $\mbx$ on the limit cycle.
	That is, we apply a small instantaneous perturbation to the limit cycle, either while it is sliding along $\Sigma$ (part c) or else just before landing (part d), in the $\mbn$ direction, and estimate the phase difference between the perturbed and unperturbed trajectories (cf.~Fig.~\ref{fig:proof-parts-b-c-d}).
	
	\paragraph{(c) Along the sliding region, the component of $\z$ normal to $\Sigma$ is zero.}
	By \eqref{eq:z-norm} the normal component of the iPRC for a point on the sliding component of the trajectory, denoted by $\mbx_\text{a} = (w_\text{a},0)$ is given by 
	\begin{equation}\label{eq:iprc-zero-normcomponent}
	\z_v(\mbx_\text{a})=\lim_{\varepsilon\to0}\frac{\phi(w_\text{a},\varepsilon)-\phi(w_\text{a},0)}{\varepsilon}.
	\end{equation}
	By $\mbx_\text{a}'=(w_\text{a},\varepsilon)$ we denote a  point that is located at a distance of $\varepsilon$ above $\mbx_\text{a}$.
	Our goal is to show $\z_v(\mbx_\text{a})=0$.
	
	The perturbed trajectory from $\mbx_\text{a}'$ is governed by the interior flow $\Phi_\text{II}$ until it reaches the sliding region at a point $\mbx_\text{b}'\in\Sigma$, after some time $\tau$.
	Meanwhile the unperturbed trajectory from $\mbx_\text{a}$ is governed by the sliding flow $\Phi_\text{I}$ until it crosses the liftoff point at $\mathcal{L}$ (Fig.~\ref{fig:proof-parts-b-c-d}, dotted line).
	
	To first order in $\varepsilon$, the time for the perturbed trajectory $\mbx'(t)$ to return to the constraint surface is
	\begin{align}\nonumber
	\tau(\varepsilon)&
	=-\frac{\varepsilon}{g^\text{int}(\mbw_\text{a},\varepsilon)}+O(\varepsilon^2)
	=-\frac{\varepsilon}{g^\text{int}(\mbw_\text{a},0)+\varepsilon D_v g^\text{int}(\mbw_\text{a},0)+O(\varepsilon^2)}+O(\varepsilon^2)\\
	&=-\frac{\varepsilon}{g^\text{int}(\mbw_\text{a},0)}+O(\varepsilon^2),\text{ as }\varepsilon\to 0.
	\end{align}
	Because $\mbx_\text{a}=(\mbw_\text{a},0)$ is in the sliding region, $g^\text{int}(\mbw_\text{a},0)<0$; we conclude that $\tau$ and $\varepsilon$ are of the same order.
	We use $(p)$ to denote terms of order $p$ in $\varepsilon$ or $\tau$.
	
	At time $\tau$ following the perturbation, the location of the perturbed trajectory is 
	\begin{align}
	\mbx_\text{b}'&=\Phi_\text{II}(\mbx_\text{a}',\tau)\\ \nonumber
	&=\mbx_\text{a}'+\tau \mbF^\text{int}(\mbx_\text{a}')+O(2)\\
	&=\mbx_\text{a}+\varepsilon\mbn+\tau\left( \mbF^\text{int}(\mbx_\text{a})+\varepsilon \mbn^\intercal D\mbF^\text{int}(\mbx_\text{a})\right)+O(2) \nonumber\\ \nonumber
	&=\mbx_\text{a}+(0,\ldots,0,\varepsilon)-\frac\varepsilon{g^
		\text{int}(\mbx_\text{a})}(f_1^\text{int}(\mbx_\text{a}),\ldots,f_{n-1}^\text{int}(\mbx_\text{a}),g^\text{int}(\mbx_\text{a}))+O(2)\\ \nonumber
	&=\mbx_\text{a}-\frac\varepsilon{g^
		\text{int}(\mbx_\text{a})}(f_1^\text{int}(\mbx_\text{a}),\ldots,f_{n-1}^\text{int}(\mbx_\text{a}),0)+O(2).
	\end{align}
	Simultaneously, the location of the unperturbed trajectory is
	\begin{align}
	\mbx_\text{b}&=\Phi_\text{I}(\mbx_\text{a},\tau)\\ \nonumber
	&=\mbx_\text{a}+\tau \mbF^\text{slide}(\mbx_\text{a})+O(2)\\ \nonumber
	&=\mbx_\text{a}-\frac{\varepsilon}{g^\text{int}(\mbx_\text{a})}(f_1^\text{int}(\mbx_\text{a}),\ldots,f_{n-1}^\text{int}(\mbx_\text{a}),0) +O(2),
	\end{align}
	since for $\mbx\in\Sigma$, we have $f^\text{slide}(\mbx)=f^\text{int}(\mbx)$ by construction.
	Comparing the difference in location of the two trajectories at time $\tau$ after the perturbation, we see that
	\begin{equation}
	||\mbx_b'-\mbx_b||=O(\varepsilon^2).
	\end{equation} 
	By assumption, the asymptotic phase function $\phi(\mbx)$ is $C^1$ with respect to displacements tangent to the constraint surface.
	Since both $\mbx_\text{b}$ and $\mbx_\text{b}'$ are on this surface, $\mbn^\intercal(\mbx_\text{b}'-\mbx_\text{b})=0,$ and $\phi(\mbx_\text{b}')=\phi(\mbx_\text{b})+O(\varepsilon^2)$.
	Therefore $\z_v(\mbx_\text{a})=0$ for points $\mbx_\text{a}$ on the sliding component of the limit cycle.
	This completes the proof of part (c).

	\paragraph{(d) The normal component of $\z$ is continuous at the landing point.}
	In order to show that the normal component of the iPRC ($\z_v$) is continuous at the landing point, we prove that $\z_v$ has a well-defined limit at the landing point and moreover, this limit equals $0$ which is the value of $\z_v$ at the landing point as proved in (c).
	To this end, consider a point on the limit cycle shortly ahead of the landing point,  $\mbx_
	\text{a}=(\mbw_\text{a},h)$ with $0<h\ll 1$ fixed,  (cf.~Fig.~\ref{fig:proof-parts-b-c-d}D).
	By \eqref{eq:z-norm}
	\begin{eqnarray}\label{eq:prc-d}
	\Eqn{\z_v(\mbx_a)=\lim_{\varepsilon\to 0}\frac{\phi(\mbw_\text{a},h+\varepsilon)-\phi(\mbw_\text{a},h)}{\varepsilon}.}
	\end{eqnarray}
	Our goal is to show $\lim_{h\to 0}\z_v(\mbx_a) = \z_v(\mbx_\text{land}) = 0$.
	
	We consider the case $\varepsilon>0$; the treatment for $\varepsilon<0$ is similar.
	For $\varepsilon>0$, when the unperturbed trajectory arrives at the constraint surface (at landing point $\mbx_\text{land}$), the perturbed trajectory is at a point $\mbx_\text{b}'$ that is still in the interior of the domain.
	Denote the unperturbed landing time $t=0$; denote the time of flight from initial point $\mbx_\text{a}$ to $\mbx_\text{land}$ by $s$.
	Through an estimate similar to that in part (c), to first order in $h$, we have
	\begin{equation}
	s(h)=-\frac{h}{g^\text{int}(\mbx_\text{land})}+O(h^2).
	\end{equation}
	
	Between $t=-s$ and $t=0$, the displacement between the perturbed trajectory ($\mbx'(t)$) and the unperturbed trajectory ($\mbx(t)$) satisfies 
	\begin{align}
	\frac{d(\mbx'-\mbx)}{dt}=D\mbF^\text{int}(\mbx(t))\cdot(\mbx'-\mbx)+O(\epsilon^2),
	\end{align}
	with initial condition $\mbx'(-s)-\mbx(-s)=\varepsilon\mbn$.
	Because the interior vector field is presumed $C^1$, for $h,s\ll 1$ we have 
	\begin{subequations}
		\begin{align*}
		\mbx'_\text{b}-\mbx_\text{land}&= \mbx'_\text{a}-\mbx_\text{a}+s\left( \varepsilon D_v \mbF^\text{int}(\mbx_\text{land}) + O(\epsilon^2)\right) + O(s^2)\\
		&= \varepsilon \mbn -h\left( \varepsilon \frac{D_v \mbF^\text{int}(\mbx_\text{land})}{g^\text{int}(\mbx_\text{land})} + O(\epsilon^2)\right)+O(h^2) \\
		&= (0,\cdots,0, \varepsilon) -\frac{h\varepsilon}{g^\text{int}(\mbx_\text{land})} (f_{1,v}^\text{int}(\mbx_\text{land}),\ldots,f_{n-1,v}^\text{int}(\mbx_\text{land}),g^\text{int}_v(\mbx_\text{land})) + O(2)\\
		&= \left(-h\varepsilon \frac{\mathbf{f}_v^\text{int}(\mbx_\text{land})}{g^\text{int}(\mbx_\text{land})}, \varepsilon-h\varepsilon\frac{g^\text{int}_v(\mbx_\text{land})}{g^\text{int}(\mbx_\text{land})}\right) + O(2).
		\end{align*}
	\end{subequations}
	Here $\mathbf{f}_v^\text{int} = (f_{1,v}^\text{int},\ldots,f_{n-1,v}^\text{int})$, where $f^{\rm int}_{k,v}$ denotes $\partial f^{\rm int}_k/\partial v$, and $O(2)$ denotes terms of order 2 in $\varepsilon$ or $h$ as in (c).
	In the rest of this proof, we drop the dependence of the functions on $\mbx_{\rm land}$ for simplicity.
	
	Since $\mbx_\text{land}$ is in the sliding region, it follows that $\mbx_b'$ is $\varepsilon-h\varepsilon\frac{g^\text{int}_v}{g^\text{int}} + O(2)$ above the sliding region.
	Through a similar estimation as in part (c), to first order in $\varepsilon$ and $h$, the time for the perturbed trajectory to arrive at the sliding region is 
	\[
	\tau(h,\varepsilon) = \frac{\varepsilon-h\varepsilon\frac{g^\text{int}_v}{g^\text{int}}}{-g^\text{int}}+O(2) = -\frac{\varepsilon}{g^\text{int}} + h\varepsilon \frac{g_v^\text{int}}{(g^\text{int})^2} +O(2).
	\]
	
	At time $\tau$, the location of the perturbed trajectory is
	\begin{equation}\label{eq:D:xland'}
	\begin{array}{rcl}
	\mbx_\text{land}' &=& \Phi_\text{II}(\mbx_b',\tau)\\
	&=& \mbx_b'+ \tau \mbF^\text{int}(\mbx_b') + O(2)\\
	&=& \mbx_\text{land} + \left(-h\varepsilon \frac{\mathbf{f}_v^\text{int}}{g^\text{int}}, \varepsilon-h\varepsilon\frac{g^\text{int}_v}{g^\text{int}}\right) +
	\tau \mbF^\text{int}(\mbx_\text{land})+O(2)\\
	&=& \mbx_\text{land} + \left(-h\varepsilon \frac{\mathbf{f}_v^\text{int}}{g^\text{int}}, \varepsilon-h\varepsilon\frac{g^\text{int}_v}{g^\text{int}}\right) +
	\left(-\frac{\varepsilon}{g^\text{int}} + h\varepsilon \frac{g_v^\text{int}}{(g^\text{int})^2} \right) (\mathbf{f}^\text{int}, g^\text{int})+O(2)\\
	&=& \mbx_\text{land} + \left(-h\varepsilon \frac{\mathbf{f}_v^\text{int}}{g^\text{int}}, 0) +
	(-\frac{\varepsilon}{g^\text{int}} + h\varepsilon \frac{g_v^\text{int}}{(g^\text{int})^2} \right) (\mathbf{f}^\text{int}, 0)+O(2).
	\end{array}
	\end{equation}
	
	Simultaneously, the location of the unperturbed trajectory is 
	\begin{equation}\label{eq:D:xc}
	\begin{array}{rcl}
	\mbx_\text{c}&=& \Phi_\text{I}(\mbx_\text{land},\tau)\\
	&=& \mbx_{\rm land}+\tau \mbF^\text{slide}(\mbx_\text{land}) +O(2)\\
	&=& \mbx_\text{land} + \left(-\frac{\varepsilon}{g^\text{int}} + h\varepsilon \frac{g_v^\text{int}}{(g^\text{int})^2} \right)(\mathbf{f}^\text{int}, 0) + O(2).
	\end{array}
	\end{equation}
	
	Comparing the difference between \eqref{eq:D:xland'} and \eqref{eq:D:xc}, we see that
	\begin{equation}
	\norm{\mbx_{\rm land}' - \mbx_c} = O(h\varepsilon).
	\end{equation}
	Recall that the asymptotic phase is assumed to be $C^1$, with respect to displacements tangent to $\Sigma$.
	Since $\mbx_{\rm land}'$ and $\mbx_c$ are on $\Sigma$, 
	it follows that 
	\[
	\phi(\mbx_{\rm land}') - \phi(\mbx_c) = O(h\varepsilon).
	\]
	Therefore, by \eqref{eq:prc-d},
	\[
	\z_v(\mbx_a) = \lim_{\varepsilon\to 0} \frac{\phi(\mbx_a')-\phi(\mbx_a)}{\varepsilon} = \lim_{\varepsilon\to 0} \frac{\phi(\mbx_{\rm land}')-\phi(\mbx_c)}{\varepsilon}=O(h)
	\]
	Consequently,
	\[
	\lim_{h\to 0} \z_v(\mbx_a) = 0
	\]
	as required.
	This completes the proof of part (d).
	
	\paragraph{(e) The tangential components of $\z$ are continuous at both landing and liftoff points.}
	We denote the tangential components of the iPRC by $\z_\mbw$, where $\mbw$ represents vectors in the $n-1$ dimensional tangent space of the hard boundary.
	The $n-1$ dimensional iPRC vector $\z_\mbw$ obeys a restricted  (\textit{i.e.}~reduced dimension) adjoint equation given in terms of $f_\mbw$, the $(n-1)\times (n-1)$ Jacobian derivative of $f$ with respect to the $n-1$ tangential coordinates ($\mbw$), and $g_\mbw$, the $1\times(n-1)$ Jacobian derivative of $g$ with respect to the tangential coordinates, and $\z_v$, the (scalar) component of $\z$ in the normal direction
	\begin{eqnarray}\label{eq:zu-interior}
	\Eqn{\frac{d\z_\mbw}{dt}=-f_\mbw(\mbw,v)^
		\intercal \z_\mbw-g_\mbw(\mbw,v)^\intercal \z_v} 
	\end{eqnarray}
	along the limit cycle in the interior domain.
	On the other hand, along the sliding component of the limit cycle that is restricted to $\{\Sigma: v=0\}$, $\z_u$ satisfies  
	\begin{eqnarray}\label{eq:zu-sliding}
	\Eqn{\frac{d\z_\mbw}{dt}=-f_\mbw(\mbw,0)^\intercal \z_\mbw.}
	\end{eqnarray}
	By part (c), $\z_v$ goes continuously to zero as the trajectory from the interior approaches the landing point.
	Therefore $\z_\mbw$ is continuous at the landing point.
	
	\begin{figure}[htpb]\centering
		\includegraphics[width=4in]{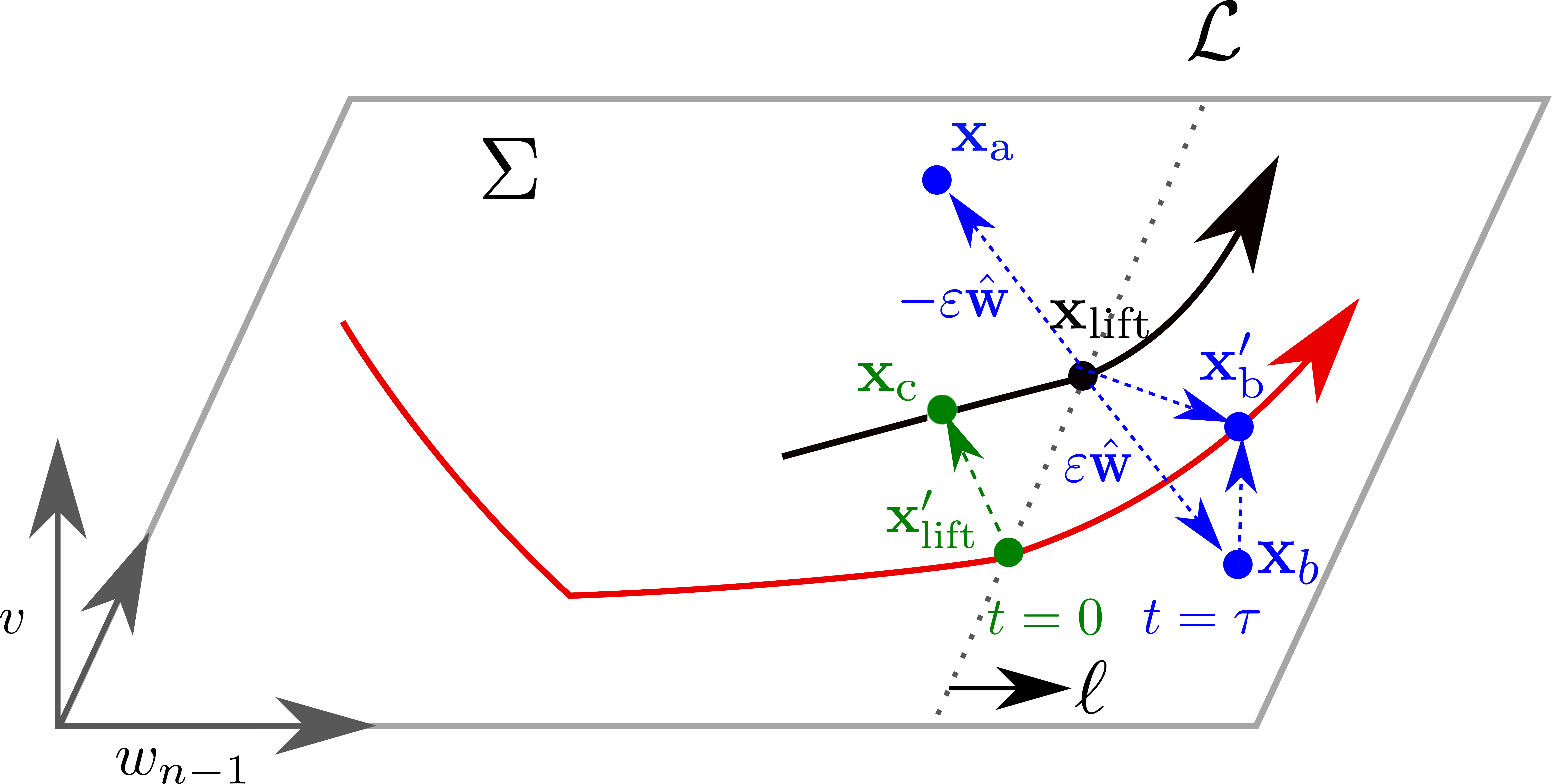}
		\caption{\label{fig:proof-partE-liftoff} Unperturbed trajectory (black) leaves the hard boundary at the liftoff point $\mbx_{\rm lift}$, in the $(\mathbf{w}, v)$ phase space.
			An instantaneous perturbation tangent to $\Sigma$ is made to the liftoff point at $t=\tau$, pushing it to $\mbx_a$ on the sliding region or to $\mbx_b$ that is outside the sliding region.
			The points $\mbx_c$ and $\mbx_{\rm lift}'$ denote the positions of the unperturbed trajectory and the perturbed trajectory at $t=0$. } 
	\end{figure}
	
	Next we prove the continuity of $\z_\mbw$ at the liftoff point $\mbx_{\rm lift}=(\mbw_{\rm lift},0)$.
	Recall that in the coordinates employed, the unit vector tangent to $\Sigma$ and normal to $\mathcal{L}$ at $\mbx_\text{lift}$ is $\ell=(0,\ldots,0,1,0)$ (cf.~Fig.~\ref{fig:proof-partE-liftoff}).
	Fix an arbitrary tangential unit vector $\hat{\mbw}$ oriented away from the sliding region (such that $\ell^\intercal\hat{\mbw}> 0$).
	The left and right limits of $\z_w$ at $\mbx_{\rm lift}$ are given by
	\begin{equation}\label{eq:lift-leftlimit}
	\z_{\hat{\mbw}}^-(\mbx_{\rm lift}) = \lim_{\varepsilon\to 0^+} \frac{\phi(\mbw_{\rm lift}-\varepsilon \hat{\mbw},0)-\phi(\mbw_{\rm lift},0)}{-\varepsilon} 
	\end{equation}
	and
	\begin{equation}\label{eq:lift-rightlimit}
	\z_{\hat{\mbw}}^+(\mbx_{\rm lift}) = \lim_{\varepsilon\to 0} \frac{\phi(\mbw_{\rm lift}+\varepsilon \hat{\mbw},0)-\phi(\mbw_{\rm lift},0^+)}{\varepsilon}.
	\end{equation}
	By $\mbx_a=(\mbw_{\rm lift}-\varepsilon \hat{\mbw},0)$ and $\mbx_b=(\mbw_{\rm lift}+\varepsilon \hat{\mbw},0)$ we denote the two points that are located at a distance of $\varepsilon$ away from $\mbx_{\rm lift}$ along the  $-\hat{\mbw}$ and $\hat{\mbw}$ directions, respectively (cf.~Fig.~\ref{fig:proof-partE-liftoff}).
	We will show that 
	\begin{equation}\label{eq:defined-limit}
	\z_{\hat{\mbw}}^-(\mbx_{\rm lift}) = \z_{\hat{\mbw}}^+(\mbx_{\rm lift}).
	\end{equation}
	The equality of these limits will establish that  $\z_{\mbw}$ is continuous at the liftoff point.
	
	First, we consider $\z_{\hat{\mbw}}^+(\mbx_\text{lift})$.
	Given $\hat{\mbw}$, there exists a unique point $\mbx_{\rm lift}'$ at the liftoff boundary $\mathcal{L}\cap\Sigma$, and a time $\tau>0$, such that the trajectory beginning from $\mbx_{\rm lift}'$ at time $0$ passes directly over $\mbx_b'$ at time $\tau$, in the sense that $\Phi_\text{II}(\mbx_\text{lift}',\tau)=(\mbw_b, h)$, where $\Phi_\text{II}$ is the flow operator in the complement of the sliding region, $h>0$ is the ``height'' of $\mbx_b'$ above $\mbx_b$, and $\mbw_b$ is the coordinate vector along the tangent space of the hard boundary.
	Let $\mbx_{\rm lift} = (\mbw_{\rm lift},0)$ and $\mbx_{\rm lift}' = (\mbw_{\rm lift}',0)$.
	By our construction, $\mbw_b = \mbw_{\rm lift} + \varepsilon\hat{\mbw}$.
	Hence, the location of the perturbed trajectory at time $\tau$ is
	\begin{equation*}
	\begin{array}{rcl}
	(\mbw_b,h)=   (\mbw_{\rm lift}+\varepsilon \hat{\mbw},h) &=&\Phi_{\rm II}(\mbx_{\rm lift}',\tau)\\ &=&  (\mbw_{\rm lift}',0) + (f^{\rm int}(\mbx_{\rm lift}'), 0) \tau + O(\tau^2)\\ &=& (\mbw_{\rm lift}'+f^{\rm int}(\mbx_{\rm lift}')\tau+O(\tau^2), O(\tau^2)).
	\end{array}
	\end{equation*}
	Hence
	\begin{equation}\label{eq:wlift}
	\mbw_{\rm lift}-\mbw_{\rm lift}' = f^{\rm int}(\mbx_{\rm lift}')\tau -\varepsilon \hat{\mbw}+O(\tau^2),
	\end{equation}
	\begin{equation}\label{eq:hlift}
	h = O(\tau^2),
	\end{equation}
	and
	\begin{equation}\label{eq:wb}
	\mbw_b - \mbw_{\rm lift}'= f^{\rm int}(\mbx_{\rm lift}')\tau+O(\tau^2)
	\end{equation}
	
	On the other hand,
	\[
	\varepsilon \hat{\mbw} + (\mbw_{\rm lift}' - \mbw_{b}) = \mbw_{\rm lift}- \mbw_{\rm lift}'.
	\]
	By \eqref{eq:wlift} and \eqref{eq:wb}, the above equation becomes 
	\[
	\varepsilon \hat{\mbw} - f^{\rm int}(\mbx_{\rm lift}')\tau =f^{\rm int}(\mbx_{\rm lift}')\tau-\varepsilon \hat{\mbw} + O(\tau^2).
	\]
	That is,
	\[\varepsilon \hat{\mbw} = f^{\rm int}(\mbx_{\rm lift}')\tau+ O(\tau^2).
	\]
	Taking the inner product of both sides with the unit vector $\ell$ (normal to $\mathcal{L}$), and noting that for sufficiently small $\epsilon$,  $\ell^\intercal f^\text{int}(x_\text{lift}')>0$ (our nondegeneracy condition), we have
	\[
	\tau=\epsilon \frac{\ell^\intercal \hat{\mbw}}{\ell^\intercal f^\text{int}(\mbx_\text{lift}')}+O(\tau^2),
	\]
	and hence $\tau=O(\varepsilon)$.
	Therefore, \eqref{eq:hlift} becomes
	\begin{equation}\label{eq:h-is-O-eps2}
	h=O(\varepsilon^2)    
	\end{equation}
	and hence the phase difference between $\mbx_b'$ and $\mbx_b$ is 
	\begin{equation}\label{eq:phase-diff-b}
	\phi(\mbx_b)-\phi(\mbx_b')=O(\varepsilon^2)
	\end{equation}
	due to the assumption that $\phi$ is Lipschitz continuous.
	
	
	Next we show \eqref{eq:defined-limit} holds using \eqref{eq:lift-leftlimit} and \eqref{eq:lift-rightlimit}.
	Let the unperturbed trajectory pass through $\mbx_\text{lift}$ at time $\tau$, and 
	let $\mbx_c$ be the location of the unperturbed trajectory at time $t=0$ (see Fig.~\ref{fig:proof-partE-liftoff}).
	Let $\Delta\mbx_c=\mbx_{\rm lift}'-\mbx_c$ and  $\Delta\mbx_b=\mbx_b' -\mbx_{\rm lift}$.
	Then by part (b),
	\[
	\Delta\mbx_b-\Delta\mbx_c= O(|\Delta\mbx_b|^2); 
	\]
	since the saltation matrix is equal to the identity matrix at the liftoff boundary.
	Since $\mbx_{\rm lift},\mbx_b,\mbx_b'$ form a right triangle, 
	\[
	|\Delta\mbx_b|^2 = \varepsilon^2+ h^2  = \varepsilon^2+ O(\varepsilon^4),
	\]
	which implies that 
	\begin{equation}\label{eq:DeltaxbDeltaxc}
	\Delta\mbx_b-\Delta\mbx_c= O(\varepsilon^2).
	\end{equation}
	
	Direct computation shows
	\begin{equation}\label{eq:xb-xlift}
	\Eqn{\phi(\mbx_b)-\phi(\mbx_{\rm lift}) &=&(\phi(\mbx_b)-\phi(\mbx_b')) + (\phi(\mbx_b')-\phi(\mbx_{\rm lift}))\\
		&=& (\phi(\mbx_{\rm lift}')-\phi(\mbx_c))+O(\varepsilon^2) \\
		&=&  D_\mbw\phi(\mbx_c)\cdot \Delta\mbx_c+O(\varepsilon^2)\\
		&=& D_\mbw\phi(\mbx_c)\cdot \Delta\mbx_b+O(\varepsilon^2)\\
		&=&  D_\mbw\phi(\mbx_c)\cdot (\mbx_b' - \mbx_{\rm lift})+O(\varepsilon^2)\\
		&=&  D_\mbw\phi(\mbx_c)\cdot (\mbx_b - \mbx_{\rm lift})+O(\varepsilon^2)\\
		&=&  D_\mbw\phi(\mbx_c)\cdot \varepsilon\hat{\mbw}+O(\varepsilon^2)}.
	\end{equation}
	To obtain the second equality, we translate the trajectories backward in time by $\tau$ beginning from $\mbx_b'$ and $\mbx_\text{lift}$, respectively; shifting both trajectories by an equal time interval does not change their phase relationship.
	The $O(\epsilon^2)$ difference arises from \eqref{eq:phase-diff-b}.
	The third equality follows from the assumption that $\phi$ is differentiable with respect to displacements tangent to the sliding region.
	The fourth equality uses \eqref{eq:DeltaxbDeltaxc}; the fifth and seventh follow from the definitions; the sixth uses \eqref{eq:h-is-O-eps2}.
	
	Recall the we assume $\phi$ to have Lipschitz continuous derivatives in the tangential directions at the boundary surface (except possibly at the landing and liftoff points).
	Under this assumption, 
	taking the limit $\varepsilon\to 0^+$ leads to $\mbx_c\to \mbx_{\rm lift}^-$ and hence
	\[\z^+_{\hat{\mbw}}(\mbx_{\rm lift}) =D_\mbw\phi(\mbx_{\rm lift}^-)\cdot \hat{\mbw}\] by \eqref{eq:lift-rightlimit}.
	On the other hand,
	\begin{equation}\label{eq:xa-xlift}
	\Eqn{\phi(\mbx_a)-\phi(\mbx_{\rm lift}) &=&
		D_\mbw\phi(\mbx_a)\cdot (\mbx_a-\mbx_{\rm lift}) +O(\varepsilon^2)\\
		&=& -D_\mbw\phi(\mbx_a)\cdot \varepsilon\hat{\mbw}+O(\varepsilon^2). }
	\end{equation}
	Taking the limit $\varepsilon\to 0+$ results in $\mbx_a\to \mbx_{\rm lift}^-$ and hence \eqref{eq:lift-leftlimit} together with \eqref{eq:xa-xlift}, implies 
	\[
	\z^-_{\hat{\mbw}}(\mbx_{\rm lift}) =D_\mbw\phi(\mbx_{\rm lift}^-)\cdot \hat{\mbw}.
	\] 
	Hence, \eqref{eq:defined-limit} holds.
	\endproof

	\newpage
	\section{Numerical Algorithms} \label{sec:algorithm}
	
	We will now describe how the results presented in \S\ref{sec:smooth-theory} and \S\ref{sec:nonsmooth-theory} can be implemented as numerical algorithms.
	MATLAB code that implements these algorithms for the example system described in \S\ref{sec:toy-model} is available: \url{https://github.com/yangyang-wang/LC_in_square}.
	
	Consider a multiple-zone Filippov system generalized from \eqref{eq:slidingVF},
	\begin{equation}\label{eq:multiple-zone-FP-simp}
	\frac{d\mbx}{dt}=F(\mbx),
	\end{equation}
	that produces a $T_0$-periodic limit cycle solution $\gamma(t)\subset\R^n$.
	Suppose $\gamma(t)$ includes  $k$ sliding components  confined to boundary surfaces denoted as $\Sigma^i\subset\R^{n-1},\,i\in\{1,...,k\}$.
	$\gamma(t)$ exits the $i$-th boundary $\Sigma^i$ at a unique liftoff point $\mbx_{\text{lift}}^i$ given that the nondegeneracy condition \eqref{eq:nondege-1} at $\mbx_{\text{lift}}^i$ is satisfied.
	We denote the  normal vector to $\Sigma^i$ at  liftoff, landing, or boundary crossing points by $n^i$.
	We denote the interior domain by $\RR^{\rm interior}$, which can now consist of multiple subdomains separated by transversal crossing boundaries, and denote the piecewise smooth vector field in $\RR^{\rm interior}$ by $F^{\rm interior}$.
	By \eqref{eq:slidingVF}, the sliding vector field on the sliding region $\RR^{\mathrm{slide}_i}\subset\Sigma^i$ is therefore 
	\begin{equation}\label{eq:multiple-zone-sliding}
	F^{\mathrm{slide}_i}(\mbx)= F^{\rm interior}(\mbx) - (n^i\cdot F^{\rm interior}(\mbx))n^i
	\end{equation} 
	
	Using this notation, the vector field \eqref{eq:multiple-zone-FP-simp} can be written as
	\begin{eqnarray}\label{eq:multiple-zone-FP}
	F(\mbx):=\left\{
	\Eqn{
		F^{\rm interior}(\mbx), & \mbx \in \RR^{\rm interior}&\\
		F^{\mathrm{slide}_i}(\mbx),  & \mbx \in\RR^{\mathrm{slide}_i}&\\
	}
	\right.
	\end{eqnarray}
	and we denote the vector field after a static perturbation by 
	\begin{eqnarray}\label{eq:multiple-zone-PFP}
	F_{\varepsilon}(\mbx):=\left\{
	\Eqn{
		F^{\rm interior}_\varepsilon(\mbx), & \mbx \in \RR^{\rm interior}&\\
		F^{\mathrm{slide}_i}_\varepsilon(\mbx),  & \mbx \in\RR^{\mathrm{slide}_i}&\\
	}
	\right.
	\end{eqnarray}
	where $i\in\{1, ..., k\}$.
	Here we assume that the regions are independent of static perturbation with size $\varepsilon$.
	
	Notice that the computation of the iSRC requires estimating the rescaling factors, for which we need to compute the iPRC or the lTRC depending on whether a global uniform rescaling \eqref{eq:src} or a piecewise uniform rescaling \eqref{eq:src-multiplescales} is needed.
	We hence first present the numerical algorithms for obtaining the iPRC in \S\ref{sec:algorithm-iprc} and the lTRC in \S\ref{sec:algorithm-lTRC}; the algorithm for solving the homogeneous variational equation for the linear shape responses of $\gamma(t)$ to instantaneous perturbations (the variational dynamics $\mbu$) is presented in \S\ref{sec:algorithm-u};
	lastly, in \S\ref{sec:algorithm-x1} we illustrate the algorithms for computing the linear shape responses of $\gamma(t)$ to sustained perturbations (the iSRC $\gamma_1$) with a uniform rescaling factor computed from the iPRC as well as with piecewise uniform rescaling factors computed from the lTRC.
	
	For simplicity, we assume the initial time is $t_0=0$.
	
	\subsection{Algorithm for Calculating the iPRC \texorpdfstring{$\z$}{z} for LCSCs}\label{sec:algorithm-iprc}
	
	It follows from Remark \ref{rem:iprc-liftoff} that the iPRC $\z$ for the LCSCs need to be solved backward in time.
	While there is no discontinuity of $\z$ at a landing point, a time-reversed version of the jump matrix at the liftoff point on the hard boundary $\Sigma^i$, denoted as $\mathcal{J}^i_{\rm lift}$, is given by
	\begin{equation}\label{eq:back-jump-def}
	\mathcal{J}^i_{\rm lift}=I-n^i {n^i}^\intercal,
	\end{equation} 
	where $I$ is the identity matrix.
	$\mathcal{J}^i_{\rm lift}$ updates $\z$ local to the liftoff point as
	\begin{eqnarray}
	\Eqn{
		{\z_{\text{lift}}^{i^-}}=\mathcal{J}^i_{\rm lift} {\z_{\text{lift}}^{i^+}}
	}
	\end{eqnarray}
	where ${\z_{\text{lift}}^{i^-}}$ and ${\z_{\text{lift}}^{i^+}}$ 
	are the iPRC just before and just after the trajectory crosses the liftoff point $x_{\text{lift},i}$ in forwards time.

	We now describe an algorithm for numerically obtaining the complete iPRC $\z$ for $\gamma(t)$, a stable limit cycle with sliding components along hard boundaries and transversal crossing boundaries as described before.
	
	\paragraph{Algorithm for $\z$}
	\begin{enumerate}
		\item Fix an initial condition $\mbx_0=\gamma(0)$ on the limit cycle, and integrate~\eqref{eq:multiple-zone-FP} to compute $\gamma(t)$ over $[0, T_0]$.
		\item Integrate the adjoint equation backward in time by defining $s=T_0-t$ and numerically solve for the fundamental matrix $\Psi(s)$ over one period $0\le s \le T_0$, where $\Psi$ satisfies
		\begin{enumerate}
			\item $\Psi(0)=I$, the identity matrix.
			\item For $s$ such that $\gamma(T_0-s)$ lies in the interior of the domain, 
			$$ \frac{d\Psi}{ds}=A^{\rm interior}(T_0-s)\Psi $$
			where $A^{\rm interior}(t)=\left(DF^{\rm interior} (\gamma(t))\right)^\intercal$ is the transpose of the Jacobian of the interior vector field $F^{\rm interior}$.
			\item For $s$ such that $\gamma(T_0-s)$ lies within a sliding component along boundary $\Sigma^i$, 
			$$ \frac{d\Psi}{ds}= A^i(T_0-s)\Psi $$
			where $A^i(t)=\left(DF^{\mathrm{slide}_i}(\gamma(t))\right)^\intercal$ is the transpose of the Jacobian of the sliding vector field $F^{\mathrm{slide}_i}$, given in \eqref{eq:multiple-zone-sliding}.
			\item At any time $t_p$ when $\gamma$ transversely crosses a switching surface with a normal vector $n_p$, 
			$$ \Psi^- = \mathcal{J} \Psi^+$$
			where $\Psi^-=\lim_{s\to (T_0-t_p)^+}\Psi(s)$ and $\Psi^+=\lim_{s\to (T_0-t_p)^-}\Psi(s)$ are the fundamental matrices just before and just after crossing the surface in forwards time. $\mathcal{J}=S^\intercal$ since $J^\intercal S=I$ as discussed in \S\ref{sec:nonsmooth-review}, where the saltation matrix at any transversal crossing point is
			$$S=I+\frac{(F_p^+-F_p^-)n_p^\intercal }{n_p^\intercal F_p^-}$$
			where $F_p^-, F_p^+$ are the vector fields just before and just after the crossing in forwards time (see \eqref{eq:salt}).
			
			\item At a liftoff point on the $i$-th hard boundary $\Sigma^i$ (in backwards time, a transition from the interior to $\Sigma^i$), update $\Psi$ as
			$$ \Psi^- = \mathcal{J}^i \Psi^+$$
			where $\mathcal{J}^i=I-n^i n^{i \intercal}$ as defined in \eqref{eq:back-jump-def}, and then switch the integration from the full Jacobian $A^{\rm interior}$ to the restricted Jacobian $A^i$.
			\item At a landing point on the $i$-th hard boundary $\Sigma^i$ (in backwards time, a transition from $\Sigma^i$ to the interior) switch integration from the restricted Jacobian $A^i$ to the full Jacobian $A^{\rm interior}$; no other change in $\Psi$ is needed.
		\end{enumerate}
		\item Diagonalize the fundamental matrix at one period $\Psi(T_0)$; it should have a single eigenvector $v$ with unit eigenvalue.
		The initial value for $\z_\text{BW}$ (represented in \emph{backwards time}) at the point $\gamma(T_0)=\gamma(0)=\mbx_0$ is given by 
		$$\z_\text{BW}(0)=\frac{v}{F(\mbx_0)\cdot v}$$
		\item The iPRC in backward time over $s\in [0, T_0]$ is given by $\z_\text{BW}(s)=\Psi(s)\z_{\rm BW}(0)$ and is $T_0$-periodic.
		Equivalently, one may repeat step (2) by replacing $\Psi(s)$ with $\z_\text{BW}(s)$ and replacing the initial condition $\Psi(0)=I$ with $\z_{\rm BW}(0)$ to solve for the complete iPRC.
		\item The iPRC in forward time is then given by $\z(t)=\z_\text{BW}(T_0-t)$ where $t\in[0,T_0]$.
		\item
		The linear shift in period in response to the static perturbation can be calculated by evaluating the integral (see \eqref{eq:T1}) 
		\[
		T_1=-\int_{0}^{T_0} \z^\intercal(t)\frac{\partial F_\varepsilon(\gamma(t))}{\partial \varepsilon}\Big|_{\varepsilon=0}dt
		\]
	\end{enumerate}

	\begin{remark}
		An alternative way (in MATLAB) to do backward integration is reversing the time span in the numerical solver; that is, integrate the adjoint equation over $[T_0, 0]$ to compute $\z(t)$.
	\end{remark}

	\subsection{Algorithm for Calculating the lTRC for LCSCs}\label{sec:algorithm-lTRC}
	
	The lTRC satisfies the same adjoint equation, \eqref{eq:prc}, as the iPRC, and hence exhibits the same jump matrix at each liftoff, landing and boundary crossing point.
	It follows that the algorithm for the iPRC from \S\ref{sec:algorithm-iprc} can mostly carry over to computing the lTRC.
	
	Suppose the domain of $\gamma(t)$ can be divided into $m$ regions $\RR^1, ..., \RR^m$, each distinguished by its own timing sensitivity properties.
	We denote the lTRC in $\RR^j$ by $\eta^j$.
	
	Below we describe the algorithm to compute $\eta^j$ in region $\RR^j$ bounded by the two local timing surfaces $\Sigma^{\rm in}$ and $\Sigma^{\rm out}$.
	Following the notations in \S\ref{sec:smooth-theory}, $t^{\rm in}$ and $t^{\rm out}$ denote the time of entry into and exit out of $\RR^j$, at locations $\mbx^{\rm in}$ and $\mbx^{\rm out}$, respectively.
	The algorithm for computing $\eta^j$ is described as follows.
	
	\paragraph{Algorithm for $\eta^j$}
	\begin{enumerate}
		\item  Compute $\gamma$, the unperturbed limit cycle, and $T_0$, its period, by integrating ~\eqref{eq:multiple-zone-FP}.
		\item Compute $t^{\rm in}, t^{\rm out}$ for region $j$.
		Evaluate  $\mbx^{\rm in}=\gamma(t^{\rm in}),\, \mbx^{\rm out}=\gamma(t^{\rm out})$ and $T_0^{j}=t^{\rm out}-t^{\rm in}$.
		\item Compute the boundary value for $\eta^j$ at the exit point $\mbx^{\rm out}$ (see \eqref{eq:ltrc0}) 
		\[
		\eta^{j}(\mbx^{\rm out})=\frac{-n^{\rm out}}{{n^{\rm out}}^\intercal F(\mbx^{\rm out})}
		\]
		where $n^{\rm out}$ is a normal vector to $\Sigma^{\rm out}$.
		\item Integrate the adjoint equation backward in time by defining $s=T_0-t$ and numerically solve for $\eta^j_{\rm BW}(s)$ (represented in backwards time) over $[T_0-t^{\rm out}, T_0-t^{\rm in}]$.
		$\eta^j_{\rm BW}(s)$ satisfies the initial condition $\eta^j_{\rm BW}(T_0-t_{\rm out})=\eta^{j}(t_{\rm out})$ computed from step (3) as well as conditions (b) through (f) from step (2) of \textbf{Algorithm for $\z$} in \S\ref{sec:algorithm-iprc}.
		\item The lTRC in forward time is then given by $\eta^j(t) = \eta^j_{\rm BW}(T_0-t)$ where $t\in [t_{\rm in}, t_{\rm out}]$.
		\item Compute $\gamma_\varepsilon$, the limit cycle under some small static perturbation $\varepsilon\ll 1$, and find $\mbx_{\varepsilon}^{\rm in}$, the coordinate of the intersection point where $\gamma_\varepsilon(t)$ crosses $\Sigma^{\rm in}$.
		The linear shift in time in region $j$ in response to the static perturbation can be calculated by evaluating the integral (see~\eqref{eq:local-time-shift})
		\[
		T^{j}_{1} = \eta^j(\mbx^{\rm in})\cdot \frac{ \mbx_{\varepsilon}^{\rm in}-\mbx^{\rm in}}{\varepsilon}+\int_{t^{\rm in}}^{t^{\rm out}}\eta^j(\gamma(t))\cdot \frac{\partial F_\varepsilon(\gamma(t))}{\partial \varepsilon}\Big|_{\varepsilon=0}dt.
		\]
	\end{enumerate}
	
	\begin{remark}
		All the local linear shifts in time sum up to the global linear shift in period, that is, 
		$T_1 = \sum_{j=1}^{j=m}T_1^j$.
	\end{remark}

	\subsection{Algorithm for Solving the Homogeneous Variational Equation for LCSCs}\label{sec:algorithm-u}
	
	Here we describe the algorithm for solving the homogeneous variational equation for linear displacement $\mbu$, the shape response to an instantaneous perturbation.
	This makes use of Theorem \ref{thm:main}, which describes different jumping behaviors of $\mbu$ at liftoff, landing, and boundary crossing points.
	Unlike the iPRC and lTRC which require integration backwards in time, the variational dynamics can be solved with forward integration.
	This makes the algorithm comparatively simpler by allowing $\gamma(t)$ and $\mbu(t)$ to be solved simultaneously.
	
	\paragraph{Algorithm for $\mbu$:}
	\begin{enumerate}
		\item Fix an initial condition $\mbx_0=\gamma(0)$ on the limit cycle and an initial condition $\mbu_0=\mbu(0)$ for the displacement at $\gamma(0)$ of the limit cycle.
		\item Integrate the original differential equation~\eqref{eq:multiple-zone-FP} and the homogeneous variational equation~\eqref{eq:var}  simultaneously forward in time  and numerically solve for $\mbu(t)$ over one period $0\le t \le T_0$, where $\mbu$ satisfies
		\begin{enumerate}
			\item $\mbu(0)=\mbu_0$.
			\item For $t$ such that $\gamma(t)$ lies in the interior of the domain, 
			$$ \frac{d\mbu}{dt}=DF^{\rm interior} (\gamma(t)) \mbu $$
			
			\item For $t$ such that $\gamma(t)$ lies within a sliding component along boundary $\Sigma^i$, 
			$$ \frac{d\mbu}{dt}= DF^{\mathrm{slide}_i}(\gamma(t))\mbu $$
			where $DF^{\mathrm{slide}_i}$ is the Jacobian of the sliding vector field $F^{\mathrm{slide}_i}$ given in \eqref{eq:multiple-zone-sliding}.
			\item At any time $t_p$ when $\gamma$ transversely crosses a switching surface with a normal vector $n_p$ separating vector field $F_p^-$ on the incoming side from vector field $F_p^+$ on the outgoing side,
			$$ \mbu^+ = S \mbu^-$$
			where $\mbu^-=\lim_{t\to t_p^-}\mbu(t)$ and $\mbu^+=\lim_{t\to t_p^+}\mbu(t)$ are the displacements just before and just after crossing the surface.
			By the definition for the saltation matrix at transversal crossing point \eqref{eq:salt}, we have $$S=I+\frac{(F_p^+-F_p^-)n_p^\intercal }{n_p^\intercal F_p^-}.$$
			\item At a landing point on the $i$-th hard boundary $\Sigma^i$, update $\mbu$ as
			$$ \mbu^+ = S^i \mbu^-$$
			where $S^i=I-n^i{n^i}^\intercal$ (recall $n^i$ is the normal vector to $\Sigma^i$) and switch integration from the full Jacobian $DF^{\rm interior}$ to the restricted Jacobian $DF^{\mathrm{slide}_i}$.
			\item At a liftoff point on the $i$-th hard boundary $\Sigma^i$, switch integration from the restricted Jacobian $DF^{\mathrm{slide}_i}$ to the full Jacobian $DF^{\rm interior}$; no other change in $\mbu$ is needed.
		\end{enumerate}
	\end{enumerate}
	
	\begin{remark}\label{rem:fund-matrix}
		The \textit{fundamental solution matrix} satisfies 
		\begin{equation*}
		\frac{d\Phi(t,0)}{dt}=DF\Phi(t,0),\, \text{with}\quad \Phi(0,0)=I
		\end{equation*} 
		and takes the initial perturbation $\mbu(0)$ to the perturbation $\mbu(t)$ at time $t$, that is, 
		\[\mbu(t)=\Phi(t,0)\mbu(0).\]
		Computing $\Phi$ therefore requires applying \textbf{Algorithm for $\mbu$} $n$ times, once for each  dimension of the state space.
		Specifically, let $\Phi(t,0)=[\phi_1(t,0)\, ...,\,\phi_n(t,0)]$.
		The $i$-th column $\phi_i(t,0)$ is the solution of the variational equation \eqref{eq:var} with the initial condition $\phi_i(0,0)=e_i$, a unit column vector with zeros everywhere except at the $i$-th row where the entry equals 1.
	\end{remark}
	\begin{remark}\label{rem:monodromy}
		Once $\Phi$ is obtained, we can obtain the monodromy matrix,
		$M= \Phi(T_0,0)$.
		It follows from the periodicity of $\gamma(t)$ that $M$ has $+1$ as an eigenvalue with eigenvector $v$ tangent to the limit cycle at $\mbx_0$; this condition provides a partial consistency check for the algorithm.
	\end{remark}

	\subsection{Algorithms for computing iSRC, the response to sustained perturbation}\label{sec:algorithm-x1}
	
	Now we discuss the calculation of iSRC $\gamma_1$, the linear shape response to a sustained perturbation.
	While $\gamma_1$ shares the same saltation as $\mbu$ at each liftoff, landing and boundary crossing point, $\gamma_1$ satisfies the nonhomogeneous version of the variational equation, \eqref{eq:src} or \eqref{eq:src-multiplescales}, where one of the nonhomogeneous terms depends on the time scaling factor, $\nu_1$ or $\nu^j_1$.
	Moreover, the initial condition for $\gamma_1$ depends on the given perturbation and hence needs to be computed in the algorithm whereas the initial value for $\mbu$ is arbitrarily preassigned.
	
	In the following, we first describe the algorithm for computing $\gamma_1$ using the global uniform rescaling and then consider using piecewise uniform rescaling.
	
	\paragraph{Algorithm for $\gamma_1$ with uniform rescaling}
	\begin{enumerate}
		\item Fix an initial condition $\mbx_0=\gamma(0)$ on the limit cycle.
		\item Compute the linear shift in period $T_1$ using \textbf{Algorithm for $\z$}, then evaluate $\nu_1=T_1/T_0$.
		\item Choose an arbitrary Poincar\'{e} section $\Sigma$ (this can be one of the switching boundaries for appropriate $\mbx_0$) that is transverse to $\gamma$ at $\mbx_0$.
		Compute  $\gamma_\varepsilon$, the limit cycle under some fixed small static perturbation, and find ${\mbx_0}_\varepsilon$, the coordinate of the intersection point where $\gamma_\varepsilon(t)$ crosses $\Sigma$.
		The initial value for $\gamma_1$ at the initial point $\mbx_0$ is then given by 
		$$\gamma_1(0) = \frac{{\mbx_{0}}_\varepsilon - \mbx_0}{\varepsilon}$$
		\item Integrate the original differential equation~\eqref{eq:multiple-zone-FP} with the initial condition $\mbx_0$ and the nonhomogeneous variational equation~\eqref{eq:src} simultaneously forward in time  and numerically solve for $\gamma_1$ over one period $0\le t \le T_0$, where $\gamma_1$ satisfies
		\begin{enumerate}
			\item[(i)] $\gamma_1(0)=({\mbx_{0}}_\varepsilon - \mbx_0)/\varepsilon$.
			\item[(ii)] For $t$ such that $\gamma(t)$ lies in the interior of the domain, 
			$$ \frac{d\gamma_1}{dt}=DF^{\rm interior}(\gamma(t)) \gamma_1 + \nu_1 F^{\rm interior}(\gamma( t)) +\frac{\partial F^{\rm interior}_\varepsilon(\gamma(t))}{\partial \varepsilon}\Big|_{\varepsilon=0} $$
			
			\item[(iii)] For $t$ such that $\gamma(t)$ lies within a sliding component along boundary $\Sigma^i$, 
			$$ \frac{d\gamma_1}{dt}= DF^{\mathrm{slide}_i}(\gamma(t)) \gamma_1 + \nu_1 F^{\mathrm{slide}_i}(\gamma( t)) + \frac{\partial F_\varepsilon^{\mathrm{slide}_i}(\gamma( t))}{\partial \varepsilon}\Big|_{\varepsilon=0} $$
			where $DF^{\mathrm{slide}_i}$ is the Jacobian of the sliding vector field $F^{\mathrm{slide}_i}$ given in \eqref{eq:multiple-zone-sliding}.
			\item[(iv)] For transversal crossings, landing points, and liftoff points, apply (d), (e) and (f), respectively, from step 2) in \textbf{Algorithm for $\mbu$} in \S\ref{sec:algorithm-u}, by replacing $\mbu$ with $\gamma_1$.
		\end{enumerate}
	\end{enumerate}
	
	Next we consider the case when $\gamma(t)$ exhibits $m$ different uniform timing sensitivities at regions $\RR^1, ..., \RR^m$, each bounded by two local timing surfaces, as discussed in \S\ref{sec:algorithm-lTRC}.
	Piecewise uniform rescaling is therefore needed to compute the shape response curve.
	The procedure for obtaining $\gamma_1$ in this case is nearly the same as described in \textbf{Algorithm for $\gamma_1$ with uniform rescaling}, except we now need to compute various rescaling factors using the lTRC.
	This hence leads to different variational equations that need to be solved.
	On the other hand, the local timing surfaces naturally serve as the Poincar\'{e} sections that are required to compute the initial values for $\gamma_1$ in the uniform rescaling case.
	
	\paragraph{Algorithm for $\gamma_1$ with piecewise uniform rescaling}
	\begin{enumerate}
		\item Take the initial condition for $\gamma(t)$ to be  $\gamma(0)=\mbx_0\in \Sigma$, where $\Sigma$ is one of the local timing surfaces.
		Compute $\gamma(t)$, the unperturbed trajectory, and $\gamma_\varepsilon(t)$, the trajectory under some static perturbation $0<\varepsilon\ll 1$, by integrating~\eqref{eq:multiple-zone-FP}.
		\item For $j\in\{1,...,m\}$, compute $T^j_0$, the time that $\gamma(t)$ spends in region $j$ and $T^j_1$, the linear shift in time in region $j$ using \textbf{Algorithm for $\eta^j$}, and then evaluate $\nu^j_1=T^j_1/T^j_0$.
		\item Compute ${\mbx_0}_\varepsilon$, the coordinate of the intersection point where $\gamma_\varepsilon(t)$ crosses $\Sigma$.
		The initial value for $\gamma_1$ at the initial point $\mbx_0$ is given by 
		$$\gamma_1(0) = \frac{{\mbx_0}_\varepsilon- \mbx_0}{\varepsilon}$$ 
		\item Integrate the original differential equation~\eqref{eq:multiple-zone-FP} with the initial condition $\mbx_0$ and the piecewise nonhomogeneous variational equation~\eqref{eq:src-multiplescales}  simultaneously forward in time  and numerically solve for $\gamma_1$ over one period $0\le t \le T_0$, where $\gamma_1$ satisfies
		\begin{enumerate}
			\item[(i)] $\gamma_1(0)=({\mbx_{0}}_\varepsilon - \mbx_0)/\varepsilon$.
			\item[(ii)] For $t$ such that $\gamma(t)$ lies in the intersection of the interior of the domain and region $\RR^j$, 
			$$ \frac{d\gamma_1}{dt}=DF^{\mathrm{interior}_j}(\gamma(t)) \gamma_1 + \nu^j_1 F^{\mathrm{interior}_j}(\gamma( t)) +\frac{\partial F_\varepsilon^{\mathrm{interior}_j}(\gamma( t))}{\partial \varepsilon}\Big|_{\varepsilon=0} $$
			where $DF^{\mathrm{interior}_j}$ is the Jacobian of the interior vector field $F^{\mathrm{interior}_j}$ in $\RR^j$.
			\item[(iii)] For $t$ such that $\gamma(t)$ lies within the intersection of a hard boundary $\Sigma^i$ and region $\RR^j$, 
			$$ \frac{d\gamma_1}{dt}= DF^{\mathrm{slide}_i}(\gamma(t)) \gamma_1 + \nu^j_1 F^{\mathrm{slide}_i}(\gamma( t)) + \frac{\partial F_\varepsilon^{\mathrm{slide}_i}(\gamma( t))}{\partial \varepsilon}\Big|_{\varepsilon=0} $$
			where $DF^{\mathrm{slide}_i}$ is the Jacobian of the sliding vector field $F^{\mathrm{slide}_i}(\mbx)=F^{\mathrm{interior}_j}(\mbx)-(n^i\cdot F^{\mathrm{interior}_j}(\mbx))n^i$ given in \eqref{eq:multiple-zone-sliding}.
			\item[(iv)] For transversal crossings, landing points, and liftoff points, apply (d), (e) and (f), respectively, from step 2) in \textbf{Algorithm for $\mbu$} in \S\ref{sec:algorithm-u}, replacing $\mbu$ with $\gamma_1$.
		\end{enumerate}
	\end{enumerate}


\newpage
	\begin{thebibliography}{99}
		
		\bibitem[Aihara and Suzuki (2010)]{AiharaSuzuki2010}
		K.~Aihara and H.~Suzuki. Theory of hybrid dynamical systems and its applications to biological and medical systems.
		\textit{Philosophical Transactions of the Royal Society A}  \textbf{368} (2010), 4893-4914.
		
		\bibitem[Barajon et al.(1992)]{barajon1992} I.~Barajon, J.~Gossard and H.~Hultborn. Induction of fos expression by activity in the spinal rhythm generator for scratching, \textit{Brain research} \textbf{588}(1) (1992), 168--172.
		
		\bibitem[Bernardo et al.(2008)]{bernardo2008} M.~Bernardo, C.~Budd, A.R.~Champneys and P.~Kowalczyk. Piecewise-smooth dynamical systems: theory and applications, \textit{Springer Science and Business Media} \textbf{163} (2008).
		
		\bibitem[Branicky(1998)]{branicky1998} M.~Branicky. Multiple Lyapunov functions and other analysis tools for switched and hybrid systems, \textit{IEEE Transactions on automatic control} \textbf{43}(4) (1998), 475--482.
		
		
		\bibitem[Brown et al.(2004)]{brown2004} E.~Brown, J.~Moehlis and P.~Holmes. On the phase reduction and response dynamics of neural oscillator populations, \textit{Neural computation} \textbf{16}(4) (2004), 673-715.
		
		\bibitem[Burden et al.(2015)]{BRS2015}
		S. A.~Burden, S.~Revzen and S. S.~Sastry. Model reduction near periodic orbits of hybrid dynamical systems. \textit{IEEE Transactions on Automatic Control}  \textbf{60}(10), (2015) 2626-2639.
		
		\bibitem[Castejon et al.(2013)]{castejon2013} O.~Castejon, A.~Guillamon and G.~Huguet. Phase-amplitude response functions for transient-state stimuli, \textit{Journal of Mathematical Neuroscience} \textbf{3}(1) (2013), 13.
		
		\bibitem[Chartrand et al.(2018)]{CGL18}
		T.~Chartrand, M.S.~Goldman and T.J.~Lewis. Synchronization of electrically coupled resonate-and-fire neurons. \textit{arXiv preprint arXiv:1801.05874} (2018).
		
		\bibitem[Chiel(2007)]{chiel2007} H.J.~Chiel. Aplysia feeding biomechanics, \textit{Scholarpedia} \textbf{2}(9) (2007), 4165.
		
		\bibitem[Coombes et al.(2012)]{coombes2012} S.~Coombes, R.~Thul and K.~Wedgwood. Nonsmooth dynamics in spiking neuron models, \textit{Physica D: Nonlinear Phenomena} \textbf{241}(22) (2012), 2042-2057.
		
		\bibitem[Dieci and Lopez(2011)]{DL11}
		L.~Dieci and L.~Lopez. Fundamental matrix solutions of piecewise smooth differential systems, \textit{Mathematics and Computers in Simulation} \textbf{81}(5) (2011), 932-953.
		
		\bibitem[Diekman et al.(2017)]{diekman2017} C.~Diekman, P.~Thomas and C.~Wilson. Eupnea, tachypnea, and autoresuscitation in a closed-loop respiratory control model, \textit{Journal of Neurophysiology} \textbf{118}(4) (2017), 2194--2215.
		
		
		\bibitem[Doedel(1981)]{doedel1981} E.J.~Doedel. Auto: a program for the automatic bifurcation analysis of autonomous systems, \textit{Congressus Numerantum} \textbf{30} (1981), 265--284.
		
		\bibitem[Doedel et al.(2009)]{doedel2009} E.J.~Doedel, A.R.~Champneys, T.F.~Fairgrieve, Y.A.~Kuznetsov, K.E.~Oldeman, R.C.~Paffenroth, B.~Sanstede, X.J.~Wang, and C.~Zhang,  Auto-07p: continuation and bifurcation software for ordinary differential equations. available from: http://cmvl.cs.concordia.ca/, (2009).
		
		\bibitem[Ermentrout(1996)]{ermentrout1996} B.~Ermentrout. Type I membranes, phase resetting curves, and synchrony, \textit{Neural Computation} \textbf{8}(5) (1996), 979--1001.
		
		\bibitem[Ermentrout(2002)]{ermentrout2002} B.~Ermentrout. Simulating, analyzing, and animating dynamical systems: a guide to XPPAUT for researchers and students, \textit{Siam} \textbf{14} (2002).
		
		\bibitem[Ermentrout and Kopell(1986)]{EK1986} B.~Ermentrout and N.~Kopell. Parabolic bursting in an excitable system coupled with a slow oscillation, \textit{SIAM Journal on Applied Mathematics} \textbf{46}(2) (1986), 233--253.
		
		\bibitem[Ermentrout and Terman(2010)]{ET2010} B.~Ermentrout and D.~Terman. Mathematical foundations of neuroscience, \textit{Springer Science \& Business Media} \textbf{35} (2010).
		
		\bibitem[Filippov(1988)]{filippov1988} A.F.~Filippov. Differential equations with discontinuous right-hand sides. \textit{Mathematics and its Applications, Kluwer Academic, Dordrecht}, (1988).
		
		\bibitem[Galvanetto and Bishop(1999)]{GB99} U.~Galvanetto and S.~Bishop. Dynamics of a simple damped oscillator undergoing stick-slip vibrations, \textit{Meccanica} \textbf{34} (1999), 337--347.
		
		\bibitem[Galvanetto(2001)]{Galvanetto2001} U.~Galvanetto. Some discontinuous bifurcations in a two-block stick-slip system, \textit{Journal of Sound and Vibration} \textbf{248}(4) (2001), 653--669.
		
		\bibitem[Gelfand et al.(2004)]{gelfand1988} I.~Gelfand, G.~Orlovsky and M.~Shik. Locomotion and scratching in tetrapods, \textit{Neural control of rhythmic movements in vertebrates} (1988), 167--199.
		
		\bibitem[Goebel et al.(2009)]{GST2009} R.~Goebel, R.~Sanfelice and A.~Teel. Hybrid dynamical systems, \textit{IEEE control systems magazine} \textbf{29}(2) (2009), 28--93.
		
		\bibitem[Govaerts and Sautois(2006)]{GS2006} W.~Govaerts and B.~Sautois. Computation of the phase response curve: a direct numerical approach, \textit{Neural Computation} \textbf{18}(4) (2006), 817--847.
		
		\bibitem[Guckenheimer and Javeed(2018)]{guckenheimer2018} J.~Guckenheimer and A.~Javeed. Locomotion: exploiting noise for state estimation, \textit{Biological cybernetics} (2018), 1--12.
		
		\bibitem[Guillamon and Huguet(2009)]{GH2009} A.~Guillamon and G.~Huguet. A computational and geometric approach to phase resetting curves and surfaces, \textit{SIAM Journal on Applied Dynamical Systems} \textbf{8}(3) (2009), 1005--1042.
		
		\bibitem[Haghverdi et al.(2003)]{HTP2003} E.~Haghverdi, P.~Tabuada and G.~Pappas. Bisimulation relations for dynamical and control systems, \textit{Electronic Notes in Theoretical Computer Science} \textbf{69} (2003), 120--136.
		
		\bibitem[Haghverdi et al.(2005)]{HTP2005} E.~Haghverdi, P.~Tabuada and G.~Pappas. Bisimulation relations for dynamical, control systems and hybrid systems, \textit{Theoretical Computer Science} \textbf{342}(2-3) (2005), 229-261.
		
		
		\bibitem[Holmes et al.(2006)]{holmes2006} P.~Holmes, R.J.~Full, D.~Koditschek and J.~Guckenheimer. The dynamics of legged locomotion: Models, analyses, and challenges, \textit{SIAM review} \textbf{48}(2) (2006), 207--304.
		
		\bibitem[Izhikevich(2000)]{izhikevich2000} E.~Izhikevich. Phase Equations for Relaxation Oscillators, \textit{SIAM Journal on Applied Mathematics} \textbf{60}(5) (2000), 1789--1804.
		
		\bibitem[Jeffrey(2018)]{Jeffrey2018} M.~Jeffrey. Hidden Dynamics: The Mathematics of Switches, Decisions and Other Discontinuous Behaviour. Springer, New York, 2018
		
		\bibitem[Jelbart and Wechselberger(2020)]{JW2020} S.~Jelbart and M.~Wechselberger. Two-stroke relaxation oscillators, \textit{Nonlinearity} \textbf{33}(5) (2020), 2364
		
		\bibitem[Jordan et al.(2007)]{JSS07} D.~Jordan, P.~Smith, and P.~Smith. Nonlinear ordinary differential equations: an introduction for scientists and engineers. (Vol. 10). \textit{Oxford University Press on Demand} (2007)
		
		\bibitem[Josic et al.(2006)]{JosicShea-BrownMoehlis2006} J.~Kresimir, E. Shea-Brown and J.Moehlis. Isochron, \textit{Scholarpedia}, \textbf{1}(8) (2006), 1361.
		
		\bibitem[Kuramoto(1975)]{kuramoto1975}
		Kuramoto, Yoshiki. "Self-entrainment of a population of coupled non-linear oscillators." \textit{International symposium on mathematical problems in theoretical physics.} Springer, Berlin, Heidelberg, 1975.
		
		\bibitem[Kuramoto(1984)]{kuramoto1984chemical} Y.~Kuramoto. Chemical oscillations, Waves and Turbulence. \textit{Springer-Verlag}, Berlin, (1984).
		
		\bibitem[Lee et al.(2009)]{Lee2009} D.~Lee, H.~Kim and S.~Sastry. Feedback linearization vs. adaptive sliding mode control for a quadrotor helicopter. \textit{International Journal of control, Automation and systems} \textbf{7}(3) (2009), 419--428.
		
		\bibitem[Leine and Nijmeijer(2013)]{LN2013} R.I.~Leine and H.~Nijmeijer. Dynamics and bifurcations of non-smooth mechanical systems, \textit{Springer Science and Business Media} \textbf{18} (2013).
		
		\bibitem[Lyttle et al.(2017)]{lyttle2017} D.~Lyttle, J.~Gill, K.~Shaw, P.~Thomas and H.~Chiel. Robustness, flexibility, and sensitivity in a multifunctional motor control model, \textit{Biological cybernetics} \textbf{111}(1) (2017), 25--47.
		
		\bibitem[Meiss(2007)]{meiss07} J.~Meiss. Differential dynamical systems. \textit{Siam} \textbf{14} (2007).
		
		\bibitem[Monga and Moehlis(2018)]{monga2018} B.~Monga and J.~Moehlis. Optimal phase control of biological oscillators using augmented phase reduction, \textit{Biological cybernetics} (2018), 1--18.
		
		\bibitem[Monga et al.(2018)]{monga2018b} B.~Monga, D.~Wilson, T.~Matchen and Jeff Moehlis. Phase reduction and phase-based optimal control for biological systems: a tutorial, \textit{Biological cybernetics} (2018), 1--36.
		
		
		\bibitem[Mortin and Stein(1989)]{mortin1989} L.~Mortin and P.~Stein. Spinal cord segments containing key elements of the central pattern generators for three forms of scratch reflex in the turtle, \textit{Journal of Neuroscience} \textbf{9}(7) (1989), 2285--2296.
		
		
		\bibitem[Park (2013)]{park2013} Y.~Park. Infinitesimal Phase Response Curves for Piecewise Smooth Dynamical Systems (M.S. thesis). Case Western Reserve University, Cleveland, OH. (2013).
		
		\bibitem[Park et al.(2017)]{Park2017} Y.~Park, S.~Heitmann and B.~Ermentrout. The Utility of Phase Models in Studying Neural Synchronization. Chapter 36 of \textit{Computational Models of Brain and Behavior}, Wiley Online Library, (2017). \textit{arXiv preprint arXiv:1707.05713}.
		
		\bibitem[Park et al.(2018)]{park2018} Y.~Park, K.M.~Shaw, H.J.~Chiel and P.J.~Thomas. The infinitesimal phase response curves of oscillators in piecewise smooth dynamical systems, \textit{European Journal of Applied Mathematics} \textbf{10} (2018), 1017.
		
		\bibitem[Pérez-Cervera et al.(2020)]{perez-cervera2020} A.~Pérez-Cervera, T.M.~Seara and G.~Huguet. Global phase-amplitude description of oscillatory dynamics via the parameterization method. \textit{arXiv preprint arXiv:2004.03647} (2020).
		
		\bibitem[Revzen and Guckenheimer(2011)]{revzen2011} S.~Revzen and J.~Guckenheimer. Finding the dimension of slow dynamics in a rhythmic system, \textit{Journal of The Royal Society Interface} \textbf{9}(70) (2011), 957--971.
		
		\bibitem[Schwemmer and Lewis(2012)]{SL2012} M.~Schwemmer and T.~Lewis. The theory of weakly coupled oscillators, \textit{Phase response curves in neuroscience. Springer, New York, NY} (2012), 3-31.
		
		\bibitem[Shaw et al.(2015)]{shaw2015} K.~Shaw, D.~Lyttle, J.~Gill, M.~Cullins, J.~McManus, H.~Lu, P.~Thomas and H.~Chiel. The significance of dynamical architecture for adaptive responses to mechanical loads during rhythmic behavior, \textit{Journal of Computational Neuroscience} \textbf{38}(1) (2015), 25--51.
		
		\bibitem[Shaw et al.(2012)]{shaw2012} K.~Shaw, Y.~Park, H.~Chiel and P.~Thomas. Phase resetting in an asymptotically phaseless system: on the phase response of limit cycles verging on a heteroclinic orbit, \textit{SIAM Journal on Applied Dynamical Systems} \textbf{11}(1) (2012), 350--391.
		
		
		\bibitem[Shirasaka et al.(2017)]{shirasaka2017} S.~Shirasaka, W.~Kurebayashi and H.~Nakao. Phase reduction theory for hybrid nonlinear oscillators, \textit{Physical Review E} \textbf{95}(1) (2017), 012212.
		
		\bibitem[Sismondo(1990)]{sismondo1990animal} E.~Sismondo. Synchronous, alternating, and phase-locked stridulation by a tropical katydid, \textit{Science} \textbf{249} (1990), 55--58.
		
		\bibitem[Slotine and Sastry(1983)]{SS83} J.~Slotine and S.~Sastry. Tracking control of non-linear systems using sliding surfaces, with application to robot manipulators. \textit{International journal of control} \textbf{38}(2) (1983), 465--492.
		
		\bibitem[Slotine(1984)]{Slotine84} J.~Slotine. Sliding controller design for non-linear systems. \textit{International Journal of control} \textbf{40}(2) (1984), 421--434.
		
		\bibitem[Somers and Kopell(1993)]{SK93} D.~Somers and N.~Kopell. Rapid synchronization through fast threshold modulation. \textit{Biological Cybernetics} \textbf{68} (1993), 393--407.
		
		\bibitem[Spardy et al.(2011a)]{spardy2011a} L.~Spardy, S.~Markin, N.~Shevtsova, B.~Prilutsky, I.~Rybak and J.~Rubin. A dynamical systems analysis of afferent control in a neuromechanical model of locomotion: I. Rhythm generation, \textit{Journal of neural engineering} \textbf{8}(6) (2011), 065003.
		
		\bibitem[Spardy et al.(2011b)]{spardy2011b} L.~Spardy, S.~Markin, N.~Shevtsova, B.~Prilutsky, I.~Rybak and J.~Rubin. A dynamical systems analysis of afferent control in a neuromechanical model of locomotion: II. Phase asymmetry, \textit{Journal of neural engineering} \textbf{8}(6) (2011), 065004.
		
		\bibitem[Sutton et al.(2004)]{sutton2004} G.P.~Sutton, E.V.~Mangan, D.M.~Neustadter, R.D.~Beer, P.E.~Crago and H.J.~Chiel. Neural control exploits changing mechanical advantage and context dependence to generate different feeding responses in Aplysia, \textit{Biological cybernetics} \textbf{91}(5) (2004), 333--345.
		
		\bibitem[Taylor et al.(2008)]{Taylor2008} S.~Taylor, R.~Gunawan, L.~Petzold and F.~Doyle III. Sensitivity measures for oscillating systems: Application to mammalian circadian gene network. \textit{IEEE transactions on automatic control} \textbf{53} Special Issue, (2008), 177--188.
		
		\bibitem[Wiggins(1994)]{Wiggins1994book}
		S.~Wiggins. Normally hyperbolic invariant manifolds in dynamical systems. \textit{Springer Science $\&$ Business Media}, \textbf{105} (1994).
		
		\bibitem[Wilson(2019)]{wilson2019} D.~Wilson. Isostable reduction of oscillators with piecewise smooth dynamics and complex Floquet multipliers, \textit{Physical Review E} \textbf{99}(2) (2019), 022210.
		
		
		\bibitem[Wilson(2020a)]{wilson2020a} D.~Wilson. A data-driven phase and isostable reduced modeling framework for oscillatory dynamical systems. \textit{Chaos: An Interdisciplinary Journal of Nonlinear Science} \textbf{30}(1) (2020), 013121.
		
		\bibitem[Wilson(2020b)]{wilson2020b} D.~Wilson. Phase-amplitude reduction far beyond the weakly perturbed paradigm. \textit{Physical Review E} \textbf{101}(2) (2020), 022220.
		
		\bibitem[Wilson and Moehlis(2015)]{WM2015} D.~Wilson and J.~Moehlis. Extending phase reduction to excitable media: theory and applications, \textit{SIAM Review} \textbf{57}(2) (2015), 201--222.
		
		\bibitem[Wilson and Ermentrout(2018)]{WE2018} D.~Wilson and B.~Ermentrout. Greater accuracy and broadened applicability of phase reduction using isostable coordinates, \textit{Journal of Mathematical Biology} \textbf{76}(1-2) (2018), 37--66.
		
		\bibitem[Wilson and Moehlis(2016)]{WM2016} D.~Wilson and J.~Moehlis. Isostable reduction of periodic orbits, \textit{Physical Review E} \textbf{94}(5) (2016), 052213.
		
		\bibitem[Winfree(1980)]{winfree1980circadian} A.~Winfree. The Geometry of Biological Time. \textit{Springer-Verlag}, New York, (1980).
		
		\bibitem[Yu et al.(1999)]{yu1999} S.N.~Yu, P.E.~Crago and H.J.~Chiel. Biomechanical properties and a kinetic simulation model of the smooth muscle I2 in the buccal mass of Aplysia, \textit{Biological cybernetics} \textbf{81}(5-6) (1999), 505--513.
		
	\end{thebibliography}
\end{document}